\documentclass{gtart}
\def\ifplaintex{\expandafter\ifx\csname documentclass\endcsname\relax}

\def\gtp{{\mathsurround=0pt\it $\cal G\mskip-2mu$eometry \&\ 
$\cal T\!\!$opology $\cal P\!$ublications}}  % GT publications

\def\recd{{\small Received:\qua\receiveddate\ifx\reviseddate\relax
\else\qquad Revised:\qua\reviseddate\fi\par}} 

\def\lognumber#1{\def\thelognumber{#1}}
\def\volumenumber#1{\def\thevolumenumber{#1}}
\def\volumeyear#1{\def\thevolumeyear{#1}}
\def\papernumber#1{\def\thepapernumber{#1}}
\def\pagenumbers#1#2{\def\startpage{#1}\def\finishpage{#2}}
\def\published#1{\def\publishdate{#1}}

\def\received#1{\def\receiveddate{#1}}
\def\revised#1{\def\reviseddate{#1}}
\def\accepted#1{\def\accepteddate{#1}}
\def\asciititle#1{\def\theasciititle{#1}}

\def\asciiauthors#1{\def\theasciiauthors{#1}}
\def\asciiaddress#1{\def\theasciiaddress{#1}}

\def\coverauthors#1{\def\thecoverauthors{#1}}
\long\def\asciiabstract#1{\long\def\theasciiabstract{#1}}

\let\\\par\let\thelognumber\relax\let\thevolumenumber\relax
\let\thepapernumber\relax\let\thevolumeyear\relax\let\startpage\relax
\let\finishpage\relax\let\publishdate\relax\let\receiveddate\relax
\let\reviseddate\relax\let\accepteddate\relax\let\theasciititle\relax
\let\theasciiauthors\relax\let\theasciiaddress\relax
\let\theasciiabstract\relax

\let\thecoverauthors\relax\let\theasciiemail\relax

\ifplaintex
\font\logobig=cmssbx10 scaled 3836
\font\logomed=cmssbx10 scaled 2557
\else
\font\logobig=cmssbx10 scaled 4200
\font\logomed=cmssbx10 scaled 2800
\fi

\long\def\makeagttitle{   %%% start of definition of \makeagttitle
\count0=\startpage
\agt\hfill      %   Journal title (top left) 
\hbox to 45truept{\vbox to 0pt{\vglue -13truept{\logomed A\kern -.37em{\logobig 
T}\kern -.38em G}\vss}\hss}
\break
{\small Volume \thevolumenumber\ (\thevolumeyear)
\startpage--\finishpage\nl
Published: \publishdate}

\vglue .25truein

{\parskip=0pt\leftskip 0pt plus
1fil\def\\{\par\smallskip}{\Large\bf\thetitle}\par\medskip} \vglue
0.05truein

{\parskip=0pt\leftskip 0pt plus 1fil\def\\{\par}{\sc\theauthors}
\par\medskip}%
 
\vglue 0.03truein 

{\small\leftskip 25truept\rightskip 25truept{\bf Abstract}\stdspace\theabstract

{\bf AMS Classification}\stdspace\theprimaryclass
\ifx\thesecondaryclass\relax\else; \thesecondaryclass\fi\par
{\bf Keywords}\stdspace \thekeywords\par}\vglue 7truept

}   %%%% end of definition of \makeagttitle

\ifplaintex
\font\phead=cmsl9 scaled 950
\font\pnum=cmbx10 scaled 913
\font\pfoot=cmsl9 scaled 950
\headline{\vbox to 0pt{\vskip -4.5mm\line{\small\phead\ifnum
\count0=\startpage ISSN 1472-2739 (on-line) 1472-2747 (printed)
\hfill {\pnum\folio}\else\ifodd\count0\def\\{ }% 
\ifx\theshorttitle\relax\thetitle\else\theshorttitle\fi\hfill{\pnum\folio}
\else\def\\{ and }{\pnum\folio}\hfill\ifx\theshortauthors\relax\theauthors
\else\theshortauthors\fi\fi\fi}\vss}}
\footline{\vbox to 0pt{\vglue 0mm\line{\small\pfoot\ifnum\count0=\startpage
\copyright\ \gtp\hfill\else
\agt, Volume \thevolumenumber\ (\thevolumeyear)\hfill\fi}\vss}}
\else
\font=cmsl9 scaled 1050
\font\lnum=cmbx10 
\font=cmsl9 scaled 1050
\makeatletter
\def\@oddhead{{\small\lhead\ifnum\count0=\startpage ISSN 1472-2739 
(on-line) 1472-2747 (printed)\hfill {\lnum\number\count0}\else\ifodd\count0
\def\\{ }\ifx\theshorttitle\relax \thetitle \else\theshorttitle\fi\hfill
{\lnum\number\count0}\else\def\\{ and }{\lnum\number\count0}
\hfill\ifx\theshortauthors\relax 
\theauthors\else\theshortauthors\fi\fi\fi}}\def\@evenhead{\@oddhead}
\def\@oddfoot{\small\lfoot\ifnum\count0=\startpage\copyright\ \gtp\hfill\else
\agt, Volume \thevolumenumber\ (\thevolumeyear)\hfill\fi}
\def\@evenfoot{\@oddfoot}
\makeatother
\fi
\let\maketitlepage\makeagttitle
\let\makeshorttitle\maketitlepage
\let\maketitle\maketitlepage

\newwrite\gtoutfile
\long\gdef\makeheadfile{  %%% start of definition of \makeheadfile
{\def\\{, }\def\s{ }
\immediate\openout\gtoutfile head.xxx
\immediate\write\gtoutfile{To: math@arxiv.org}
\immediate\write\gtoutfile{Subject: put OR rep NNNNN:ppppp}
\immediate\write\gtoutfile{--text follows this line--}
\immediate\write\gtoutfile{Proxy-for: \ifx\theasciiauthors\relax
\theauthors\else\theasciiauthors\fi\s<\ifx\theasciiemail\relax\theemail\else\theasciiemail\fi>}
\immediate\write\gtoutfile{\noexpand\\}
\immediate\write\gtoutfile{Authors: \ifx\theasciiauthors\relax
\theauthors\else\theasciiauthors\fi}
{\def\\{ }\immediate\write\gtoutfile{Title: \ifx\theasciititle\relax
\thetitle\else\theasciititle\fi}}
\immediate\write\gtoutfile{Subj-class: GT or SG, GR etc}
\immediate\write\gtoutfile{MSC-class: \theprimaryclass\ifx\thesecondaryclass\relax\else, \thesecondaryclass\fi}
\immediate\write\gtoutfile{Journal-ref: Algebr. Geom. Topol. \thevolumenumber\s
(\thevolumeyear) \startpage-\finishpage}
\immediate\write\gtoutfile{Comments: Published by Algebraic and
Geometric Topology at}
\immediate\write\gtoutfile{\s\s\s  http://www.maths.warwick.ac.uk/agt/AGTVol\thevolumenumber/agt-\thevolumenumber-\thepapernumber.abs.html}
\immediate\write\gtoutfile{\noexpand\\}
\immediate\write\gtoutfile{}
\ifx\theasciiabstract\relax
\immediate\write\gtoutfile{\theabstract}\else
\immediate\write\gtoutfile{\theasciiabstract}\fi
\immediate\write\gtoutfile{}
\immediate\write\gtoutfile{\noexpand\\}
\immediate\write\gtoutfile{}
\immediate\closeout\gtoutfile}}  %%% end of definition of \makeheadfile

\def\maketitlepage{\makeagttitle\makeheadfile}
\let\makeshorttitle\maketitlepage
\let\maketitle\maketitlepage

\def\ifplaintex{\expandafter\ifx\csname documentclass\endcsname\relax}

\def\gtp{{\mathsurround=0pt\it $\cal G\mskip-2mu$eometry \&\ 
$\cal T\!\!$opology $\cal P\!$ublications}}  % GT publications

\def\recd{{\small Received:\qua\receiveddate\ifx\reviseddate\relax
\else\qquad Revised:\qua\reviseddate\fi\par}} 

\def\lognumber#1{\def\thelognumber{#1}}
\def\volumenumber#1{\def\thevolumenumber{#1}}
\def\volumeyear#1{\def\thevolumeyear{#1}}
\def\papernumber#1{\def\thepapernumber{#1}}
\def\pagenumbers#1#2{\def\startpage{#1}\def\finishpage{#2}}
\def\published#1{\def\publishdate{#1}}

\def\received#1{\def\receiveddate{#1}}
\def\revised#1{\def\reviseddate{#1}}
\def\accepted#1{\def\accepteddate{#1}}
\def\asciititle#1{\def\theasciititle{#1}}

\def\asciiauthors#1{\def\theasciiauthors{#1}}
\def\asciiaddress#1{\def\theasciiaddress{#1}}

\def\coverauthors#1{\def\thecoverauthors{#1}}
\long\def\asciiabstract#1{\long\def\theasciiabstract{#1}}

\let\\\par\let\thelognumber\relax\let\thevolumenumber\relax
\let\thepapernumber\relax\let\thevolumeyear\relax\let\startpage\relax
\let\finishpage\relax\let\publishdate\relax\let\receiveddate\relax
\let\reviseddate\relax\let\accepteddate\relax\let\theasciititle\relax
\let\theasciiauthors\relax\let\theasciiaddress\relax
\let\theasciiabstract\relax

\let\thecoverauthors\relax\let\theasciiemail\relax

\ifplaintex
\font\logobig=cmssbx10 scaled 3836
\font\logomed=cmssbx10 scaled 2557
\else
\font\logobig=cmssbx10 scaled 4200
\font\logomed=cmssbx10 scaled 2800
\fi

\long\def\makeagttitle{   %%% start of definition of \makeagttitle
\count0=\startpage
\agt\hfill      %   Journal title (top left) 
\hbox to 45truept{\vbox to 0pt{\vglue -13truept{\logomed A\kern -.37em{\logobig 
T}\kern -.38em G}\vss}\hss}
\break
{\small Volume \thevolumenumber\ (\thevolumeyear)
\startpage--\finishpage\nl
Published: \publishdate}

\vglue .25truein

{\parskip=0pt\leftskip 0pt plus
1fil\def\\{\par\smallskip}{\Large\bf\thetitle}\par\medskip} \vglue
0.05truein

{\parskip=0pt\leftskip 0pt plus 1fil\def\\{\par}{\sc\theauthors}
\par\medskip}%
 
\vglue 0.03truein 

{\small\leftskip 25truept\rightskip 25truept{\bf Abstract}\stdspace\theabstract

{\bf AMS Classification}\stdspace\theprimaryclass
\ifx\thesecondaryclass\relax\else; \thesecondaryclass\fi\par
{\bf Keywords}\stdspace \thekeywords\par}\vglue 7truept

}   %%%% end of definition of \makeagttitle

\ifplaintex
\font\phead=cmsl9 scaled 950
\font\pnum=cmbx10 scaled 913
\font\pfoot=cmsl9 scaled 950
\headline{\vbox to 0pt{\vskip -4.5mm\line{\small\phead\ifnum
\count0=\startpage ISSN 1472-2739 (on-line) 1472-2747 (printed)
\hfill {\pnum\folio}\else\ifodd\count0\def\\{ }% 
\ifx\theshorttitle\relax\thetitle\else\theshorttitle\fi\hfill{\pnum\folio}
\else\def\\{ and }{\pnum\folio}\hfill\ifx\theshortauthors\relax\theauthors
\else\theshortauthors\fi\fi\fi}\vss}}
\footline{\vbox to 0pt{\vglue 0mm\line{\small\pfoot\ifnum\count0=\startpage
\copyright\ \gtp\hfill\else
\agt, Volume \thevolumenumber\ (\thevolumeyear)\hfill\fi}\vss}}
\else
\font=cmsl9 scaled 1050
\font\lnum=cmbx10 
\font=cmsl9 scaled 1050
\makeatletter
\def\@oddhead{{\small\lhead\ifnum\count0=\startpage ISSN 1472-2739 
(on-line) 1472-2747 (printed)\hfill {\lnum\number\count0}\else\ifodd\count0
\def\\{ }\ifx\theshorttitle\relax \thetitle \else\theshorttitle\fi\hfill
{\lnum\number\count0}\else\def\\{ and }{\lnum\number\count0}
\hfill\ifx\theshortauthors\relax 
\theauthors\else\theshortauthors\fi\fi\fi}}\def\@evenhead{\@oddhead}
\def\@oddfoot{\small\lfoot\ifnum\count0=\startpage\copyright\ \gtp\hfill\else
\agt, Volume \thevolumenumber\ (\thevolumeyear)\hfill\fi}
\def\@evenfoot{\@oddfoot}
\makeatother
\fi
\let\maketitlepage\makeagttitle
\let\makeshorttitle\maketitlepage
\let\maketitle\maketitlepage

\newwrite\gtoutfile
\long\gdef\makeheadfile{  %%% start of definition of \makeheadfile
{\def\\{, }\def\s{ }
\immediate\openout\gtoutfile head.xxx
\immediate\write\gtoutfile{To: math@arxiv.org}
\immediate\write\gtoutfile{Subject: put OR rep NNNNN:ppppp}
\immediate\write\gtoutfile{--text follows this line--}
\immediate\write\gtoutfile{Proxy-for: \ifx\theasciiauthors\relax
\theauthors\else\theasciiauthors\fi\s<\ifx\theasciiemail\relax\theemail\else\theasciiemail\fi>}
\immediate\write\gtoutfile{\noexpand\\}
\immediate\write\gtoutfile{Authors: \ifx\theasciiauthors\relax
\theauthors\else\theasciiauthors\fi}
{\def\\{ }\immediate\write\gtoutfile{Title: \ifx\theasciititle\relax
\thetitle\else\theasciititle\fi}}
\immediate\write\gtoutfile{Subj-class: GT or SG, GR etc}
\immediate\write\gtoutfile{MSC-class: \theprimaryclass\ifx\thesecondaryclass\relax\else, \thesecondaryclass\fi}
\immediate\write\gtoutfile{Journal-ref: Algebr. Geom. Topol. \thevolumenumber\s
(\thevolumeyear) \startpage-\finishpage}
\immediate\write\gtoutfile{Comments: Published by Algebraic and
Geometric Topology at}
\immediate\write\gtoutfile{\s\s\s  http://www.maths.warwick.ac.uk/agt/AGTVol\thevolumenumber/agt-\thevolumenumber-\thepapernumber.abs.html}
\immediate\write\gtoutfile{\noexpand\\}
\immediate\write\gtoutfile{}
\ifx\theasciiabstract\relax
\immediate\write\gtoutfile{\theabstract}\else
\immediate\write\gtoutfile{\theasciiabstract}\fi
\immediate\write\gtoutfile{}
\immediate\write\gtoutfile{\noexpand\\}
\immediate\write\gtoutfile{}
\immediate\closeout\gtoutfile}}  %%% end of definition of \makeheadfile

\def\maketitlepage{\makeagttitle\makeheadfile}
\let\makeshorttitle\maketitlepage
\let\maketitle\maketitlepage

\lognumber{32}
\volumenumber{1}
\volumeyear{2001}
\papernumber{32}
\pagenumbers{627}{686}
\received{9 April 2001}
\revised{28 August 2001}
\accepted{5 September 2001}
\published{30 October 2001}

\usepackage{amsmath}
\usepackage{amssymb}

\usepackage{texdraw}
\usepackage{xypic}

\newcommand{\mA}{\mathcal A}
\newcommand{\mC}{\mathcal C}
\newcommand{\mD}{\mathcal D}
\newcommand{\mF}{\mathcal F}

\newcommand{\mR}{\mathcal R}
\newcommand{\mT}{\mathcal T}
\newcommand{\mS}{\mathcal S}
\newcommand{\mV}{\mathcal V}
\newcommand{\omV}{\overline{\mathcal V}}

\newcommand{\frg}{\mathfrak g}
\newcommand{\frsl}{\mathfrak sl}

\newcommand{\ep}{\epsilon}
\newcommand{\vep}{\varepsilon}

\newcommand{\Z}{\mathbb Z}
\newcommand{\R}{\mathbb R}
\newcommand{\C}{\mathbb C}
\newcommand{\Q}{\mathbb Q}
\newcommand{\I}{\mathbb I}

\newcommand{\card}{\text{\rm card}}
\newcommand{\col}{\text{\rm col}}
\newcommand{\diag}{\text{\rm diag}}
\newcommand{\Dim}{\text{\rm Dim}}
\newcommand{\e}{\text{\rm e}}
\newcommand{\emb}{\text{\rm emb}}
\newcommand{\End}{\text{\rm{End}}}

\newcommand{\Hom}{\text{\rm{Hom}}}
\newcommand{\id}{\text{\rm{id}}}

\newcommand{\interior}{\text{\rm int}}

\newcommand{\lk}{\text{\rm lk}}
\newcommand{\Mod}{\text{\rm Mod}}
\newcommand{\ns}{\text{\rm n}}
\newcommand{\nul}{\text{\rm null}}
\newcommand{\nvf}{\text{\rm nvf}}

\newcommand{\os}{\text{\rm o}}

\newcommand{\s}{\text{\rm s}}
\newcommand{\sign}{\text{\rm{sign}}}

\newcommand{\Span}{\text{\rm{Span}}}
\newcommand{\tr}{\text{\rm tr}}

\newcommand{\um}{\underline{m}}

\newcommand{\musum}{\sum_{\underline{\mu} \in \{ \pm 1 \}^{n} }}
\newcommand{\msum}{\sum_{\underline{m}=\underline{0}}^{\underline{\alpha}-\underline{1}}}
\newcommand{\rhoalphasum}{\sum_{j=1}^{n} \frac{\rho_{j}}{\alpha_{j}}}

\newcommand{\kvi}{\frac{i \pi}{4}}

\newcommand{\ra}{\rightarrow}
\newcommand{\Lra}{\Longrightarrow}
\newcommand{\sm}{\setminus}

\newtheorem{thm}{Theorem}[section]

\newtheorem{cor}[thm]{Corollary}
\newtheorem{lem}[thm]{Lemma}

\theoremstyle{definition}
\newtheorem{conv}[thm]{Conventions}
\newtheorem{rem}[thm]{Remark}

\newcommand{\refthm}[1]{Theorem~\ref{#1}}
    \newcommand{\refconv}[1]{Conventions~\ref{#1}}
    \newcommand{\refcor}[1]{Corollary~\ref{#1}}
    \newcommand{\reflem}[1]{Lemma~\ref{#1}}
    \newcommand{\refrem}[1]{Remark~\ref{#1}}

\newcommand{\HS}{\noindent \hfill{$\sq$}}

\begin{document}

\title{Reshetikhin--Turaev invariants of Seifert\\$3$--manifolds and a rational surgery formula}
\shorttitle{Reshetikhin--Turaev invariants of Seifert $3$--manifolds}
\asciititle{Reshetikhin-Turaev invariants of Seifert 3-manifolds and a rational surgery formula}

\author{S\o ren Kold Hansen}
\coverauthors{S\noexpand\o ren Kold Hansen}
\asciiauthors{Soren Kold Hansen}
\address{Institut de Recherche Math\'{e}matique Avanc\'{e}e, Universit\'{e} Louis Pasteur - C.N.R.S.\\7 rue Ren\'{e} Descartes, 67084 Strasbourg, France}
\asciiaddress{Institut de Recherche Mathematique Avancee, Universite Louis Pasteur - C.N.R.S.\\7 rue Rene Descartes, 67084 Strasbourg, France}

\email{hansen@math.u-strasbg.fr}

\begin{abstract}
We calculate the RT--invariants of all
oriented Seifert manifolds directly from surgery presentations.
We work in the general framework of an
arbitrary modular category as in \cite{Turaev}, and the invariants
are expressed in terms of the $S$-- and $T$--matrices
of the modular category.
In another direction we derive a rational surgery formula,
which states how the
RT--invariants behave under rational surgery along framed links in arbitrary
closed oriented $3$--manifolds with embedded colored ribbon graphs.
The surgery formula is used to give another derivation of the RT--invariants of
Seifert manifolds with orientable base.
\end{abstract}

\asciiabstract{We calculate the RT-invariants of all oriented Seifert
manifolds directly from surgery presentations.  We work in the general
framework of an arbitrary modular category as in [V. G. Turaev,
Quantum invariants of knots and $3$--manifolds, de Gruyter
Stud. Math. 18, Walter de Gruyter (1994)], and the
invariants are expressed in terms of the S- and T-matrices of
the modular category.  In another direction we derive a rational
surgery formula, which states how the RT-invariants behave under
rational surgery along framed links in arbitrary closed oriented
3-manifolds with embedded colored ribbon graphs.  The surgery
formula is used to give another derivation of the RT-invariants of
Seifert manifolds with orientable base.}




\makeshorttitle

\def\undersmile#1{\lower5.7pt\hbox{$\smallsmile$}\kern-0.6em #1}

\section{Introduction}

A major challenge in the theory of quantum invariants of links and
$3$--manifolds, notably the Jones polynomial of links in $S^{3}$ \cite{Jones},
is to determine relationships between these
invariants and classical invariants.
In 1988 Witten \cite{Witten}
gave a sort of an answer by his interpretation of the Jones
polynomial (and its generalizations) in terms of quantum field theory. Witten
not only gave a description of the Jones polynomial in terms of $3$--dimensional
topology/geometry, but he also initiated the era of quantum invariants of
$3$--manifolds by defining invariants $Z_{k}^{G}(M,L) \in \C$ of an
arbitrary closed oriented $3$--manifold $M$ with an embedded colored link $L$ by
quantizing the Chern--Simons field theory associated to a
simply connected compact simple Lie group $G$,
$k$ being an arbitrary positive integer, called the (quantum) level.
The invariant $Z_{k}^{G}(M,L)$ is given by a Feynman path integral 
over the (infinite dimensional) space of gauge equivalence
classes of connections in a $G$--bundle over $M$.
This integral should be understood in a formal way since,
at the moment of
writing, it seems that no mathematically rigorous definition is known,
cf.\ \cite[Sect.~20.2.A]{JohnsonLapidus}.

By using stationary phase approximation techniques together with path 
integral arguments Witten was able \cite{Witten}
to express the leading asymptotics 
of $Z_{k}^{G}(M)$ as $k \ra \infty$ in terms of such 
topological/geometric invariants as Chern--Simons invariants, Reidemeister
torsions and spectral flows, so here we see a way to extract topological
information from the invariants $Z_{k}^{G}(M)$ (here $L=\emptyset$).
Furthermore, a full asymptotic expansion of $Z_k^G(M)$ as $k 
\rightarrow \infty $ is expected on the basis of a full perturbative 
analysis of the Feynman path integral, see e.g.\ \cite{AxelrodSinger1},
\cite{AxelrodSinger2}.

Reshetikhin and Turaev \cite{ReshetikhinTuraev2}
constructed invariants $\tau_{r}^{\frsl_{2}(\C)}(M,L) \in \C$
by a mathematical approach 
via representations of
a quantum group $U_{q}(\frsl_{2}(\C))$, $q=\exp(2\pi i/r)$, $r$
an integer $\geq 2$. Shortly
afterwards, quantum invariants $\tau_{r}^{\frg}(M,L) \in \C$ associated to
other complex simple Lie algebras $\frg$ were constructed using 
representations of the quantum groups $U_{q}(\frg)$, $q=\exp(2\pi i/r)$ a
`nice' root of unity, see \cite{TuraevWenzl1}. Both in Witten's approach and in the
approach of Reshetikhin and Turaev the invariants are part of a 
so-called topological quantum field theory (TQFT). This implies that the
invariants are defined for compact oriented $3$--dimensional cobordisms 
(perhaps with some extra structure on the boundary),
and satisfy certain cut-and-paste axioms, see \cite{Atiyah}, \cite{Quinn},
\cite{Turaev}. The TQFT of Reshetikhin and Turaev can from
an algebraic point of view be given a more general formulation by using
so-called modular categories \cite{Turaev}. The representation theory of
$U_{q}(\frg)$, $\frg$ an arbitrary complex simple Lie algebra, induces such a
modular category if $q=\exp(2\pi i/r)$ is chosen properly,
see \cite{TuraevWenzl1}, the appendix in \cite{TuraevWenzl2}, \cite{Kirillov},
\cite{BakalovKirillov}, \cite{Sawin}, and \cite{Le}.

It is believed that the TQFT's of Witten and Reshetikhin--Turaev coincide.
In particular it is conjectured, that Witten's leading asymptotics for 
$Z_{k}^{G}(M)$
should be valid for the function $r \mapsto \tau_{r}^{\frg}(M)$ in the 
limit $r \ra \infty$ and furthermore, that this function should have a full 
asymptotic expansion.
In this paper we initiate a verification of
this conjecture for oriented Seifert manifolds by deriving formulas for the
RT--invariants of these manifolds. In a subsequent paper
\cite{Hansen2} we then use these formulas to calculate the large $r$ 
asymptotics
of the RT--invariants and thereby prove the so-called asymptotic expansion
conjecture for such manifolds in the $\frsl_{2}(\C)$--case.
The precise formulation of this conjecture,
which is a combination of  Witten's leading asymptotics
and the existence of a full asymptotic expansion of a certain type, was 
proposed  by Andersen in \cite{Andersen},
where he proved it for mapping 
tori of finite order diffeomorphisms of orientable surfaces of
genus at least two
using the gauge theory definition of the quantum invariants.

In the following a Seifert manifold means an oriented Seifert manifold.
Calculations of quantum invariants of lens spaces and other Seifert manifolds
have been done by several people \cite{Andersen}, \cite{Garoufalidis},
\cite{Jeffrey}, \cite{LawrenceRozansky}, \cite{Neil},
\cite{Rozansky}, \cite{Takata1}, \cite{Takata2}, \cite{Turaev} and probably many
more. The papers \cite{Garoufalidis},
\cite{Jeffrey}, \cite{LawrenceRozansky}, \cite{Neil},
\cite{Takata1} and \cite{Takata2} calculate and study the quantum invariants
of lens spaces and other Seifert manifolds with base equal to $S^{2}$. In
\cite{Takata1}, \cite{Takata2} the so-called $P\frsl_{n}(\C)$--invariants are
calculated. The $P\frsl_{n}(\C)$--invariant of a closed oriented $3$--manifold $M$
associated with an integer $r >n$ coprime to $n$ is a factor of
$\tau_{r}^{\frsl_{n}(\C)}(M)$.
Neil \cite{Neil} calculates the $\frsl_{2}(\C)$--invariants
based on Lickorish skein theoretical approach \cite{Lickorish2}.
In \cite{LawrenceRozansky} the $SU(2)$--invariants of certain Seifert manifolds
with base $S^{2}$ are calculated and studied. The class of Seifert manifolds
considered includes the
Seifert manifolds which are integral homology spheres.
Rozansky \cite{Rozansky} derives
a formula for the $SU(2)$--invariants of all Seifert manifolds with
orientable base.
The papers \cite{Garoufalidis},
\cite{Jeffrey}, \cite{LawrenceRozansky} and
\cite{Rozansky} are based on Witten's approach
to the invariants.
Andersen \cite{Andersen} calculates quantum $G$--invariants of all
mapping tori of finite order diffeomorphisms of orientable surfaces of
genus at least two, where $G$ is an arbitrary simply connected compact simple Lie group.
The mapping tori of finite order diffeomorphisms 
of an orientable surface $\Sigma_{g}$ of genus $g$
are precisely the Seifert manifolds with 
base $\Sigma_{g}$ and Seifert Euler number equal to zero.
Turaev has calculated
the RT--invariants associated to an arbitrary
unimodal modular category of all graph
manifolds, cf.\ \cite[Sect.~X.9]{Turaev}.
These manifolds include the Seifert manifolds with orientable base
(but not the ones with non-orientable base).

In this paper we extend the above results in two directions. Firstly,
we calculate the RT--invariants of all Seifert manifolds. In particular we
calculate the invariants of Seifert manifolds with non-orientable base. This
case has to the authors knowledge not been considered before in the literature.
Secondly, our calculations are done for arbitrary modular categories, 
cf.\ \refthm{invariants}.
We present three different calculations of the RT--invariants of Seifert manifolds
with different levels of generality. In our first approach
we calculate the invariants of all Seifert manifolds directly from surgery
presentations only using the theory of RT--invariants of closed
oriented $3$--manifolds without refering to the underlying TQFT.
In our second approach we calculate the RT--invariants of all Seifert manifolds
with orientable base using a rational surgery formula
for the RT--invariants, \refthm{surgery-formula},
derived in this paper.
In these two approaches we work in the framework of an arbitrary modular
category. In our third approach we use a formula for the RT--invariants of
graph manifolds due to Turaev, see \cite[Theorem X.9.3.1]{Turaev}. This formula
is valid for all modular categories satisfying a special condition called
unimodality. As mentioned above the graph
manifolds include the Seifert manifolds with orientable
base. We show that Turaev's formula
specializes to our formula for the invariants of these Seifert manifolds.

The rational surgery formula, \refthm{surgery-formula}, states how the
RT--invariants behave under rational surgeries along framed links in arbitrary
closed oriented $3$--manifolds with embedded colored ribbon graphs.
This formula
generalizes the defining formula for the RT--invariants of closed oriented
$3$--manifolds with embedded colored ribbon graphs (which is a
surgery formula for surgeries on $S^{3}$ with 
embedded colored ribbon graphs along framed links).
The surgery formula 
has the very same form as the surgery formulas presented
in the Chern--Simons TQFT of Witten, see \cite[Sect.~4]{Witten},
\cite{LawrenceRozansky}, \cite{Rozansky}.

In the final part of the paper we analyse more carefully the
$\frsl_{2}(\C)$--case.
In the general formulas for the RT--invariants of
the Seifert manifolds, see \refthm{invariants}, a certain factor of so-called $S$-- and
$T$--matrices is present. In the $\frsl_{2}(\C)$--case the
$S$-- and $T$--matrices can be identified (up to normalization)
with the values of a certain
representation $\mR$ of $SL(2,\Z)$
in the standard generators of $SL(2,\Z)$. This representation
has been carefully studied by Jeffrey in \cite{Jeffrey},
where an explicit formula for $\mR(A)$ in terms of the
entries of $A \in SL(2,\Z)$ is given. We use this formula
to give expressions for the RT--invariants
of the Seifert manifolds in terms of the Seifert invariants, see \refthm{CS11}.
\refthm{CS11} generalizes results in the literature,
in particular the formulas for the RT--invariants
of Seifert manifolds with orientable base given in \cite{Rozansky}.

The paper is organized as follows.
In Sect.~\ref{sec-Seifert-manifolds} we recall the definition and
classification of Seifert manifolds
\cite{Seifert1}, \cite{Seifert2}. We also present surgery presentations of
the Seifert manifolds due to Montesinos \cite{Montesinos}.
In Sect.~\ref{sec-Modular-categories} we give a short introduction to
the modular categories. This
is a preliminary section intended to fix notation
used throughout in the paper.
In Sect.~\ref{sec-The-Reshetikhin} we calculate the RT--invariants of
all Seifert manifolds directly
from surgery presentations.  
In Sect.~\ref{sec-A-generalized} we derive
the rational surgery formula for the RT--invariants of closed oriented
$3$--manifolds with embedded colored ribbon graphs.
In Sect.~\ref{sec-A-second} we calculate the RT--invariants
of the Seifert manifolds with orientable base using the surgery formula.
In Sect.~\ref{sec-A-third} we show that Turaev's formula for
the RT--invariants of graph manifolds specializes to our formula
for the RT--invariants of Seifert manifolds with orientable base.
In Sect.~\ref{sec-The-case}
we analyse the $\frsl_{2}(\C)$--case in greater detail.
Besides we have added two appendices, one comparing different normalizations
of the RT--invariants used in the literature and one 
discussing different definitions of framed links
in arbitrary closed oriented $3$--manifolds.

\rk{Acknowledgements} Much of this work were done while the author was
supported by a Marie Curie Fellowship of the European Commission (CEE
$\text{N}^{o}$ HPMF-CT-1999-00231).
I acknowledge hospitality of
l'Institut de Recherche Math\'{e}matique Avanc\'{e}e, Universit\'{e}
Louis Pasteur and C.N.R.S., Strasbourg, while being a Marie Curie Fellow.
The author thanks V. Turaev for having given valuable
comments to earlier versions of this paper.
He also thanks the referee for many helpful comments.
Parts of this paper are contained in the
authors thesis \cite{Hansen1}. I would here like to thank my Ph.D.\ advisor
J.\ E.\ Andersen for helpful discussions about these parts.

\section{Seifert manifolds}\label{sec-Seifert-manifolds}

Seifert manifolds were invented by H.~Seifert in \cite{Seifert1}.
For an english translation, see \cite{Seifert2}.
We consider only oriented Seifert manifolds in this section as in the
rest of the paper. These will be denoted Seifert manifolds (as in the
introduction).

\rk{Oriented Seifert manifolds and their classification}
Let $\nu,\mu$ be coprime integers with $\mu > 0$, and let
$\rho \co B^{2} \ra B^{2}$ be the rotation
by the angle $2\pi (\nu/\mu)$ in the anti-clockwise direction,
where $B^{2} \subseteq \C$
is the standard oriented unit disk.
The (oriented) {\it fibered solid torus} $T(\mu,\nu)$
is the oriented space $B^{2} \times [0,1]/_{R}$,
where $R$ identifies $(x,1)$ with $(\rho(x),0)$, $x \in B^{2}$,
and the orientation is given by the
orientation of $B^{2}$ followed by the orientation of $[0,1]$. 
By this identification the lines (fibers)
$\{ x\} \times [0,1]$ of $B^{2} \times [0,1]$,
$x \in B^{2} \sm \{ 0 \}$, are decomposed into classes,
such that each class contains
exactly $\mu$ lines, which match together to give one fiber of $T(\mu,\nu)$.
The image of $\{0 \} \times [0,1]$ in $T(\mu,\nu)$ is also a fiber, called the
`middle fiber'. 
The pair $(\mu,\nu)$ is an
invariant of $T(\mu,\nu)$ if we normalize to
$0 \leq \nu < \mu$. 
The following definition is Seifert's definition of a fibered
space \cite{Seifert1} adapted to the oriented case.
{\it A Seifert manifold
is a closed connected and oriented
$3$--manifold $M$, which
can be decomposed into a collection of disjoint simple closed
curves, called fibers, such that
each fiber $H$ has a neighborhood $N$,
called a fiber neighborhood, which is
homeomorphic to a fibered solid torus $T(\mu,\nu)$
by an orientation and fiber preserving homeomorphism mapping
$H$ to the
middle fiber of $T(\mu,\nu)$.}
By \cite[Lemma 2]{Seifert2}, the numbers $\mu,\nu$
are invariants of the fiber $H$, called the (oriented)
{\it fiber invariants} of $H$.
If $\mu>1$, we call $H$ an {\it exceptional fiber}; if
$\mu=1$, an {\it ordinary fiber}.
In a fiber neighborhood of a fiber $H$
all fibers except possibly $H$ are ordinary fibers,
so there are only finitely many (possibly zero) exceptional fibers in a
Seifert manifold. For a Seifert manifold $M$, the {\it base}
is the quotient space of $M$ obtained by identifying each 
fiber to a point. The base is a closed connected surface, orientable or
non-orientable. The genus of the non-orientable $\# ^{g} \R \text{P}^{2}$
is $g$.
Two Seifert manifolds
are {\it equivalent} if there is a fiber and orientation preserving
homeomorphism between them.
We have the following classification result due to Seifert.

\begin{thm}\label{SeifertClassification}
{\em \cite{Seifert1}}\qua An equivalence class of Seifert manifolds is determined
by a system of invariants
$$
\left( \ep;g\, | \, b;(\alpha_{1},\beta_{1}),\ldots,(\alpha_{r},\beta_{r}) \right).
$$
Here $\ep=\os$ if the
base is orientable and $\ep=\ns$ if not, and the non-negative integer
$g$ is the genus of the base.
Moreover, $r \geq 0$ is the number of exceptional fibers, and
$(\alpha_{i},\beta_{i})$ are the {\em(}oriented{\em )} Seifert invariants of the 
$i$'th exceptional fiber.
The invariant $b$ can take any value in $\Z$ {\em (}$-b$ is the Euler number
of the locally trivial $S^{1}$--bundle $( \ep;g \, | \, b)${\em )}.

An oriented Seifert manifold $M$ belonging to the class determined
by the invariants
$( \ep;g \, | \, b;(\alpha_{1},\beta_{1}),\ldots,(\alpha_{r},\beta_{r}) )$
belongs after reversing its orientation to the class determined by the invariants
$( \ep;g \, | \, -r-b; (\alpha_{1},\alpha_{1}-\beta_{1}),\ldots,(\alpha_{r},\alpha_{r}-\beta_{r}) )$,
$\ep=\os,\ns$.
\end{thm}

The Seifert invariants 
$(\alpha_{i},\beta_{i})$ of the $i$'th exceptional
fiber are the unique integers such that
$\alpha_{i}=\mu_{i}$, $\beta_{i}\nu_{i} \equiv 1 \pmod{\mu_{i}}$
and $0 < \beta_{i} < \alpha_{i}$, where
$\mu_{i},\nu_{i}$ are the fiber invariants of that fiber.
One can obtain $( \ep;g \, | \, b)$ from
$( \ep;g \, | \, b;(\alpha_{1},\beta_{1}),\ldots,(\alpha_{r},\beta_{r}) )$
by cutting out fiber neighborhoods of the exceptional fibers and gluing in
ordinary solid tori (i.e.\ $T(1,0)$'s) by certain fiber preserving
homeomorphisms,
see \cite[Sect.~7]{Seifert2}, \cite[Sect.~4.2]{Montesinos} for details.
The {\it Seifert Euler number} of the Seifert fibration
$( \ep;g \, | \, b;(\alpha_{1},\beta_{1}),\ldots,(\alpha_{r},\beta_{r}) )$
is the rational number
$e=-\left(b+\sum_{j=1}^{r} \beta_{j}/\alpha_{j}\right)$.
(The reason for the choice of sign of $e$ is the following.
Let $\ep \in \{ \os,\ns\}$ and let $X$ be a closed surface of genus $g$,
orientable if $\ep=\os$ and non-orientable if $\ep=\ns$. Then
$( \ep;g \, | \, -\chi(X))$ is the unit tangent bundle of $X$, where
$\chi(X)$ is the Euler characteristic of $X$. More generally, 
$( \ep;g \, | \, b)$ is a locally trivial $S^{1}$--bundle over the surface
$X$.
The number $-b$ is the Euler number of this bundle and
is an obstruction to the existence
of a section of $( \ep;g \, | \, b)$, see \cite[Chap.~1]{Montesinos}.
The Seifert Euler number is a natural generalization of
$-b$ when extending the above notions to orbifolds,
see \cite{Thurston}, \cite{Scott}, \cite{Montesinos}.)

\begin{figure}[ht!]\small

\begin{center}
\begin{texdraw}
\drawdim{cm}

\setunitscale 0.6

\linewd 0.02 \setgray 0

\move(0 4)

\move(0 0) \lellip rx:3 ry:1.6

\linewd 0.2 \setgray 1

\move(4 0) \larc r:1.5 sd:218 ed:200
\move(4 0) \larc r:3 sd:213 ed:204
\move(-4 0) \larc r:1.5 sd:-20 ed:150
\move(-4 0) \larc r:2 sd:-23 ed:150
\move(-4 0) \larc r:3.5 sd:-20 ed:155
\move(-4 0) \larc r:3 sd:-23 ed:155

\linewd 0.02 \setgray 0

\move(-4 0) \larc r:1.5 sd:-20 ed:150
\move(-4 0) \larc r:2 sd:-23 ed:150
\move(-4 0) \larc r:3.5 sd:-20 ed:155
\move(-4 0) \larc r:3 sd:-23 ed:155

\move(-5.88 0.67) \clvec(-5.85 0.9)(-5.6 0.95)(-5.45 0.8)
\move(-5.88 -0.67) \clvec(-5.85 -0.9)(-5.6 -0.95)(-5.45 -0.8)

\move(-5.3 0.75) \clvec(-5.35 0.55)(-5.6 0.5)(-5.75 0.7)
\move(-5.3 -0.75) \clvec(-5.35 -0.55)(-5.6 -0.5)(-5.75 -0.7)

\move(-7.3 1.17) \clvec(-7.25 1.4)(-7 1.45)(-6.85 1.3)
\move(-7.3 -1.17) \clvec(-7.25 -1.4)(-7 -1.45)(-6.85 -1.3)

\move(-6.7 1.29) \clvec(-6.75 1.09)(-7 1.04)(-7.15 1.24)
\move(-6.7 -1.29) \clvec(-6.75 -1.09)(-7 -1.04)(-7.15 -1.24)

\move(-4 0) \larc r:1.5 sd:210 ed:320
\move(-4 0) \larc r:2 sd:210 ed:318
\move(-4 0) \larc r:3.5 sd:205 ed:330
\move(-4 0) \larc r:3 sd:205 ed:328
\move(-4 0) \larc r:1.5 sd:160 ed:200
\move(-4 0) \larc r:2 sd:160 ed:200
\move(-4 0) \larc r:3.5 sd:160 ed:200
\move(-4 0) \larc r:3 sd:160 ed:200

\move(4 0) \larc r:1.5 sd:218 ed:200
\move(4 0) \larc r:3 sd:213 ed:204

\move(-6.8 -0.1) \htext{$\cdots$}
\move(1.5 -0.1) \htext{$\cdots$}
\move(-6.85 -0.8) \htext{$\undersmile{g}$}
\move(-6.65 2.8) \htext{$0$}
\move(-5.75 1.6) \htext{$0$}
\move(-7.95 0.3) \htext{$0$}
\move(-6.45 0.15) \htext{$0$}
\move(-3.95 -0.3) \htext{$-b$}
\move(5.1 1) \htext{$\frac{\alpha_{1}}{\beta_{1}}$}
\move(6.2 2.2) \htext{$\frac{\alpha_{n}}{\beta_{n}}$}

\move(0 -4)

\end{texdraw}
\end{center}

\caption{Surgery presentation of $(\os ;g \, | \, b;(\alpha_{1},\beta_{1}),\ldots,(\alpha_{n},\beta_{n}))$}\label{fig-A1}
\end{figure}

\rk{Surgery presentations}
Any closed connected oriented $3$--manifold can be
obtained by Dehn-surgery on $S^{3}$ along a labelled link,
the labels being the rational surgery coefficients,
cf.\ \cite{Lickorish1}, \cite{Wallace}.
We use the standard
convention for surgery coefficients,
see e.g.\ \cite[Chap.~9]{Rolfsen1}, \cite{Rolfsen2}.
In particular integer labelled links
in $S^{3}$ can be identified with framed links with the framing
indexes equal to
the labels. If $M$ is a $3$--manifold given by
surgery on $S^{3}$ along a labelled link $L$ we call $L$ a
surgery presentation of $M$.
According to \cite[Fig.~12 p.~146]{Montesinos}, the manifold
$( \ep ;g \, | \, b;(\alpha_{1},\beta_{1}),\ldots,(\alpha_{n},\beta_{n}))$ 
has a surgery
presentation as shown in Fig.~\ref{fig-A1} if $\ep=\os$ and as
shown in Fig.~\ref{fig-A2} if $\ep=\ns$.
The $\undersmile{g}$ indicate $g$ repetitions.

\begin{figure}[ht!]\small

\begin{center}
\begin{texdraw}
\drawdim{cm}

\setunitscale 0.6

\linewd 0.02 \setgray 0

\move(0 4)

\move(0 0) \lellip rx:3 ry:1.6
\move(-7.5 0) \lellip rx:0.5 ry:0.3
\move(-5.5 0) \lellip rx:0.5 ry:0.3

\linewd 0.2 \setgray 1
\move(-4 0) \larc r:1.5 sd:-20 ed:180
\move(-4 0) \larc r:3.5 sd:-20 ed:180
\move(4 0) \larc r:1.5 sd:215 ed:200
\move(4 0) \larc r:3 sd:212 ed:205

\linewd 0.02 \setgray 0

\move(-4 0) \larc r:1.5 sd:-20 ed:180
\move(-4 0) \larc r:3.5 sd:-20 ed:180

\move(-4 0) \larc r:1.5 sd:200 ed:320
\move(-4 0) \larc r:3.5 sd:188 ed:330

\move(4 0) \larc r:1.5 sd:215 ed:200
\move(4 0) \larc r:3 sd:212 ed:205

\move(-6.8 -0.1) \htext{$\cdots$}
\move(1.5 -0.1) \htext{$\cdots$}
\move(-6.8 -0.9) \htext{$\undersmile{g}$}
\move(-8.3 0.25) \htext{$2$}
\move(-6.3 0.25) \htext{$2$}
\move(-7.4 2.2) \htext{$\frac{1}{2}$}
\move(-5.75 1) \htext{$\frac{1}{2}$}
\move(-3.95 -0.3) \htext{$-b$}
\move(5.1 1) \htext{$\frac{\alpha_{1}}{\beta_{1}}$}
\move(6.2 2.2) \htext{$\frac{\alpha_{n}}{\beta_{n}}$}

\move(0 -4)

\end{texdraw}
\end{center}

\caption{Surgery presentation of $(\ns ; g \, | \, b;(\alpha_{1},\beta_{1}),\ldots,(\alpha_{n},\beta_{n}))$}\label{fig-A2}
\end{figure}

\rk{Non-normalized Seifert invariants}
The so-called {\it non-normalized Seifert
invariants}, see \cite{Neumann}, \cite{JankinsNeumann} or \cite{NeumannRaymond},
are sometimes more convenient to use in specific calculations.
Let $(\alpha_{j},\beta_{j})$ be a pair of coprime
integers with $\alpha_{j}>0$, $j=1,2,\ldots,n$. Then the Seifert manifold with
non-normalized Seifert invariants 
$\{\ep; g;(\alpha_{1},\beta_{1}),\ldots,(\alpha_{n},\beta_{n}) \}$
is given by a surgery presentation as shown in
Fig.~\ref{fig-A1} with $b=0$ if $\ep=\os$ and
as shown in Fig.~\ref{fig-A2} with $b=0$ if $\ep=\ns$.
It follows that
these non-normalized invariants are not unique.
In fact, by \cite[Theorem 1.5 and Theorem 1.8]{JankinsNeumann}, the sets
$\{\ep; g;(\alpha_{1},\beta_{1}),\ldots,(\alpha_{n},\beta_{n}) \}$ and
$\{\ep'; g';(\alpha_{1}',\beta_{1}'),\ldots,(\alpha_{m}',\beta_{m}') \}$
are two pairs of non-normalized Seifert invariants of the same
Seifert manifold $M$ if and only if $\ep=\ep'$, $g=g'$ (trivial),
$\sum_{i=1}^{n} \beta_{i}/\alpha_{i} = \sum_{j=1}^{m} \beta_{j}'/\alpha_{j}'$,
and disregarding any $\beta_{i}/\alpha_{i}$ and $\beta_{j}'/\alpha_{j}'$
which are integers, the remaining $\beta_{i}/\alpha_{i} \pmod{1}$ are a
permutation of the remaining $\beta_{j}'/\alpha_{j}' \pmod{1}$.
It follows that any Seifert manifold $M$ has a
unique set of non-normalized Seifert invariants (up to permutation of the
indicis) of the form
$\{ \ep; g;(1,\beta_{0}),(\alpha_{1},\beta_{1}),\ldots,(\alpha_{r},\beta_{r}) \}$
with $0<\beta_{i}<\alpha_{i}$, $i=1,\ldots,r$, so
$M=(\ep;g \, | \, \beta_{0};(\alpha_{1},\beta_{1}),\ldots,$ $(\alpha_{r},\beta_{r}))$
in the terminology of \refthm{SeifertClassification}.
This implies that the Seifert Euler number of a Seifert manifold
with non-normalized Seifert invariants 
$\{\ep; g;(\alpha_{1},\beta_{1}),\ldots,(\alpha_{n},\beta_{n}) \}$ is given by
$-\sum_{i=1}^{n} \beta_{i}/\alpha_{i}$.

\begin{rem}\label{generalizedSeifert}
\cite{JankinsNeumann}
operates with a generalization of oriented Seifert fibrations in which the
pairs $(\alpha_{j},\beta_{j})$ are allowed to be equal to $(0,\pm 1)$. However, up
to an orientation preserving homeomorphism, these generalized fibrations are
Seifert manifolds as defined above or connected sums of the form
$\#_{i=1}^{k} (S^{1} \times S^{2}) \# \#_{i=1}^{n} L(p_{i},q_{i})$,
cf.\ \cite[Theorem 5.1]{JankinsNeumann}. Since the RT--invariants behave
nicely with respect to connected sums and since the lens spaces
are (ordinary) Seifert manifolds, see the proof of \refcor{lens-spaces}, we will
continue by only considering the Seifert manifolds
in \refthm{SeifertClassification}.
\end{rem}

\section{Modular categories and $3$--manifold invariants}\label{sec-Modular-categories}

This is a preliminary section in which we recall concepts and notation from
\cite{Turaev} used throughout in this paper.
All monoidal categories in the following are assumed strict.

\rk{Ribbon categories and invariants of colored ribbon graphs}
A {\it ribbon category} $\mV$ is a monoidal category
with a braiding $c$ and a twist $\theta$ and with 
a duality $(*,b,d)$ compatible with these structures.
In $\mV$ one has a well-defined trace $\tr=\tr_{\mV}$
of morphisms and thereby a well-defined dimension $\dim=\dim_{\mV}$
of objects. These take values in the commutative semigroup
$K=K_{\mV}=\End_{\mV}(\I)$, where $\I$ is the unit object (the
multiplication being given by the composition of morphisms).

By a ($\mV$--){\it colored ribbon graph} we mean a ribbon graph
$\Omega$ with an object of $\mV$ attached to each band and annulus of $\Omega$
and with a compatible morphism of $\mV$ attached to each coupon of $\Omega$.
We let $F=F_{\mV}$ be the operator invariant of
$\mV$--colored ribbon graphs in $\R^{3}$ of Reshetikhin and Turaev,
see \cite{ReshetikhinTuraev1},
\cite{ReshetikhinTuraev2}, \cite[Chap.~I]{Turaev}.

We use the graphical calculus for morphisms of the
ribbon category $\mV$, see \cite[Sect.~I.1.6]{Turaev},
\cite[Chap.~XIV]{Kassel}.
In this calculus one
represents a morphism $f$ of $\mV$ by a colored
ribbon graph $\Omega$ mapped by $F$ to $f$ if such a
ribbon graph exists. We then write
$\Omega \doteq f$. 
We present ribbon graphs in figures according to the usual rules,
cf.\ \cite[Chap.~I]{Turaev}.
In particular we draw only the oriented
cores of the annuli and bands,
and we are careful to drawing all loops corresponding to twists in the ribbons.
Analogous to the framing numbers in figures showing framed links
we will sometimes indicate a certain number of twists in an annulus component
of a ribbon graph
by an integer instead of drawing the loops corresponding to these twists.
In figures showing colored ribbon graphs these numbers will be put into
parentheses to distinguish them from colors.

\rk{Modular categories}
A {\it monoidal Ab--category} is a monoidal category
with all morphism sets   
equipped with an additive abelian group structure making
the composition and tensor product bilinear (cf.\ \cite{MacLane};
Ab--categories are also called pre-abelian categories).

Let $\mV$ be a {\it ribbon Ab--category}, i.e.\ a ribbon category
such that the underlying monoidal category is a
monoidal Ab--category. In particular, the semigroup $K=K_{\mV}$
is a commutative unital ring, called the {\it ground ring} of $\mV$.
For any pair of objects $V$, $W$ of $\mV$, the abelian group $\Hom_{\mV}(V,W)$
acquires the structure of a left $K$--module by $kf =k \otimes f$, $k \in K$,
$f \in \Hom_{\mV}(V,W)$, which makes composition and the tensor product of
morphisms $K$--bilinear.
An object $V$ of $\mV$ is called {\it simple} if
$k \mapsto k \id_{V}$ is a bijection $K \ra \End_{\mV}(V)$.
In particular the unit
object $\I$ is simple.
An object $V$ of $\mV$ is {\it dominated}
by a family $\{ V_{i} \}_{i \in I}$ if there
exists a finite set of
morphisms $\{ f_{r} \co V_{i(r)} \ra V,\; g_{r} \co V \ra V_{i(r)} \}_{r}$
with $i(r) \in I$ such that $\id_{V} = \sum_{r} f_{r}g_{r}$.

A {\it modular category} is a tuple $\left( \mV, \{V_{i}\}_{i \in I} \right)$,
where $\mV$ is a ribbon Ab--category and $\{V_{i}\}_{i \in I}$
is a finite set of simple objects closed
under duals (i.e.\ for any $i \in I$ there exists $i^{*} \in I$ such that
$V_{i^{*}}$ is isomorphic to the dual of $V_{i}$) and dominating all objects
of $\mV$,
such that
$V_{0}=\I$ for a distinguished
element $0 \in I$, and such that the so-called
$S$--matrix $S=\left( S_{i,j} \right)_{i,j \in I}$ is invertible over
$K$. Here $S_{i,j}=\tr\left(c_{V_{j},V_{i}} \circ c_{V_{i},V_{j}} \right)$
is the invariant of the standard Hopf link
with framing $0$ and with one component colored
by $V_{i}$ and the other colored by $V_{j}$. 
The invertibility of $S$ implies that $i \mapsto i^{*}$ is an involution
in $I$.

Since $V_{i}$ is a simple object, $\theta_{V_{i}} \co V_{i} \ra V_{i}$
is equal to $v_{i}\id_{V_{i}}$ for a $v_{i} \in K$, $i \in I$.
The $T$--matrix $T=\left( T_{i,j} \right)_{i,j \in I}$ is given by
$T_{i,j} = \delta_{i,j} v_{i}$,
where $\delta_{i,j}$ is the Kronecker delta equal to $1$
if $i=j$
and to $0$ otherwise. In Fig.~\ref{fig-ST} we give a graphical description
of the entries of the $S$-- and $T$--matrices. In this and other figures
we indicate the object $V_{i}$ by $i$. Moreover,
we put $\dim(i)=\dim(V_{i})$, $i \in I$.  
We have used the identity $F(\bar{\Omega})=\tr(F(\Omega))$,
where $\bar{\Omega}$ is the
closure of a colored ribbon graph $\Omega$,
cf.\ \cite[Corollary I.2.7.2]{Turaev}.

\begin{figure}[ht!]\small

\begin{center}
\begin{texdraw}
\drawdim{cm}

\setunitscale 0.6

\linewd 0.02 \setgray 0

\move(10 0)

\move(4.3 0) \lellip rx:1.2 ry:0.6

\linewd 0.2 \setgray 1

\move(4.3 -0.2) \lvec(4.3 2.5)

\linewd 0.02 \setgray 0

\move(4.3 -0.2) \lvec(4.3 2.5)
\move(4.3 -0.8) \lvec(4.3 -2.5)
\move(4.3 -2.5) \lvec(4.15 -2.2)
\move(4.3 -2.5) \lvec(4.45 -2.2)

\move(5.5 -0.1) \lvec(5.32 0.17)
\move(5.5 -0.1) \lvec(5.65 0.17)

\move(4.4 -1.9) \htext{$j$}
\move(5.1 -1) \htext{$k$}

\move(5.9 -0.2) \htext{$ \doteq (\dim(j))^{-1}S_{k,j}\id_{V_{j}}$}

\move(-1 2.5) \lvec(-1 0) 
\clvec(-1 -0.7)(-2 -0.7)(-2 0)
\clvec(-2 0.5)(-1.15 0.6)(-1.1 0)
\move(-1 -0.5) \lvec(-1 -2.5)
\move(-1 -2.5) \lvec(-1.15 -2.2)
\move(-1 -2.5) \lvec(-0.85 -2.2)

\move(-0.7 -0.2) \htext{$\doteq v_{j}\id_{V_{j}}\; ,$}

\move(-0.9 -1.9) \htext{$j$}

\end{texdraw}
\end{center}

\nocolon
\caption{}\label{fig-ST}
\end{figure}

A {\it rank} of the modular category $\left( \mV ,\{V_{i} \}_{i \in I} \right)$
is an element $\mD=\mD_{\mV} \in K$ such that
$\mD^{2} = \sum_{i \in I} \left( \dim(i) \right)^{2}$.
A modular category does
not need to have a rank, but, as pointed out in \cite[p.~76]{Turaev}, we can always
formally change $\mV$ to a modular category with the same objects as $\mV$ and
with a rank.
We let $\Delta =\Delta_{\mV} = \sum_{i \in I} v_{i}^{-1} \left( \dim(i) \right)^{2}$.
For a modular category with a rank $\mD$ we have
\begin{equation}\label{eq:Ssquared}
S^{2} = \mD^{2} J
\end{equation}
by \cite[Formula (II.3.8.a)]{Turaev}, where
$J_{i,j}=\delta_{i^{*},j}$, $i,j \in I$.

\rk{The RT--invariants of $3$--manifolds}
We identify as usual an oriented framed link in $S^{3}=\R^{3} \cup \{ \infty \}$ 
with a ribbon graph in $S^{3}$ (actually in $\R^{3}$)
consisting solely of directed annuli, cf.\ \cite{ReshetikhinTuraev2}, \cite{Turaev}.
If $L$ is a framed link in $S^{3}$ and $B^{4}$ is the closed $4$--ball,
oriented as the unit ball in $\C^{2}$, then we get a smooth closed
connected oriented $4$--manifold $W_{L}$ by adding $2$--handles to $B^{4}$
along the components of $L$ in $S^{3}=\partial B^{4}$ using the framing
of $L$, see \cite{Kirby}. The manifold $M=M_{L}=\partial W_{L}$, oriented using
the `outward first' convention for boundaries, is the result of surgery on $S^{3}$
along $L$. Let $\Omega$ be a colored ribbon graph inside $M$
and let $\Gamma(L,\lambda)$ be the colored ribbon graph
obtained by fixing an orientation in $L$ and coloring the $i$'th
component of $L$ by $V_{\lambda(L_{i})}$. The
RT--invariant of the pair $(M,\Omega)$ based on
$\left(\mV,\{V_{i}\}_{i \in I},\mD \right)$ is given by
\begin{eqnarray}\label{eq:A4}
\tau_{(\mV,\mD)}(M,\Omega) &=& \Delta^{\sigma(L)} \mD^{-\sigma(L)-m-1} \\
 & & \hspace{.4in} \times \sum_{\lambda \in \col(L)} \left( \prod_{i=1}^{m} \dim(\lambda(L_{i})) \right) F(\Gamma(L,\lambda) \cup \Omega), \nonumber
\end{eqnarray}
cf.\ \cite[p.~82]{Turaev}, where, as usual, we identify $\Omega$ with
a colored ribbon graph in $S^{3} \sm L$. Here $m$ is the number
of components of $L$, $\sigma(L)$ is the signature of $W_{L}$, i.e.\ the
signature of the intersection form on $H_{2}(W_{L};\R)$, and
$\col(L)$ is the set of mappings from the set of components of
$L$ to $I$. The signature $\sigma(L)$
is also equal to the signature of the linking matrix of $L$.

\rk{The mirror of a modular category}
The {\it mirror} of a modular category
$\left( \mV ,\{V_{i} \}_{i \in I} \right)$
is a ribbon Ab--category $\omV$ with the same underlying monoidal Ab--category and
the same duality as $\mV$. If $\theta$ and $c$ are the twist and braiding
of $\mV$, then the twist $\bar{\theta}$ and braiding $\bar{c}$ of $\omV$
are defined by $\bar{\theta}_{V}=(\theta_{V})^{-1}$ and
$\bar{c}_{V,W}=(c_{W,V})^{-1}$ for any objects $V$, $W$ of $\mV$,
cf.\ \cite[Sect.~I.1.4]{Turaev}. By
\cite[Exercise II.1.9.2]{Turaev}, $\left( \omV, \{ V_{i} \}_{i \in I} \right)$
is a modular category with $S$--matrix $\bar{S}= \left( S_{i^{*},j} \right)_{i,j \in I}$,
where $S=\left( S_{i,j}\right)_{i,j \in I}$ is the $S$--matrix of $\mV$.
Note that $\mD$ is a rank of
$\omV$ if and only if $\mD$ is a rank of $\mV$, since the dimensions of
any object of $\mV$ with respect to $\mV$ and $\omV$ are equal, cf.\
\cite[Corollary I.2.8.5]{Turaev}. By \cite[Formula (II.2.4.a)]{Turaev}
we have
\begin{equation}\label{eq:Dsquared}
\Delta_{\mV}\Delta_{\omV}=\mD^{2}.
\end{equation}

\bigskip

We end this section by recalling the notion of a {\it unimodal modular category}
also called a {\it unimodular category}, cf.\ \cite[Sect.~VI.2]{Turaev}.
Moreover we give two small lemmas needed in the calculations of the RT--invariants
of Seifert manifolds with non-orientable base. 

Let $\left( \mV ,\{V_{i} \}_{i \in I} \right)$ be a modular category.
An element $i \in I$ is called {\it self-dual} if $i=i^{*}$. For such an element we
have a $K$--module isomorphism $\Hom_{\mV}(V \otimes V, \I) \cong K$, $V=V_{i}$.
The map $x \mapsto x(\id_{V} \otimes \theta_{V})c_{V,V}$ is a $K$--module
endomorphism of $\Hom_{\mV}(V \otimes V, \I)$, so is a multiplication by a certain
$\vep_{i} \in K$. By the definition of the braiding and twist we have
$(\vep_{i})^{2} =1$. In particular $\vep_{i} \in \{ \pm 1 \}$ if $K$ is a field.
The modular category $\left( \mV ,\{V_{i} \}_{i \in I} \right)$ is called
{\it unimodal} if $\vep_{i}=1$ for every self-dual $i \in I$. By copying a part of the
proof of \cite[Lemma VI.2.2]{Turaev} we get:

\begin{lem}\label{selfdual}
Let $\left( \mV ,\{V_{i} \}_{i \in I} \right)$ be a modular category and let
$i \in I$ be self-dual. Moreover, let $V=V_{i}$ and let $\ep_{i} \in K$
be as above. Then
\begin{equation}\label{eq:selfdual}
d_{V}(\omega \otimes \id_{V}) = \vep_{i} d_{V}^{-}(\id_{V} \otimes \omega)
\end{equation}
for any isomorphism $\omega \co V \ra V^{*}$, where $d_{V}^{-}$ is the operator
invariant $F_{\mV}$ of the left-oriented cap $\curvearrowleft$ colored with
$V$.\HS
\end{lem}

Let $(A,R,v,\{V_{i} \}_{i \in I})$ be a modular Hopf algebra over a
commutative unital ring $K$, cf.\ \cite[Chap.~XI]{Turaev}. If we write the
universal $R$--matrix as
$R=\sum_{j} \alpha_{j} \otimes \beta_{j} \in A^{\otimes 2}$, the element $u$ is
given by $u=\sum_{j} s(\beta_{j})\alpha_{j} \in A$, where $s$ is the antipode
of the underlying Hopf algebra. Let $\left( \mV ,\{V_{i} \}_{i \in I} \right)$
be the modular category induced by $(A,R,v,\{V_{i} \}_{i \in I})$, 
cf.\ \cite[Chap.~XI]{Turaev}.

\begin{lem}\label{modular-Hopf-algebra}
Let $i \in I$ be self-dual, let $V=V_{i}$ and let $\ep_{i} \in K_{\mV}=K$
be as above. For any isomorphism $\omega \co V \ra V^{*}$,
the composition

\centerline{
\xymatrix{
  {V} \ar[r]^{\omega} & {V^{*}} \ar[r]^{(\omega^{-1})^{*}} & { V^{**} } \ar[r]^{G} & {V}
}
}

is given by multiplication with $\vep_{i} uv$, where $G^{-1}$ is the canonical
$K$--module isomorphism between the finitely generated projective $K$--module $V$ and
its double dual $V^{**}$.
\end{lem}

\begin{proof}
Since $\mV$ is a ribbon category, we have a canonical $A$--module isomorphism
$\alpha_{V} \co V \ra V^{**}$ given by
$$
\alpha_{V} = (d_{V}^{-} \otimes \id_{V^{**}})(\id_{V} \otimes b_{V^{*}}),
$$
cf.\ \cite[Corollary I.2.6.1]{Turaev}. Let $Q \co V \ra V$ be multiplication
by $uv$.
Then $\alpha_{V} = G^{-1} \circ Q$. To see this, write
$b_{V^{*}}(1) = \sum_{k} g_{k} \otimes g^{k} \in V^{*} \otimes V^{**}$.
This element is characterized by the following property:
For any $\chi \in V^{*}$, $y \in V^{**}$ we have
$$
y(\chi)=\sum_{k} y(g_{k})g^{k}(\chi).
$$
Now let $x \in V$ and get
$$
\alpha_{V}(x)=\sum_{k} d_{V}^{-}(x \otimes g_{k}) \otimes g^{k}.
$$
By using that $d_{V}^{-}=d_{V}c_{V,V^{*}}(\theta_{V} \otimes \id_{V^{*}})$
we get
$$
\alpha_{V}(x)=\sum_{k} g_{k}(uv \cdot x)g^{k} \in V^{**}.
$$
If $\chi \in V^{*}$ we therefore have
$$
\alpha_{V}(x)(\chi)=\sum_{k} g_{k}(uv \cdot x)g^{k}(\chi ) = G^{-1} \circ Q(x)(\chi ).
$$
If $f \co U \ra W$ is a morphism in $\mV$, then the dual
morphism $f^{*} \co W^{*} \ra U^{*}$ is given by
$f^{*}=(d_{W} \otimes \id_{U^{*}})(\id_{W^{*}} \otimes f \otimes \id_{U^{*}})(\id_{W^{*}} \otimes b_{U})$.
By using the graphical calculus together with (\ref{eq:selfdual}) one
immediately gets that $(\omega^{-1})^{*} \circ \omega =\vep_{i} \alpha_{V}$.
\end{proof}

\section{The Reshetikhin--Turaev invariants of Seifert manifolds}\label{sec-The-Reshetikhin}

In this section we calculate the RT--invariants of all oriented Seifert
manifolds.
Throughout, $\left(\mV,\{V_{i}\}_{i \in I} \right)$ is a
fixed modular category with a fixed rank $\mD$. We let $F=F_{\mV}$,
$\Delta=\Delta_{\mV}$, $K=K_{\mV}$, and $\tau=\tau_{(\mV,\mD)}$.

\rk{Notation} For the next theorem and for later use we introduce some notation.
Let $y(i,j) \in K$ be the scalar such that $F(T_{ij})=y(i,j)\id_{V_{j}}$,
where $T_{ij}$ is the colored ribbon tangle in Fig.~\ref{fig-A8}. That is,
$y(i,j)=(\dim(j))^{-1}\tr(F(T_{ij}))$.
We put
\begin{equation}\label{eq:vector}
\kappa(j) = \sum_{i \in I} \dim(i)y(i,j), \hspace{.2in} j \in I.
\end{equation}
For every self-dual element $i \in I$, let $\vep_{i} \in K$ be
as in the last part of Sect.~\ref{sec-Modular-categories}.

\begin{figure}[ht!]\small

\begin{center}
\begin{texdraw}
\drawdim{cm}

\setunitscale 0.6

\linewd 0.02 \setgray 0

\move(8 0)

\move(-4 8)

\move(-4 7) \lvec(-4 5)

\move(-3.5 4) \clvec(-3.5 3.5)(-4 3.4)(-4 3)

\move(-3.5 2) \clvec(-3.5 1.5)(-4 1.4)(-4 1)

\move(-3.5 5) \clvec(-3.5 4.5)(-4 4.4)(-4 4)

\move(-3.5 3) \clvec(-3.5 2.5)(-4 2.4)(-4 2)

\move(-4 1) \lvec(-4 -1)

\move(-4 -0.3) \lvec(-4.15 0)
\move(-4 -0.3) \lvec(-3.85 0)

\linewd 0.2 \setgray 1

\move(-4 4.9) \clvec(-4 4.5)(-3.5 4.4)(-3.5 4.1)

\move(-4 2.9) \clvec(-4 2.5)(-3.5 2.4)(-3.5 2.1)

\move(-4 3.9) \clvec(-4 3.5)(-3.5 3.4)(-3.5 3.1)

\move(-4 1.9) \clvec(-4 1.5)(-3.5 1.4)(-3.5 1.1)

\linewd 0.02 \setgray 0

\move(-4 5) \clvec(-4 4.5)(-3.5 4.4)(-3.5 4)

\move(-4 3) \clvec(-4 2.5)(-3.5 2.4)(-3.5 2)

\move(-4 4) \clvec(-4 3.5)(-3.5 3.4)(-3.5 3)

\move(-4 2) \clvec(-4 1.5)(-3.5 1.4)(-3.5 1)

\move(-2.5 5) \larc r:1 sd:0 ed:180
\move(-2.5 1) \larc r:1 sd:180 ed:0

\move(-1.5 5) \lvec(-1.5 1)

\move(-1.5 1.95) \lvec(-1.65 1.65)
\move(-1.5 1.95) \lvec(-1.35 1.65)

\move(-3.8 -0.6) \htext{$j$}
\move(-1.3 1.8) \htext{$i$}

\move(-3 -2) \htext{$T_{ij}$}

\move(0 -1.6) \htext{,}

\linewd 0.02 \setgray 0

\move(5 8)

\move(5 7) \lvec(5 5)

\move(5.5 4) \clvec(5.5 3.5)(5 3.4)(5 3)

\move(5.5 2) \clvec(5.5 1.5)(5 1.4)(5 1)

\move(5.5 5) \clvec(5.5 4.5)(5 4.4)(5 4)

\move(5.5 3) \clvec(5.5 2.5)(5 2.4)(5 2)

\move(5 1) \lvec(5 -1)

\move(5 -0.3) \lvec(4.85 0)
\move(5 -0.3) \lvec(5.15 0)

\linewd 0.2 \setgray 1

\move(5 4.9) \clvec(5 4.5)(5.5 4.4)(5.5 4.1)

\move(5 2.9) \clvec(5 2.5)(5.5 2.4)(5.5 2.1)

\move(5 3.9) \clvec(5 3.5)(5.5 3.4)(5.5 3.1)

\move(5 1.9) \clvec(5 1.5)(5.5 1.4)(5.5 1.1)

\linewd 0.02 \setgray 0

\move(5 5) \clvec(5 4.5)(5.5 4.4)(5.5 4)

\move(5 3) \clvec(5 2.5)(5.5 2.4)(5.5 2)

\move(5 4) \clvec(5 3.5)(5.5 3.4)(5.5 3)

\move(5 2) \clvec(5 1.5)(5.5 1.4)(5.5 1)

\move(6.5 5) \larc r:1 sd:0 ed:180
\move(6.5 1) \larc r:1 sd:180 ed:0

\move(7.5 5) \lvec(7.5 1)

\move(7.5 1.95) \lvec(7.35 1.65)
\move(7.5 1.95) \lvec(7.65 1.65)

\move(5.2 -0.6) \htext{$j$}
\move(7.7 1.2) \htext{$i$}

\move(1.0 2.5) \htext{$\sum_{i \in I} \dim(i)$}

\move(7.9 2.5) \htext{$\hspace{.1in} \doteq \hspace{.1in} \kappa(j) \id_{V_{j}}$}

\end{texdraw}
\end{center}

\nocolon
\caption{}\label{fig-A8}
\end{figure}

The group $SL(2,\Z)$ is generated by two matrices
\begin{equation}\label{eq:generators}
\Xi = \left( \begin{array}{cc}
		0 & -1 \\
		1 & 0
		\end{array}
\right),\hspace{.2in}
\Theta = \left( \begin{array}{cc}
                1 & 1 \\
		0 & 1
		\end{array}
\right). 
\end{equation}
For a tuple of integers $\mC=(a_{1},\ldots,a_{n})$
we let
\begin{equation}\label{eq:Bmatrix}
B_{k}^{\mC} = \left( \begin{array}{cc}
                \alpha_{k}^{\mC} & \rho_{k}^{\mC} \\
		\beta_{k}^{\mC} & \sigma_{k}^{\mC}
		\end{array}
\right) = \Theta^{a_{k}}\Xi\Theta^{a_{k-1}}\Xi \ldots \Theta^{a_{1}}\Xi 
,\hspace{.2in} k=1,2,\ldots,n
\end{equation}
and let $B^{\mC}=B_{n}^{\mC}$. Moreover, we put
\begin{equation}\label{eq:A3}
G^{\mC}=T^{a_{n}}ST^{a_{n-1}}S \cdots ST^{a_{1}}S.
\end{equation}
A continued fraction expansion
$$
\frac{p}{q}=a_{n}-\frac{1}{a_{n-1}-\dfrac{1}{\cdots -\dfrac{1}{a_{1}}}}, \hspace{.2in} a_{i} \in \Z,
$$
$p,q \in \Z$ not both equal to zero, is abbreviated$(a_{1},\ldots,a_{n})$.
Given pairs $(\alpha_{j},\beta_{j})$ of coprime integers we
let $\mC_{j}=(a_{1}^{(j)},a_{2}^{(j)},\ldots,a_{m_{j}}^{(j)})$ be a continued
fraction expansion of $\alpha_{j}/\beta_{j}$,
$j=1,2,\ldots,n$.

\begin{thm}\label{invariants}
The RT--invariant $\tau$ of 
$M=(\os;g \, | \, b;(\alpha_{1},\beta_{1}),\ldots,(\alpha_{n},\beta_{n}))$
is
\begin{equation}\label{eq:osurgery}
\tau(M)= (\Delta\mD^{-1})^{\sigma_{\os}} \mD^{2g-2-\sum_{j=1}^{n} m_{j}} \sum_{j \in I} v_{j}^{-b} \dim(j)^{2-n-2g} \left( \prod_{i=1}^{n} (SG^{\mC_{i}})_{j,0} \right),
\end{equation}
where
\begin{equation}\label{eq:osignature0}
\sigma_{\os}=\sign(e) + \sum_{j=1}^{n} \sum_{l=1}^{m_{j}} \sign(\alpha_{l}^{\mC_{j}}\beta_{l}^{\mC_{j}}).
\end{equation}
Here $e=-\left( b+\sum_{j=1}^{n} \frac{\beta_{j}}{\alpha_{j}} \right)$ is the Seifert Euler
number. 

The RT--invariant $\tau$ of the Seifert manifold $M$ with
non-normalized Seifert invariants
$\{ \os;g;(\alpha_{1},\beta_{1}),\ldots,(\alpha_{n},\beta_{n})\}$
is given by the same expression with the exceptions, that the factor $v_{j}^{-b}$
has to be removed and $e=-\sum_{j=1}^{n} \frac{\beta_{j}}{\alpha_{j}}$.

The RT--invariant $\tau$ of
$M=(\ns ;g \, | \, b;(\alpha_{1},\beta_{1}),\ldots,(\alpha_{n},\beta_{n}))$
is
\begin{eqnarray}\label{eq:nsurgery}
\tau(M) &=& (\Delta\mD^{-1})^{\sigma_{\ns}} \mD^{g-2-\sum_{j=1}^{n} m_{j}} \\
 & & \hspace{.4in} \times \sum_{j \in I} \left(\vep_{j}\right)^{g} \delta_{j,j^{*}} v_{j}^{-b} \dim(j)^{2-n-g} \left( \prod_{i=1}^{n} (SG^{\mC_{i}})_{j,0} \right), \nonumber 
\end{eqnarray}
where $\delta_{j,k}$ is the Kronecker delta equal to $1$ if $j=k$ and
to $0$ otherwise, and
\begin{equation}\label{eq:nsignature0}
\sigma_{\ns}  = \sum_{j=1}^{n} \sum_{l=1}^{m_{j}} \sign(\alpha_{l}^{\mC_{j}}\beta_{l}^{\mC_{j}}). 
\end{equation}

The RT--invariant $\tau$ of the Seifert manifold $M$ with
non-normalized Seifert invariants
$\{\ns;g;(\alpha_{1},\beta_{1}),\ldots,(\alpha_{n},\beta_{n})\}$
is given by the same expression with the exception, that the factor $v_{j}^{-b}$
has to be removed.
\end{thm}

The theorem is also valid in case $n=0$. In this case one just
has to put all sums $\sum_{j=1}^{n}$
equal to zero and all products $\prod_{i=1}^{n}$ equal to $1$.
Note that $\ep_{j}^{g}=1$ if $g$ is even and $\ep_{j}^{g}=\ep_{j}$
if $g$ is odd since $\ep_{j}^{2}=1$.

\rk{Preliminaries} Before giving the proof of \refthm{invariants}
we make some preliminary remarks.

{\bf 1)} Let $\mC=(a_{1},\ldots,a_{n}) \in \Z^{n}$ and consider the matrices in
(\ref{eq:Bmatrix}). By \cite[Proposition 2.5]{Jeffrey}
we have that $(a_{1},\ldots,a_{k})$ is a continued
fraction expansion of $\alpha_{k}^{\mC}/\beta_{k}^{\mC}$,
$k=1,2,\ldots,n$, and that $\beta_{k}^{\mC}=\alpha_{k-1}^{\mC}$,
$k=2,3,\ldots,n$. Note that $\alpha_{1}^{\mC}=a_{1}$ and
$\beta_{1}^{\mC}=1$.

\begin{figure}[ht!]\small

\begin{center}
\begin{texdraw}
\drawdim{cm}

\setunitscale 0.6

\linewd 0.02

\move(0 3)
\move(0 -2)

\move(-9 0) \larc r:3 sd:140 ed:220

\move(-12.5 1.0) \htext{$\frac{p}{q}$}
\move(-12.3 -2.3) \htext{$L_{i}$}

\move(-10 -0.1) \htext{$\sim$}

\move(0 0) \larc r:1 sd:-45 ed:130
\move(0 0) \larc r:1 sd:150 ed:300
\move(-1.5 0) \larc r:1 sd:-35 ed:100
\move(-1.5 0) \larc r:1 sd:260 ed:310

\move(-3 -0.1) \htext{$\cdots$}

\move(-4.5 0) \larc r:1 sd:70 ed:130
\move(-4.5 0) \larc r:1 sd:150 ed:300
\move(-6 0) \larc r:1 sd:-35 ed:130
\move(-6 0) \larc r:1 sd:150 ed:310
\move(-7.5 0) \larc r:1 sd:-35 ed:310

\move(-7.8 1.1) \htext{$a_{1}$}
\move(-6.3 1.1) \htext{$a_{2}$}
\move(-4.8 1.1) \htext{$a_{3}$}
\move(-2.5 1.1) \htext{$a_{n-2}$}
\move(-0.8 1.1) \htext{$a_{n-1}$}

\move(3.5 0) \larc r:3 sd:140 ed:160
\move(3.5 0) \larc r:3 sd:170 ed:220

\move(0.3 1.8) \htext{$a_{n}$}
\move(0.5 -2.3) \htext{$L_{i}$}

\lpatt(0.067 0.1)

\move(3.5 0) \larc r:3 sd:125 ed:140
\move(3.5 0) \larc r:3 sd:220 ed:235

\move(-9 0) \larc r:3 sd:125 ed:140
\move(-9 0) \larc r:3 sd:220 ed:235

\end{texdraw}
\end{center}

\nocolon
\caption{}\label{fig-A4}
\end{figure} 

{\bf 2)} Two labelled links in $S^{3}$ are (surgery)
equivalent if surgeries on $S^{3}$ along these labelled links result in
$3$--manifolds which are isomorphic as oriented $3$--manifolds.
Correspondingly we talk about equivalent surgery presentations. 
We have the following well-known fact \cite[p.~273]{Rolfsen1}:
Let $(a_{1},\ldots,a_{n})$
be a continued fraction expansion of
$p/q \in \Q$ and let $L$ be a labelled link with a component
$L_{i}$ with surgery coefficient $p/q$. Then this link is surgery
equivalent to a link obtained from $L$ by 
changing the surgery coefficient of $L_{i}$
to $a_{n}$ and shackling $L_{i}$ with an integer labelled Hopf chain with $n-1$
components with labels $a_{1},\ldots,a_{n-1}$ as shown
in Fig.~\ref{fig-A4}. For a proof of this, simply use standard
surgery modifications, cf.\ \cite[Sect.~9.H]{Rolfsen1}, \cite{Rolfsen2}.
(Begin by unknotting $L_{i}$ in the presentation in the right-hand
side of Fig.~\ref{fig-A4} and get rid of the Hopf chain, see the proof of
\cite[Proposition 17.3]{PrasolovSosinski}. Finally recover the
original $L_{i}$ by knotting.)
Alternatively, see \cite[Appendix]{KirbyMelvin2}.

\begin{figure}[ht!]\small

\begin{center}
\begin{texdraw}
\drawdim{cm}

\setunitscale 0.6

\linewd 0.02 \setgray 0

\move(0 0) \lellip rx:1.2 ry:0.6
\move(-0.15 -0.6) \lvec(0.15 -0.75)
\move(-0.15 -0.6) \lvec(0.15 -0.45)
\move(0.9 -0.9) \htext{$k$}

\linewd 0.2 \setgray 1

\move(-0.5 -0.2) \lvec(-0.5 2.5)
\move(0.5 -0.2) \lvec(0.5 2.5)

\linewd 0.02 \setgray 0

\move(-0.5 -0.2) \lvec(-0.5 2.5)
\move(0.5 -0.2) \lvec(0.5 2.5)

\move(-0.5 -0.8) \lvec(-0.5 -2.5)
\move(0.5 -0.8) \lvec(0.5 -2.5)

\move(-0.5 -2.5) \lvec(-0.65 -2.2)
\move(-0.5 -2.5) \lvec(-0.35 -2.2)
\move(0.5 2.5) \lvec(0.35 2.2)
\move(0.5 2.5) \lvec(0.65 2.2)

\move(-6.8 -0.3) \htext{$\dim(i) \sum_{k \in I} \dim(k)$}
\move(1.5 -0.3) \htext{$\doteq \hspace{.05in} \delta_{i,j} \mD^{2}$}

\move(-0.9 -1.9) \htext{$j$}
\move(0.6 1.65) \htext{$i$}

\move(4.8 0.8) \larc r:0.5 sd:180 ed:0
\move(4.8 -0.8) \larc r:0.5 sd:0 ed:180

\move(4.3 0.8) \lvec(4.3 2.5)
\move(5.3 0.8) \lvec(5.3 2.5)
\move(4.3 -0.8) \lvec(4.3 -2.5)
\move(5.3 -0.8) \lvec(5.3 -2.5)

\move(4.3 -2.5) \lvec(4.15 -2.2)
\move(4.3 -2.5) \lvec(4.45 -2.2)

\move(5.3 2.5) \lvec(5.15 2.2)
\move(5.3 2.5) \lvec(5.45 2.2)

\move(5.4 1.65) \htext{$i$}
\move(3.9 -1.9) \htext{$i$}

\end{texdraw}
\end{center}

\nocolon
\caption{}\label{fig-A9}
\end{figure}

\begin{figure}[ht!]\small

\begin{center}
\begin{texdraw}
\drawdim{cm}

\setunitscale 0.6

\linewd 0.02 \setgray 0

\move(0 0) \lellip rx:1.2 ry:0.6
\move(-0.15 -0.6) \lvec(0.15 -0.75)
\move(-0.15 -0.6) \lvec(0.15 -0.45)
\move(1.1 -0.8) \htext{$k$}

\linewd 0.2 \setgray 1

\move(-0.5 -0.2) \lvec(-0.5 3)
\move(0.5 -0.2) \lvec(0.5 3)

\linewd 0.02 \setgray 0

\move(-0.5 -0.2) \lvec(-0.5 3)
\move(0.5 -0.2) \lvec(0.5 3)

\move(-0.5 -0.8) \lvec(-0.5 -3)
\move(0.5 -0.8) \lvec(0.5 -3)

\move(-0.5 -3) \lvec(-0.65 -2.7)
\move(-0.5 -3) \lvec(-0.35 -2.7)
\move(0.5 3) \lvec(0.35 2.7)
\move(0.5 3) \lvec(0.65 2.7)

\move(1.5 -0.3) \htext{$\doteq \hspace{.05in} \sum_{l}$}

\move(4.2 0) \lellip rx:0.8 ry:0.4

\linewd 0.2 \setgray 1

\move(4.2 1.4) \lvec(4.2 0)

\linewd 0.02 \setgray 0

\move(4.2 1.4) \lvec(4.2 0)

\move(4.2 -1.4) \lvec(4.2 -0.6)

\move(3.5 1.4) \lvec(4.9 1.4) \lvec(4.9 2.3) \lvec(3.5 2.3) \lvec(3.5 1.4)

\move(3.5 -1.4) \lvec(4.9 -1.4) \lvec(4.9 -2.3) \lvec(3.5 -2.3) \lvec(3.5 -1.4)

\move(3.8 2.3) \lvec(3.8 3)
\move(4.6 2.3) \lvec(4.6 3)
\move(4.6 3) \lvec(4.45 2.7)
\move(4.6 3) \lvec(4.75 2.7)
\move(3.8 2.3) \lvec(3.65 2.6)
\move(3.8 2.3) \lvec(3.95 2.6)

\move(4.2 -1.4) \lvec(4.05 -1.1)
\move(4.2 -1.4) \lvec(4.35 -1.1)

\move(5.0 0) \lvec(4.8 0.1)
\move(5.0 0) \lvec(5.1 0.2)

\move(3.8 -2.3) \lvec(3.8 -3)
\move(3.8 -3) \lvec(3.65 -2.7)
\move(3.8 -3) \lvec(3.95 -2.7)

\move(4.6 -2.3) \lvec(4.6 -3)
\move(4.6 -2.3) \lvec(4.45 -2.6)
\move(4.6 -2.3) \lvec(4.75 -2.6)

\move(-1 -2.2) \htext{$j$}
\move(0.6 1.95) \htext{$i$}

\move(4.4 -1.2) \htext{$i(l)$}

\move(3.9 1.6) \htext{$f_{l}$}
\move(3.9 -2) \htext{$g_{l}$}

\move(5.15 0) \htext{$k$}

\move(4.8 2.5) \htext{$i$}
\move(4.8 -2.9) \htext{$i$}

\move(3.3 2.5) \htext{$j$}
\move(3.3 -2.9) \htext{$j$}

\end{texdraw}
\end{center}

\nocolon
\caption{}\label{fig-A9extra}
\end{figure}

{\bf 3)} The identity in Fig.~\ref{fig-A9} is due to Turaev,
cf.\ \cite[Exercise II.3.10.2]{Turaev}. For the sake of completeness
we give a proof of it here.

\begin{proof}[Proof of the identity in Fig.~\ref{fig-A9}]
The axiom of domination for a modular category, see
Sect.~\ref{sec-Modular-categories}, implies that
we for arbitrary $i,j \in I$ can write the identity endomorphism of
$V_{j} \otimes V_{i}^{*}$ as a finite sum
$$
\id_{V_{j} \otimes V_{i}^{*}}=\sum_{l} f_{l}g_{l},
$$
where $f_{l} \co V_{i(l)} \ra V_{j} \otimes V_{i}^{*}$ and
$g_{l} \co V_{j} \otimes V_{i}^{*} \ra V_{i(l)}$ are certain
morphisms. By this we get the identity in Fig.~\ref{fig-A9extra}.
According to \cite[Lemma II.3.2.3]{Turaev} we have $\dim(k)=d_{0}^{-1}d_{k}$,
where the elements $d_{i} \in K$, $i \in I$, are defined by
(\ref{eq:d}), cf.\ \cite[p.~87]{Turaev}.
By using this and \cite[Lemma II.3.2.2 (i)]{Turaev} we get
the identity shown in Fig.~\ref{fig-A9extra1}, where
$x=\sum_{u \in I} d_{u}\dim(u)$. Since
$V_{i}^{*} \cong V_{i^{*}}$ and
$\Hom(\I,V_{j} \otimes V_{i^{*}})=0$
unless $i=j$, cf.\ \cite[Lemma II.3.5]{Turaev}, we get the
result for $i \neq j$.
Assume $i=j$. By \cite[Lemma II.3.5]{Turaev},
the $K$--module 
$\Hom(\I,V_{i} \otimes V_{i}^{*})$
is generated by $b_{V_{i}}$, where
$b$ is part of
the duality of the modular category. Similarly,
$\Hom(V_{i} \otimes V_{i}^{*}, \I)$
is generated by $d_{V_{i}}^{-}$, where
$d_{V_{i}}^{-} \co V_{i} \otimes V_{i}^{*} \ra \I$ is the operator invariant
of the left-oriented cap $\curvearrowleft$ colored with $V_{i}$.
We can therefore write
$\id_{V_{i} \otimes V_{i}^{*}}=fg+\sum_{l:i(l) \neq 0} f_{l}g_{l}$, where
$f=ab_{V_{i}}$ and $g=a'd_{V_{i}}^{-}$, $a,a' \in K$.
By this we get
$$
d_{V_{i}}^{-}b_{V_{i}}=d_{V_{i}}^{-}\id_{V_{i} \otimes V_{i}^{*}}b_{V_{i}}=aa'(d_{V_{i}}^{-}b_{V_{i}})^{2},
$$
since $\Hom(V_{r},V_{s})=0$ for any distinct $r,s \in I$,
cf.\ \cite[Lemma II.1.5]{Turaev}.
Since $d_{V_{i}}^{-}b_{V_{i}}=\dim(i)$ we get
$aa'=(\dim(i))^{-1}$.
Combining this with the identity in
Fig.~\ref{fig-A9extra1} and the fact that $xd_{0}^{-1}=\mD^{2}$,
cf.\ \cite[p.~89]{Turaev}, finally brings us to the identity
in Fig.~\ref{fig-A9}.
\end{proof}

The above proof does not
use the existence of a rank. In case we don't have a rank,
the identity in Fig.~\ref{fig-A9} is still valid if we replace $\mD^{2}$ by
$\sum_{u \in I} \left( \dim(u) \right)^{2}$. One should also note that
the orientation of the annulus component with color $k$ does
not play any role. This follows by the usual argument since we sum over
all colors $k$.

\begin{figure}[ht!]\small

\begin{center}
\begin{texdraw}
\drawdim{cm}

\setunitscale 0.6

\linewd 0.02 \setgray 0

\move(0 3.5)

\move(0 0) \lellip rx:1.2 ry:0.6
\move(-0.15 -0.6) \lvec(0.15 -0.75)
\move(-0.15 -0.6) \lvec(0.15 -0.45)
\move(1.1 -0.8) \htext{$k$}

\move(-0.5 -0.2) \lvec(-0.5 2.8)
\move(0.5 -0.2) \lvec(0.5 2.8)

\move(-0.5 -0.8) \lvec(-0.5 -2.8)
\move(0.5 -0.8) \lvec(0.5 -2.8)

\move(-0.5 -2.8) \lvec(-0.65 -2.5)
\move(-0.5 -2.8) \lvec(-0.35 -2.5)
\move(0.5 2.8) \lvec(0.35 2.5)
\move(0.5 2.8) \lvec(0.65 2.5)

\move(1.5 -0.3) \htext{$\doteq \hspace{.05in} xd_{0}^{-1} \sum_{l:i(l)=0}$}

\linewd 0.02 \setgray 0

\move(6.5 1.1) \lvec(6.5 -1.1)

\move(5.8 1.1) \lvec(7.2 1.1) \lvec(7.2 2) \lvec(5.8 2) \lvec(5.8 1.1)

\move(5.8 -1.1) \lvec(7.2 -1.1) \lvec(7.2 -2) \lvec(5.8 -2) \lvec(5.8 -1.1)

\move(6.1 2) \lvec(6.1 2.8)
\move(6.9 2) \lvec(6.9 2.8)
\move(6.9 2.8) \lvec(6.75 2.5)
\move(6.9 2.8) \lvec(7.05 2.5)
\move(6.1 2) \lvec(5.95 2.3)
\move(6.1 2) \lvec(6.25 2.3)

\move(6.5 -1.1) \lvec(6.35 -0.8)
\move(6.5 -1.1) \lvec(6.65 -0.8)

\move(6.1 -2) \lvec(6.1 -2.8)
\move(6.1 -2.8) \lvec(5.95 -2.5)
\move(6.1 -2.8) \lvec(6.25 -2.5)

\move(6.9 -2) \lvec(6.9 -2.8)
\move(6.9 -2) \lvec(6.75 -2.3)
\move(6.9 -2) \lvec(7.05 -2.3)

\move(-4.9 -0.3) \htext{$\sum_{k \in I} \dim(k)$}

\move(-1 -2.1) \htext{$j$}
\move(0.6 1.85) \htext{$i$}

\move(6.7 -0.9) \htext{$0$}

\move(6.2 1.3) \htext{$f_{l}$}
\move(6.2 -1.7) \htext{$g_{l}$}

\move(7.1 2.2) \htext{$i$}
\move(7.1 -2.6) \htext{$i$}

\move(5.6 2.2) \htext{$j$}
\move(5.6 -2.6) \htext{$j$}

\end{texdraw}
\end{center}

\nocolon
\caption{}\label{fig-A9extra1}
\end{figure}

\begin{proof}[Proof of \refthm{invariants}]
Let
$M=(\os;g \, | \, b;(\alpha_{1},\beta_{1}),\ldots,(\alpha_{n},\beta_{n}))$
and let
$L$ be the link obtained from the
link in Fig.~\ref{fig-A1} by replacing the component with surgery coefficient
$\alpha_{j}/\beta_{j}$ by a chain according to $\mC_{j}$
as in Fig.~\ref{fig-A4}, $j=1,\ldots,n$.
Note that $L$ has $m=2g+1+\sum_{j=1}^{n} m_{j}$ components.
By (\ref{eq:A4}) and the identities in Fig.~\ref{fig-ST} we have
\begin{eqnarray*}
\tau(M) &=& \Delta^{\sigma(L)} \mD^{-\sigma(L)-m-1} \sum_{j \in I} \dim(j)^{1-n} \left( \prod_{i=1}^{n} (SG^{\mC_{i}})_{j,0} \right) \\
 & & \times\kern -5pt \sum_{u_{1},\ldots,u_{g},n_{1},\ldots,n_{g} \in I} \left( \prod_{l=1}^{g} \dim(u_{l})\dim(n_{l}) \right) F(\Gamma(j,u_{1},n_{1},\ldots,u_{g},n_{g})), 
\end{eqnarray*}
where
$\Gamma(j,u_{1},n_{1},\ldots,u_{g},n_{g})$ is the colored ribbon
graph shown in Fig.~\ref{fig-A6}.
If $g=0$ we have to replace the sum
$\sum_{u_{1},\ldots,u_{g},n_{1},\ldots,n_{g} \in I}$
by $v_{j}^{-b}\dim(j)$ here and can go directly
to the calculation of $\sigma(L)$. Assume $g>0$.  
By using the identity in Fig.~\ref{fig-A9} with the component colored with
$u_{j}$ in Fig.~\ref{fig-A6} equal to the component colored with $k$ in
Fig.~\ref{fig-A9}, $j=1,2,\ldots,g$, we get\eject
\begin{eqnarray*}
& & \sum_{u_{1},\ldots,u_{g},n_{1},\ldots,n_{g} \in I} \left( \prod_{l=1}^{g} \dim(u_{l})\dim(n_{l}) \right) F(\Gamma(j,u_{1},n_{1},\ldots,u_{g},n_{g})) \\
 & & \hspace{1.0in}= \mD^{2g} \sum_{n_{1},\ldots,n_{g} \in I} F(\Gamma(j,n_{1},\ldots,n_{g})),
\end{eqnarray*}
where $\Gamma(j,n_{1},\ldots,n_{g})$ is the colored ribbon tangle shown in
Fig.~\ref{fig-A10}.
The expression (\ref{eq:osurgery}) now follows by
the fact that $\sigma(L)=\sigma_{\os}$, see below, and by
\begin{eqnarray*}
&&\sum_{n_{1},\ldots,n_{g} \in I} F(\Gamma(j,n_{1},\ldots,n_{g})) \\
 & & \hspace{0.3in} = v_{j}^{-b} \sum_{n_{1},\ldots,n_{g} \in I} \left( \prod_{l=1}^{g} S_{n_{l},j}S_{n_{l}^{*},j}\dim(j)^{-2} \right) \dim(j) =\dim(j)^{1-2g} v_{j}^{-b}\mD^{2g},
\end{eqnarray*}
where the first equality follows by the identities
in Fig.~\ref{fig-ST} and the last equality follows
by (\ref{eq:Ssquared}) and the facts that $S$ is symmetric
and satisfies $S_{i,j}=S_{i^{*},j^{*}}$, $i,j \in I$, 
cf.\ \cite[Formula (II.3.3.a)]{Turaev}.

\begin{figure}[ht!]\small

\begin{center}
\begin{texdraw}
\drawdim{cm}

\setunitscale 0.6

\linewd 0.02 \setgray 0

\move(0 4.5)

\move(-1.5 0) \lellip rx:3 ry:1.6

\move(-4.5 -0.1) \lvec(-4.65 0.2)
\move(-4.5 -0.1) \lvec(-4.35 0.2)

\linewd 0.2 \setgray 1

\move(-4 0) \larc r:1.5 sd:308 ed:150
\move(-4 0) \larc r:2 sd:318 ed:150
\move(-4 0) \larc r:3.7 sd:338 ed:157
\move(-4 0) \larc r:3.2 sd:333 ed:157

\linewd 0.02 \setgray 0

\move(-4 0) \larc r:1.5 sd:308 ed:150
\move(-4 0) \larc r:2 sd:318 ed:150
\move(-4 0) \larc r:3.7 sd:338 ed:157
\move(-4 0) \larc r:3.2 sd:333 ed:157

\move(-2 0) \lvec(-2.15 -0.3)
\move(-2 0) \lvec(-1.85 -0.3)
\move(-0.3 0) \lvec(-0.45 -0.3)
\move(-0.3 0) \lvec(-0.15 -0.3)

\move(-5.88 0.67) \clvec(-5.85 0.9)(-5.6 0.95)(-5.45 0.8)
\move(-5.88 -0.67) \clvec(-5.85 -0.9)(-5.6 -0.95)(-5.45 -0.8)

\move(-5.3 0.75) \clvec(-5.35 0.55)(-5.6 0.5)(-5.75 0.7)
\move(-5.3 -0.75) \clvec(-5.35 -0.55)(-5.6 -0.5)(-5.75 -0.7)

\move(-7.5 1.17) \clvec(-7.45 1.4)(-7.2 1.45)(-7.05 1.3)
\move(-7.5 -1.17) \clvec(-7.45 -1.4)(-7.2 -1.45)(-7.05 -1.3)

\move(-6.92 1.29) \clvec(-6.95 1.09)(-7.2 1.04)(-7.35 1.24)
\move(-6.92 -1.29) \clvec(-6.95 -1.09)(-7.2 -1.04)(-7.35 -1.24)

\move(-4 0) \larc r:1.5 sd:210 ed:294
\move(-4 0) \larc r:2 sd:210 ed:306
\move(-4 0) \larc r:3.7 sd:203 ed:333
\move(-4 0) \larc r:3.2 sd:203 ed:326
\move(-4 0) \larc r:1.5 sd:160 ed:200
\move(-4 0) \larc r:2 sd:160 ed:200
\move(-4 0) \larc r:3.7 sd:160 ed:200
\move(-4 0) \larc r:3.2 sd:161 ed:199

\move(-5.5 0) \lvec(-5.65 -0.3)
\move(-5.5 0) \lvec(-5.35 -0.3)
\move(-7.2 0) \lvec(-7.35 -0.3)
\move(-7.2 0) \lvec(-7.05 -0.3)

\move(-6.9 -0.1) \htext{$\cdots$}

\move(-4.2 3.9) \htext{$n_{g}$}

\move(-4.1 2.1) \htext{$n_{1}$}

\move(-8.4 -0.5) \htext{$u_{g}$}

\move(-6.6 -0.55) \htext{$u_{1}$}
\move(1.6 -0.3) \htext{$(-b)$}
\move(1.4 -1) \htext{$j$}

\move(0 -4)

\end{texdraw}
\end{center}

\caption{Oriented colored surgery presentation of $(\os;g \, | \, b)$}\label{fig-A6}
\end{figure}

\begin{figure}[ht!]
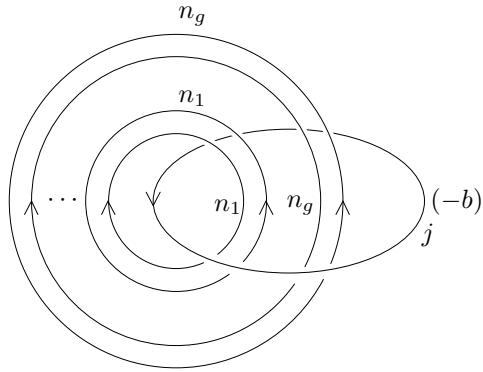
\small

\begin{center}
\begin{texdraw}
\drawdim{cm}

\setunitscale 0.6

\linewd 0.02 \setgray 0

\move(0 4.5)

\move(-1.5 0) \lellip rx:3 ry:1.6

\move(-4.5 -0.1) \lvec(-4.65 0.2)
\move(-4.5 -0.1) \lvec(-4.35 0.2)

\linewd 0.2 \setgray 1

\move(-4 0) \larc r:1.5 sd:308 ed:150
\move(-4 0) \larc r:2 sd:318 ed:150
\move(-4 0) \larc r:3.7 sd:338 ed:155
\move(-4 0) \larc r:3.2 sd:333 ed:155

\linewd 0.02 \setgray 0

\move(-4 0) \larc r:1.5 sd:308 ed:294
\move(-4 0) \larc r:2 sd:318 ed:306
\move(-4 0) \larc r:3.7 sd:338 ed:332
\move(-4 0) \larc r:3.2 sd:333 ed:325

\move(-3.2 -0.3) \htext{$n_{1}$}

\move(-1.6 -0.3) \htext{$n_{g}$}

\move(-2 0) \lvec(-2.15 -0.3)
\move(-2 0) \lvec(-1.85 -0.3)
\move(-0.3 0) \lvec(-0.45 -0.3)
\move(-0.3 0) \lvec(-0.15 -0.3)

\move(-5.5 0) \lvec(-5.65 -0.3)
\move(-5.5 0) \lvec(-5.35 -0.3)
\move(-7.2 0) \lvec(-7.35 -0.3)
\move(-7.2 0) \lvec(-7.05 -0.3)

\move(-6.9 -0.1) \htext{$\cdots$}

\move(-4 3.9) \htext{$n_{g}$}

\move(-4 2.1) \htext{$n_{1}$}

\move(1.6 -0.3) \htext{$(-b)$}
\move(1.4 -1) \htext{$j$}

\move(0 -4)

\end{texdraw}
\end{center}

\caption{The colored ribbon tangle $\Gamma(j,n_{1},\ldots,n_{g})$}\label{fig-A10}
\end{figure}

Let us show that $\sigma(L)=\sigma_{\os}$. To this end let us use the notation
\begin{equation}\label{eq:undermatrix}
A(x_{1},x_{2},\ldots,x_{k}) = \left( \begin{array}{ccccc}
  x_{1}   & 1      & 0      & \cdots  & 0      \\
  1       & x_{2}  & 1      & \cdots  & 0      \\
  0       & 1      & x_{3}  & \cdots  & 0      \\
  \cdots  & \cdots & \cdots & \cdots  & \cdots \\
  0       & 0      & 0      & \cdots  & x_{k}
\end{array}
\right),
\end{equation}
so $A(x_{1},x_{2},\ldots,x_{k})_{ij}$ is $x_{i}$ if $i=j$,
$1$ if $|i-j|=1$, and $0$ elsewhere, $i,j \in \{1,2,\ldots,k\}$.
The linking matrix of $L$ is given by $\left( \begin{array}{cc}
                0 & 0 \\
		0 & A
		\end{array}
\right)$,
where the zeroes refer to the first $2g$ rows and columns and
\begin{equation}\label{eq:matrix}
A=\left( \begin{array}{ccccc}
 -b          & e_{1}  & e_{1}   & \cdots  & e_{1}  \\
 e_{1}^{t} & A_{1}  & 0       & \cdots  & 0      \\
 e_{1}^{t} & 0      & A_{2} & \cdots  & 0      \\
 \cdots      & \cdots & \cdots  & \cdots  & \cdots \\
 e_{1}^{t} & 0      & 0       & \cdots  & A_{n}
\end{array}
\right).
\end{equation}
Here $A_{j} = A(a_{m_{j}}^{(j)},a_{m_{j}-1}^{(j)},\ldots,a_{1}^{(j)})$,
i.e.\ the linking matrix of the $j$'th chain,
and $e_{1}=(1,0,\cdots,0)$. We write $w^{t}$ for a vector $w$ considered as a
column vector. We will calculate the signature of $A$ by reducing
$A$ using combined row and column operations. A main problem is
to avoid dividing by zero. Let us consider $A_{1}$. Write
$m=m_{1}$, $a_{i}=a_{i}^{(1)}$, $p_{i}=\alpha_{i}^{\mC_{1}}$, and
$q_{i}=\beta_{i}^{\mC_{1}}$ to shorten notation.
Assume first that $p_{i} \neq 0$ for all $i=1,2,\ldots,m$ (or
equivalently that $q_{i} \neq 0$ for all $i=1,2,\ldots,m$ since
$q_{1}=1$, $p_{m} = \pm \alpha_{1} \neq 0$, and $q_{i}=p_{i-1}$,
$i=2,3,\ldots,m$). Since $(a_{1},a_{2},\ldots,a_{i})$ is
a continued fraction expansion of $p_{i}/q_{i}$, $i=1,2,\ldots,m$,
we can reduce $A$ to
\begin{equation}\label{eq:A'}
A'=\left( \begin{array}{ccccc}
 -b-\frac{q_{m}}{p_{m}} & 0      & e_{1}   & \cdots  & e_{1}    \\
 0              & A_{1}' & 0       & \cdots  & 0        \\
 e_{1}^{t}      & 0      & A_{2}   & \cdots  & 0        \\
 \cdots         & \cdots & \cdots  & \cdots  & \cdots   \\
 e_{1}^{t}      & 0      & 0       & \cdots  & A_{n}
\end{array}
\right),
\end{equation}
where $A_{1}'=\diag(p_{m}/q_{m},p_{m-1}/q_{m-1},\ldots,p_{1}/q_{1})$.

Next assume that $p_{i} = 0$ for at least one $i \in \{1,2,\ldots,m\}$.
Let $k$ be the smallest element in $\{1,2,\ldots,m\}$ such that
$p_{k}=0$. Choose a non-negative integer $l$ such that  
$a_{k+1}=a_{k+2}= \cdots=a_{k+l}=0$ and $a_{k+l+1} \neq 0$ or $k+l=m$.
Let us first consider the case $k+l < m$ and let $a=a_{k+l+1}$.
If $k>1$ we reduce $A$ to a matrix $A'$ which is equal to $A$, except
that $A_{1}$ is changed to
$A_{1}' =\left( \begin{array}{cc}
                C & 0 \\
		0 & D
		\end{array}
\right)$, where
$D=\diag(p_{k-1}/q_{k-1},p_{k-2}/q_{k-2},\ldots,p_{1}/q_{1})$
and $C=A(a_{m},a_{m-1},\ldots,a_{k+l+2},a,0,\ldots,0)$.
If $k=1$, let $A'=A$ and $A_{1}'=A_{1}$.
Next reduce $A'$ to a matrix $A''$ equal to $A'$, except that
$A_{1}'$ is changed to
\begin{equation}\label{eq:A''}
A_{1}''=\left( \begin{array}{ccccc}
 E       & 0      & \cdots   & 0       & 0      \\
 0       & G      & \cdots   & 0       & 0      \\
 \cdots  & \cdots & \cdots   & \cdots  & \cdots \\
 0       & 0      & \cdots   & G       & 0      \\
 0       & 0      & \cdots   & 0       & D
\end{array}
\right),
\end{equation}
where $G=\diag(2,-1/2)$, and where
$E=A(a_{m},a_{m-1},\ldots,a_{k+l+2},a,0)$ if $l$ is even
and $E=A(a_{m},a_{m-1},\ldots,a_{k+l+2},a)$ if $l$ is odd.
(The row and column with $D$ is not present if $k=1$.
Note that $A''=A'$ if $l=0$.)
Assume that $k+l+1=m$. Then, if $l$ is even, we reduce $A''$ futher to
a matrix equal to the right-hand side of (\ref{eq:A'}) with $A_{1}'$
replaced by a matrix $A_{1}'''$ equal to $A_{1}''$ with
$E$ replaced by $\diag(a,-1/a)$. If $l$ is odd we let $A_{1}'''=A_{1}''$.
Since $p_{k}=0$ we have that $B_{k}^{\mC_{1}}=
\pm \left( \begin{array}{cc}
                0 & -1 \\
		1 & d
		\end{array}
\right)= \pm \Xi \Theta^{d}$
for a $d \in \Z$, so $B_{k+i}^{\mC_{1}}=\pm \Xi^{i+1} \Theta^{d}$,
$i=1,2,\ldots,l$. Since $\Xi^{2}=-1$ we therefore have $q_{k+i}=0$ for $i$ odd
and $p_{k+i}=0$ for $i$ even, $i \in\{0,1,\ldots,l\}$. In particular
$q_{m}=p_{k+l}=0$ for $l$ even. For $l$ odd we have $p_{m}/q_{m}=a$.
From this we also see that the signature of $A_{1}'''$ is equal to
$\sum_{j=1}^{m} \sign(p_{j}q_{j})$.
If $k+l+1<m$ we continue the diagonalization by reducing $E$ in the same manner
as we have reduced $A_{1}$ above. If $l$ is odd, 
$p_{k+l+1}/q_{k+l+1}=a$ and the lower right block in $A_{1}''$,
i.e.\ $\diag(G,\ldots,G,D)$, has signature $\sum_{j=1}^{k+l} \sign(p_{j}q_{j})$.
If $l$ is even, $p_{k+l}=0$ and therefore $q_{k+l+1}=0$ and
$p_{k+l+2}/q_{k+l+2}=a_{k+l+2}$. In this case we begin by reducing 
$A''$ to a matrix equal to $A''$ with $E$ replaced by
$\left( \begin{array}{cc}
                E' & 0 \\
		0 & F
		\end{array}
\right)$ in $A_{1}''$, where $E'=A(a_{m},a_{m-1},\ldots,a_{k+l+2})$ and
$F=\diag(a,-1/a)$. Note that the lower right block in the reduced
$A_{1}''$, i.e.\ $\diag(F,G,\ldots,G,D)$, has
signature $\sum_{j=1}^{k+l+1} \sign(p_{j}q_{j})$. 

The only case left to consider is when $k+l=m$. 
In this case 
$B^{\mC_{1}}=B_{m}^{\mC_{1}}=\Xi^{l}B_{k}^{\mC_{1}}=\pm \Xi^{l+1} \Theta^{d}$,
so $l$ is odd since $p_{m} =\pm \alpha_{1} \neq 0$. 
But then $\beta_{1}=\pm q_{m}=0$, so this case is
only relevant in case of non-normalized Seifert invariants.
For letting the above calculation also work in this case,
let us assume for the moment that $k+l=m$ and $l$
is odd. 
We then reduce
$A$ to a matrix $H$ equal to $A$, except that $A_{1}$ is replaced
by a matrix $H_{1}$ equal to the right-hand side of (\ref{eq:A''})
with $E$ replaced by $J=A(0,0)$.
Finally, we reduce $H$ to a matrix equal to the right-hand side of
(\ref{eq:A'}) with $A_{1}'$ replaced by a matrix $H_{2}$
equal to $H_{1}$ with $J$ replaced by $G$.
Note that the signature of $H_{2}$ is equal to $\sum_{j=1}^{m} \sign(p_{j}q_{j})$.

This ends the reduction involving $A_{1}$. 
We can now continue as above reducing the parts in $A$ involving $A_{j}$,
$j=2,3,\ldots,n$, and get the result.

The RT--invariant $\tau$ of the Seifert manifold
with non-normalized Seifert invariants
$\{\os ;g;(\alpha_{1},\beta_{1}),\ldots,(\alpha_{n},\beta_{n})\}$
is calculated as above by letting $b$ be equal to zero
everywhere, since this manifold
has a surgery presentation as in Fig.~\ref{fig-A1} with $-b$ changed to $0$.

\begin{figure}[ht!]\small

\begin{center}
\begin{texdraw}
\drawdim{cm}

\setunitscale 0.6

\linewd 0.02 \setgray 0

\move(0 5.5)

\move(0 2.6) \lvec(0 -2.6)

\move(0 0.15) \lvec(-0.15 -0.15)
\move(0 0.15) \lvec(0.15 -0.15)

\move(2 2.6) \larc r:2 sd:0 ed:180
\move(2 -2.6) \larc r:2 sd:180 ed:0

\move(3.2 2.6) \lvec(4.8 2.6) \lvec(4.8 1.4) \lvec(3.2 1.4) \lvec(3.2 2.6)

\move(4 1) \lvec(4 1.4)

\move(4 -1.4) \lvec(4 -1)

\move(3.2 -2.6) \lvec(4.8 -2.6) \lvec(4.8 -1.4) \lvec(3.2 -1.4) \lvec(3.2 -2.6)

\lpatt(0.1 0.1) \move(4 1) \lvec(4 -1)

\lpatt()

\move(3.5 1.7) \htext{$T_{i_{1}j}$}
\move(3.5 -2.3) \htext{$T_{i_{g}j}$}

\move(0 4.7) \htext{$(-b-2g)$}

\move(3.8 -3.9) \htext{$j$}
 
\move(0 -5)

\end{texdraw}
\end{center}

\caption{Oriented colored surgery presentation of $(\ns ;g \, | \, b)$}\label{fig-A7}
\end{figure}

Next let us calculate $\tau(M)$ for
$M=(\ns;g \, | \, b;(\alpha_{1},\beta_{1}),\ldots,(\alpha_{n},\beta_{n}))$.
To obtain a surgery presentation of $M$ with only
integral surgery coefficients we make two left-handed twists about every
component with surgery coefficient $1/2$ in the surgery presentation in
Fig.~\ref{fig-A2}. The
components with surgery coefficients $\alpha_{j}/\beta_{j}$
are replaced by chains 
according to the continued fraction expansions
$\mC_{j}$, $j=1,\ldots,n$, as before.
Fig.~\ref{fig-A7} shows a colored oriented version
of the new surgery diagram in the case where there are no exceptional fibers.
The
coupons $T_{i_{l}j}$ represent colored ribbon tangles shown in 
Fig.~\ref{fig-A8}. 
The resulting framed link $L$ has $m=g+1+\sum_{j=1}^{n} m_{j}$ components.
By (\ref{eq:A4}) and the identities in Fig.~\ref{fig-ST} we have
\begin{eqnarray*}
\tau(M) &=& \Delta^{\sigma(L)} \mD^{-\sigma(L)-m-1} \sum_{j \in I} \dim(j)^{1-n} \left( \prod_{i=1}^{n} (SG^{\mC_{i}})_{j,0} \right) \\
 & & \hspace{.7in} \times \sum_{i_{1},\ldots,i_{g} \in I} \left( \prod_{l=1}^{g} \dim(i_{l}) \right) F(R(j,i_{1},\ldots,i_{g})), 
\end{eqnarray*}
where $R(j,i_{1},\ldots,i_{g})$ is the colored ribbon
graph in Fig.~\ref{fig-A7}.
The expression (\ref{eq:nsurgery}) now follows 
by the fact that $\sigma(L)=\sigma_{\ns}$, see below, and by
\reflem{vector} together with the identity
\begin{displaymath}
\sum_{i_{1},\ldots,i_{g} \in I} \left( \prod_{l=1}^{g} \dim(i_{l}) \right) F(R(j,i_{1},\ldots,i_{g})) = v_{j}^{-b-2g} \dim(j) \kappa(j)^{g},
\end{displaymath}
which follows by combining Figures \ref{fig-ST} and \ref{fig-A8}.

Let us show that $\sigma(L)=\sigma_{\ns}$.
The linking matrix of $L$ is given by
\begin{displaymath}
A=\left( \begin{array}{cccccc}
 0      & w^{t}        & 0      & 0       & 0       & 0        \\
 w      & -b-2g        & e_{1}  & e_{1}   & \cdots  & e_{1}    \\
 0      & e_{1}^{t}  & A_{1}  & 0       & \cdots  & 0        \\
 0      & e_{1}^{t}  & 0      & A_{2} & \cdots  & 0        \\
 \cdots & \cdots       & \cdots & \cdots  & \cdots  & \cdots   \\
 0      & e_{1}^{t}  & 0      & 0       & \cdots  & A_{n}
\end{array}
\right),
\end{displaymath}
where $w=(-2,-2,\ldots,-2)$ is a vector of length $g$
and $A_{j}$ is given as in the case of oriented base, i.e.\
$A_{j} = A(a_{m_{j}}^{(j)},a_{m_{j}-1}^{(j)},\ldots,a_{1}^{(j)})$, see (\ref{eq:undermatrix}).
By doing the same combined row and column operations as in the case of oriented base
we reduce $A$ to
\begin{displaymath}
A'=\left( \begin{array}{ccccc}
 D      & 0      & 0        & \cdots  & 0      \\
 0      & A_{1}' & 0        & \cdots  & 0      \\
 0      & 0      & A_{2}' & \cdots  & 0      \\
 \cdots & \cdots & \cdots   & \cdots  & \cdots \\
 0      & 0      & 0        & \cdots  & A_{n}'
\end{array}
\right),
\end{displaymath}
where $A_{j}'$ is a diagonal matrix with signature
$\sum_{k=1}^{m_{j}} \sign(\alpha_{k}^{\mC_{j}}\beta_{k}^{\mC_{j}})$,
$j=1,\ldots,n$, and $D=\left( \begin{array}{cc}
 0 & w^{t} \\
 w & e-2g
\end{array}
\right),$
where, as usual,
$e=-\left( b+\sum_{j=1}^{n} \beta_{j}/\alpha_{j} \right)$
is the Seifert Euler number. By this and the fact that the signature of $D$ is zero 
it follows that $\sigma(L)$ is equal to $\sigma_{\ns}$
in (\ref{eq:nsignature0}).

The RT--invariant $\tau$ of the Seifert manifold
with non-normalized Seifert invariants
$\{\ns ;g;(\alpha_{1},\beta_{1}),\ldots,(\alpha_{n},\beta_{n})\}$
is calculated as above by letting $b$ be equal to zero everywhere,
since this manifold
has a surgery presentation as in Fig.~\ref{fig-A2} with $-b$ changed to $0$.
\end{proof}

The following lemma was first proved by the author
in the case, where the modular categories are induced
by the quantum groups associated to $\frsl_{2}(\C)$.
This was done by a rather long $R$--matrix calculation.
After having presented the result to V. Turaev, he found
the proof below using a geometric computation
which works for an arbitrary modular
category.

\begin{lem}\label{vector}
For all $j \in I$,
\begin{displaymath}
\kappa(j)=\sum_{i \in I} \dim(i) y(i,j) =  \vep_{j} \mD^{2} v_{j}^{2} \delta_{j,j^{*}} \left( \dim(j) \right)^{-1}.
\end{displaymath}
\end{lem}

\begin{figure}[ht!]
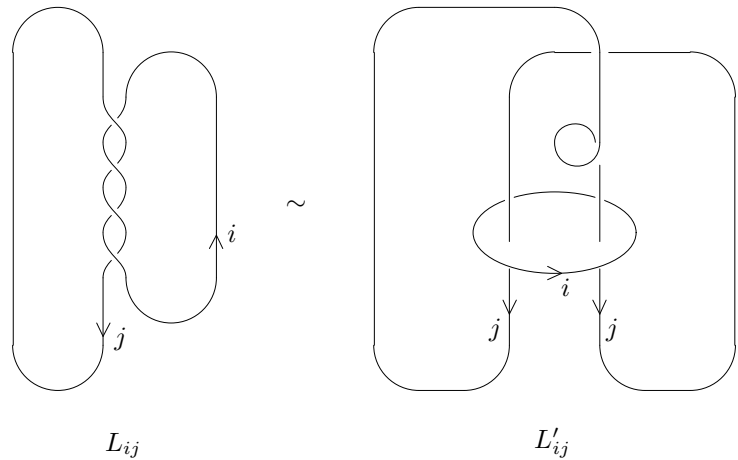
\small

\begin{center}
\begin{texdraw}
\drawdim{cm}

\setunitscale 0.6

\linewd 0.02 \setgray 0

\move(-13 5.5)

\move(4 -6)

\move(-10 4) \lvec(-10 3)

\move(-9.5 2) \clvec(-9.5 1.5)(-10 1.4)(-10 1)

\move(-9.5 0) \clvec(-9.5 -0.5)(-10 -0.6)(-10 -1)

\move(-9.5 3) \clvec(-9.5 2.5)(-10 2.4)(-10 2)

\move(-9.5 1) \clvec(-9.5 0.5)(-10 0.4)(-10 0)

\move(-10 -1) \lvec(-10 -2.5)

\move(-10 -2.3) \lvec(-10.15 -2)
\move(-10 -2.3) \lvec(-9.85 -2)

\linewd 0.2 \setgray 1

\move(-10 2.9) \clvec(-10 2.5)(-9.5 2.4)(-9.5 2.1)

\move(-10 0.9) \clvec(-10 0.5)(-9.5 0.4)(-9.5 0.1)

\move(-10 1.9) \clvec(-10 1.5)(-9.5 1.4)(-9.5 1.1)

\move(-10 -0.1) \clvec(-10 -0.5)(-9.5 -0.6)(-9.5 -0.9)

\linewd 0.02 \setgray 0

\move(-10 3) \clvec(-10 2.5)(-9.5 2.4)(-9.5 2)

\move(-10 1) \clvec(-10 0.5)(-9.5 0.4)(-9.5 0)

\move(-10 2) \clvec(-10 1.5)(-9.5 1.4)(-9.5 1)

\move(-10 0) \clvec(-10 -0.5)(-9.5 -0.6)(-9.5 -1)

\move(-8.5 3) \larc r:1 sd:0 ed:180
\move(-8.5 -1) \larc r:1 sd:180 ed:0

\move(-11 4) \larc r:1 sd:0 ed:180
\move(-11 -2.5) \larc r:1 sd:180 ed:0

\move(-12 4) \lvec(-12 -2.5)

\move(-7.5 3) \lvec(-7.5 -1)

\move(-7.5 -0.05) \lvec(-7.65 -0.35)
\move(-7.5 -0.05) \lvec(-7.35 -0.35)

\move(-9.8 -2.6) \htext{$j$}
\move(-7.3 -0.2) \htext{$i$}

\move(-10 -5) \htext{$L_{ij}$}

\linewd 0.02 \setgray 0

\move(0 3.5)

\move(0 0) \lellip rx:1.8 ry:0.9
\move(0.15 -0.9) \lvec(-0.15 -1.05)
\move(0.15 -0.9) \lvec(-0.15 -0.75)
\move(0.1 -1.4) \htext{$i$}

\linewd 0.2 \setgray 1

\move(-1 -0.2) \lvec(-1 3)
\move(1 -0.2) \lvec(1 4)

\linewd 0.02 \setgray 0

\move(-1 -0.2) \lvec(-1 3)
\move(1 -0.2) \lvec(1 1.5)

\move(1 2) \clvec(1 1.3)(0 1.3)(0 2)
\clvec(0 2.5)(0.85 2.6)(0.9 2)

\move(1 2) \lvec(1 4)

\move(-1 -0.8) \lvec(-1 -2.5)
\move(1 -0.8) \lvec(1 -2.5)

\move(2 -2.5) \larc r:1 sd:180 ed:270

\move(3 -2.5) \larc r:1 sd:270 ed:0

\move(2 -3.5) \lvec(3 -3.5)

\move(-3 -2.5) \larc r:1 sd:180 ed:270

\move(-2 -2.5) \larc r:1 sd:270 ed:0

\move(-3 -3.5) \lvec(-2 -3.5)

\move(4 -2.5) \lvec(4 3)

\move(3 3) \larc r:1 sd:0 ed:90

\move(-4 -2.5) \lvec(-4 4)

\move(0 4) \larc r:1 sd:0 ed:90

\move(-3 4) \larc r:1 sd:90 ed:180

\move(3 4) \lvec(1.2 4)

\move(0 3) \larc r:1 sd:90 ed:180

\move(-3 5) \lvec(0 5)

\move(0 4) \lvec(0.8 4)

\move(-1 -1.8) \lvec(-1.15 -1.5)

\move(-1 -1.8) \lvec(-0.85 -1.5)

\move(1 -1.8) \lvec(1.15 -1.5)

\move(1 -1.8) \lvec(0.85 -1.5)

\move(-1.5 -2.4) \htext{$j$}
\move(1.1 -2.4) \htext{$j$}

\move(-0.5 -5) \htext{$L_{ij}'$}

\move(-6 0.5) \htext{$\sim$}

\end{texdraw}
\end{center}
\vglue -20pt
\caption{A fundamental isotopy}\label{fig-isotopy}
\end{figure}

\begin{proof}
Let $L_{ij}$ be the closure of the ribbon tangle $T_{ij}$.
We have the isotopy shown in Fig.~\ref{fig-isotopy}.
Let $\omega_{j} \co V_{j} \ra \left( V_{j^{*}} \right)^{*}$
be an isomorphism and use this and its inverse to reverse
the orientation of one
of the two strings passing through the component with color
$i$ in $L_{ij}'$ (the link in the right-hand side of
Fig.~\ref{fig-isotopy}).
This enables us to use the identity in Fig.~\ref{fig-A9}, which
gives us the identity in Fig.~\ref{fig-kappaproof}.
The result now follows by applying \reflem{selfdual} together with
$F(L_{ij})=\tr(F(T_{ij}))=y(i,j)\dim(j)$.
\end{proof}

\begin{figure}[ht!]\small

\begin{center}
\begin{texdraw}
\drawdim{cm}

\setunitscale 0.6

\linewd 0.02 \setgray 0

\move(0 3)

\move(0 -3)

\move(0 0.6) \larc r:2 sd:0 ed:180

\move(0 -0.6) \larc r:2 sd:180 ed:0

\move(-3 0.6) \lvec(-1 0.6) \lvec(-1 -0.6) \lvec(-3 -0.6) \lvec(-3 0.6)

\move(3 0.6) \lvec(1 0.6) \lvec(1 -0.6) \lvec(3 -0.6) \lvec(3 0.6)

\move(1.6 -0.35) \htext{$\omega_{j}^{-1}$}

\move(-2.2 -0.3) \htext{$\omega_{j}$}

\move(2 0.6) \lvec(1.85 0.9)

\move(2 0.6) \lvec(2.15 0.9)

\move(2 -0.6) \lvec(1.85 -0.9)

\move(2 -0.6) \lvec(2.15 -0.9)

\move(2.2 1) \htext{$j$}

\move(2.2 -1.3) \htext{$j$}

\move(-16.2 -0.3) \htext{$\sum_{i \in I} \dim(i) F(L_{ij}) \hspace{.1in} \doteq \hspace{.1in} v_{j}^{2} \left( \dim(j) \right)^{-1} \mD^{2} \delta_{j,j^{*}}$}

\end{texdraw}
\end{center}
\vglue -20pt
\nocolon
\caption{}\label{fig-kappaproof}
\end{figure}

It follows that the lemma is also true in case we don't have
a rank if one replaces $\mD^{2}$ with 
$\sum_{u \in I} \left( \dim(u) \right)^{2}$.
The result in \reflem{vector} is independent of how we direct
the component with color $i$ in $T_{ij}$. This follows
by the usual argument since we sum over all colors $i$.
If we reverse the direction of the component with color $j$ we
get $\kappa(j^{*})$ instead of $\kappa(j)$, since the operator
invariant $F$ of a colored ribbon graph is unchanged by changing
the direction of an annulus component if one at the same time
changes the color of that component to the dual color.
Observe however that $\kappa(j^{*})=\kappa(j)$ since $j^{**}=j$.

\begin{rem}\label{signatures}
In this remark we give some alternative expressions
for the signatures (\ref{eq:osignature0}) and (\ref{eq:nsignature0}).
Similar formulas have been obtained
in \cite{FreedGompf} and \cite{Jeffrey} for the case $g=0$ (so $\ep=\os$)
in connection with calculations of framing corrections of Witten's
$3$--manifold invariants of lens spaces and other Seifert manifolds
with base $S^{2}$.
To this end we use the Rademacher Phi function $\Phi$, which 
is defined 
on $PSL(2,\Z)=SL(2,\Z)/\{\pm 1 \}$ by
\begin{equation}\label{eq:C5}
\Phi \left( \left[ \begin{array}{cc}
                  p & r \\
                  q & s
                  \end{array}
\right] \right) = \left\{ \begin{array}{ll}
                  \frac{p+s}{q} - 12(\sign (q))\s (s,|q|) & ,q \neq 0 \\
                  \frac{r}{s} & ,q=0.
                  \end{array}
\right.
\end{equation}
Here, for $q > 0$, the Dedekind sum $\s (s,q)$ is given by
\begin{equation}\label{eq:C6}
\s (s,q)= \frac{1}{4q} \sum_{j=1}^{q-1} \cot\frac{\pi j}{q} \cot \frac{\pi s j}{q}
\end{equation}
for $q>1$ and $\s (s,1)=0$, $s \in \Z$.
We refer to \cite{RademacherGrosswald} for a comprehensive description of this
function and also to \cite{KirbyMelvin2} for a detailed account of the presence of
the Rademacher Phi function and the related Dedekind sums in topological
settings. 
By \cite[Formula (2.20)]{Jeffrey} we have
\begin{equation}\label{eq:Jeffrey}
\sum_{l=1}^{m-1} \sign(\alpha_{l}^{\mC}\beta_{l}^{\mC}) = \frac{1}{3} \left( \sum_{l=1}^{m} a_{l} -\Phi(B^{\mC}) \right)
\end{equation}
for any sequence of integers $\mC=(a_{1},\ldots,a_{m})$.
Formula (\ref{eq:osignature0}) can therefore be changed to
\begin{equation}\label{eq:osignature}
\sigma_{\os}=\sign(e) + \sum_{j=1}^{n} \sign(\alpha_{j}\beta_{j}) +\frac{1}{3} \sum_{j=1}^{n} \left( \sum_{l=1}^{m_{j}} a_{l}^{(j)} -\Phi(B^{\mC_{j}}) \right), 
\end{equation}
where the second sum of course can be put equal to $n$
if we work with normalized Seifert
invariants ($\alpha_{j}>\beta_{j}>0$). 
We can choose the $\mC_{j}$ so that $|a_{l}^{j}| \geq 2$
for $l=1,2,\ldots,m_{j}-1$ and $j=1,2,\ldots,n$.
In this case we have that
$\sign(\alpha_{l}^{\mC_{j}}\beta_{l}^{\mC_{j}})=\sign(a_{l}^{(j)})$,
$l=1,2,\ldots,m_{j}-1$ and $j=1,2,\ldots,n$, so
\begin{equation}\label{eq:osignature1}
\sigma_{\os}=\sign(e) + \sum_{j=1}^{n} \sign(\alpha_{j}\beta_{j}) + \sum_{j=1}^{n} \sum_{l=1}^{m_{j}-1} \sign(a_{l}^{(j)}).
\end{equation}
The formula (\ref{eq:osignature1}) generalizes
\cite[Formula (2.7)]{FreedGompf}.
The expressions (\ref{eq:osignature}) and (\ref{eq:osignature1})
also hold for the signature $\sigma_{\ns}$ in (\ref{eq:nsignature0})
if we remove the term $\sign(e)$.
\end{rem}

We end this section by specializing to lens spaces.
Let $p,q$ be coprime integers. The lens space
$L(p,q)$ is given
by surgery on $S^{3}$ along the 
unknot with surgery coefficient $-p/q$. 
(Recall here that $L(p,-q)$
is diffeomorphic to $L(p,q)$ via an orientation reversing
diffeomorphism.) From this surgery description we
can directly calculate the RT--invariants of $L(p,q)$ by
using a continued fraction expansion of $-p/q$ as in
the proof of \refthm{invariants}. We choose instead to
calculate the invariants by identifying $L(p,q)$ with
certain Seifert fibrations, see the proof below.

In the following corollary we include the possibilities
$L(0,1)=S^{1} \times S^{2}$ and $L(1,q)=S^{3}$, $q \in \Z$.
(Of course we immediately get from (\ref{eq:A4}) that
$\tau(S^{3})=\mD^{-1}$ and $\tau(S^{1} \times S^{2})=1$, since
$S^{3}$ and $S^{1} \times S^{2}$ are
given by surgeries on $S^{3}$ along the empty framed link and the
unknot with framing $0$ respectively.)

\begin{cor}\label{lens-spaces}
Let $p,q$ be a pair of coprime integers and
let $(a_{1},\ldots,a_{m-1})$ be
a continued fraction expansion of $-p/q$ if $q \neq 0$.
If $q=0$ we put $m=3$ and $a_{1}=a_{2}=0$. 
Then the RT--invariant $\tau$
of the lens space $L(p,q)$ is
\begin{equation}\label{eq:lens spaces}
\tau(L(p,q))=(\Delta\mD^{-1})^{\sigma_{\os}} \mD^{-m} G^{\mC}_{0,0},
\end{equation}
where $\mC=(a_{1},\ldots,a_{m-1},0)$ and
$$
\sigma_{\os}=\sum_{l=1}^{m-1} \sign(\alpha_{l}^{\mC}\beta_{l}^{\mC})=\frac{1}{3} \left( \sum_{l=1}^{m-1} a_{l} - \Phi(B^{\mC}) \right). 
$$ 
\end{cor}

\begin{proof}
Lens spaces are Seifert manifolds with base $S^{2}$ and zero, one or
two exceptional fibers. In fact, let $M$ be the Seifert manifold with
non-normalized Seifert invariants
$\{\os;0;(\alpha_{1},\beta_{1}),(\alpha_{2},\beta_{2})\}$,
and let $\alpha_{1}',\beta_{2}'$ be integers such that
$\alpha_{2}\beta_{2}'-\beta_{2}\alpha_{2}'=1$. By 
\cite[Theorem 4.4]{JankinsNeumann}, $M$
is isomorphic to $L(p,q)$ as oriented manifold, where
$p=\alpha_{1}\beta_{2}+\alpha_{2}\beta_{1}$ and
$q=\alpha_{1}\beta_{2}'+\alpha_{2}'\beta_{1}$. In particular
$L(p,q)$ is isomorphic (as oriented manifold) to   
$\{\os;0;(|q|,\sign(q)p),(1,0)\}=\{\os;0;(|q|,\sign(q)p)\}$, see \refrem{lens-spaces1} i).
If $q=0$we have $m=3$ and $a_{1}=a_{2}=0$ by assumption, so
by (\ref{eq:Ssquared})
the right-hand side of (\ref{eq:lens spaces})
is equal to $\mD^{-1}$ as it should be.
If $q \neq 0$ we have
\begin{eqnarray*}
\tau(L(p,q)) &=& (\Delta\mD^{-1})^{\sigma_{\os}} \mD^{-m-2} \sum_{j \in I} \dim(j) (SG^{\mC})_{j,0} \\
 &=& (\Delta\mD^{-1})^{\sigma_{\os}} \mD^{-m-2} (S^{2}G^{\mC})_{0,0}
\end{eqnarray*}
by \refthm{invariants}, since 
$\mC$ is a continued fraction
expansion of $q/p$ (also for $p=0$). The formula
(\ref{eq:lens spaces}) then follows by (\ref{eq:Ssquared}).
The formula for $\sigma_{\os}$ follows
from (\ref{eq:osignature0}) and (\ref{eq:osignature}).
\end{proof}

\begin{rem}\label{lens-spaces1}
{\bf i)} The manifold $\{\os;0;(0,\pm 1)\}$
is a Seifert fibration in the extended sense of \cite{JankinsNeumann},
see \refrem{generalizedSeifert}.
\cite[Theorem 4.4]{JankinsNeumann} is valid for
these more general Seifert fibrations.
Note also that it follows
from this theorem, that a given lens space can have several distinct 
Seifert fibered structures.

{\bf ii)} If $\mC$ is chosen so that $|a_{j}| \geq 2$ for
$j=1,\ldots,m-1$ (which is possible if $|p/q| >1$ by
\cite[Lemma 3.1]{Jeffrey}), then $\sigma_{\os}=\sum_{j=1}^{m-1} \sign(a_{j})$ by
(\ref{eq:osignature1}).

{\bf iii)} In the special case $(p,q)=(n,1)$, $|n| \geq 2$,
the above coincides with the result obtained immediately by (\ref{eq:A4}),
cf.\ \cite[p.~81]{Turaev} (by the conventions
used here our $L(n,1)$ is equal to $-L(n,1)$ in \cite{Turaev}).
The manifold $L(b,1)$, $b \in \Z$, is isomorphic
(as oriented manifold) to the Seifert manifold
with (normalized) Seifert invariants
$(\os;0\; | \; b)$.
\end{rem}

\section{A rational surgery formula for the Reshetikhin--Turaev invariant}\label{sec-A-generalized}

In this section, as in the previous, 
$\left( \mV,\{V_{i}\}_{i \in I} \right)$
is a fixed modular category with a fixed rank $\mD$.
We will use notation introduced above \refthm{invariants}.
Moreover, $\Phi$ is the Rademacher function, see (\ref{eq:C5}).

The surgery formula, we are going to derive, concerns rational surgery
along framed links in arbitrary closed oriented $3$--manifolds.
Before giving the result in the general case, let us first consider
rational surgery along links in $S^{3}$. By using the surgery equivalence
described in Fig.~\ref{fig-A4}, the identities in
Fig.~\ref{fig-ST}, and the method used in the proof of
\refthm{invariants} to calculate signatures we obtain:

\begin{thm}\label{surgery-formula3sphere}
Let $L$ be a link in $S^{3}$ with $m$ components and let $M$ be the
$3$--manifold given by surgery on $S^{3}$ along $L$ with
surgery coefficient $p_{i}/q_{i} \in \Q$ attached to the $i$'th
component, $i=1,2,\ldots,m$ {\em (}so we assume
$q_{i} \neq 0$, $i=1,2,\ldots,m$, see the comments to {\em (\ref{eq:b4}))}.
Moreover, let $\Omega$ be a
colored ribbon graph in $M$ {\em (}also identified with a colored
ribbon graph in $S^{3} \sm L${\em )}.
Let $L_{0}$ be $L$ considered as a framed link
with all components given the framing $0$. Finally, let
$\mC_{i}=(a_{1}^{(i)},\ldots,a_{m_{i}}^{(i)})$ be a
continued fraction expansion of $p_{i}/q_{i}$, $i=1,2,\ldots,m$.
Then
\begin{eqnarray*}
\tau(M,\Omega) &=& (\Delta \mD^{-1})^{\sigma+\sum_{i=1}^{m} c_{i}}\mD^{-\sum_{i=1}^{m} m_{i}} \\
 && \hspace{.4in} \times \sum_{\lambda \in \col(L)} \tau(S^{3},\Gamma(L_{0},\lambda) \cup \Omega) \left( \prod_{i=1}^{m} G^{\mC_{i}}_{\lambda(L_{i}),0} \right),
\end{eqnarray*}
where
$c_{i}=\frac{1}{3}\left( \sum_{j=1}^{m_{i}} a_{j}^{(i)} - \Phi(B^{\mC_{i}}) \right)$,
$i=1,\ldots,m$, and $\sigma$ is the signature of the linking matrix of $L$
{\em (}with the surgery coefficients $p_{1}/q_{1},\ldots,p_{m}/q_{m}$ on the
diagonal {\em )}.\HS
\end{thm}

We have used
(\ref{eq:Jeffrey}). Note that 
$\tau(S^{3},\Gamma(L_{0},\lambda) \cup \Omega)=\mD^{-1}F(\Gamma(L_{0},\lambda) \cup \Omega)$.
\refthm{surgery-formula3sphere} is a generalization of the defining
surgery formula (\ref{eq:A4}). This follows by the facts that
if $a \in \Z$, then
$\Phi(\Theta^{a} \Sigma)=a$ and $(T^{a}S)_{j,0}=v_{j}^{a}\dim(j)$.

In the case of surgery on arbitrary closed oriented $3$--manifolds
along framed links  we do not have a preferred
framing as above, i.e.\ we can not identify a framing of a
link component with an integer in a canonical way, see Appendix B.
Here, by a framed link in a closed oriented $3$--manifold $M$, we mean
a pair $(L,Q)$, where
$Q=\amalg_{i=1}^{m} Q_{i} \co \amalg_{i=1}^{m} (B^{2} \times S^{1}) \ra M$
is an embedding (or more precisely an isotopy class of such embeddings)
and $L$ is the image by $Q$ of $\amalg_{i=1}^{m} (0 \times S^{1})$.
For other definitions of framed links in $3$--manifolds and how these
relate to this definition we refer to Appendix B.
To establish a surgery formula as above in this more general setting we will
need the machinery of the TQFT of Reshetikhin and Turaev.
We have to be precise with orientations because the TQFT--calculations
are sensitive to these orientations. We will use the following conventions.

\begin{conv}\label{conventions}
The space $B^{2} \times S^{1}$ is the standard
solid torus in $\R^{3}$ with the orientation induced by the standard right-handed
orientation of $\R^{3}$. Here $S^{1}$
is the standard unit circle in the $xz$--plane with centre $0$ and oriented
counterclockwise, i.e., $e_{3}$ is a positively oriented tangent vector in
the tangent space $T_{e_{1}}S^{1} \subseteq \R^{3}$, $e_{i}$ being
the $i$'th standard unit vector in $\R^{3}$, see Fig.~\ref{fig-torus}.
For a framed link $(L,Q)$ as above we will always assume that
each copy of $B^{2} \times S^{1}$ is this oriented standard solid torus,
and that $Q$ is orientation preserving after giving the
image of $Q$ the orientation induced by that of $M$ (we
can always obtain this by composing some of the $Q_{i}$ by
$g \times \id_{S^{1}}$ if necessarily, where $g \co B^{2} \ra B^{2}$ is an orientation
reversing homeomorphism).
Moreover, we orient $L$ so that $Q_{i}$ restricted to
$S^{1} \times \{0\}$ is orientation preserving for each $i$.
The oriented meridian $\alpha$ and longitude $\beta$,
see Fig.~\ref{fig-torus},
represent a basis (over $\Lambda$) of $H_{1}(\Sigma_{(1;)};\Lambda) = \Lambda \oplus \Lambda$,
$\Lambda=\Z,\R$, $\Sigma_{(1;)}=S^{1} \times S^{1}$.
(For the notation $\Sigma_{(1;)}$, see Sect.~\ref{sec-The-TQFT}.) 
We identify elements of
$H_{1}(\Sigma_{(1;)};\Lambda)$ with $2$--columns via
$x[\alpha]+y[\beta] \longleftrightarrow 
\left( \begin{array}{c}
		x \\
		y
		\end{array}
\right)$. The endomorphisms of $H_{1}(\Sigma_{(1;)};\Lambda)$
are identified with $2$x$2$--matrices with entries in $\Lambda$
acting on the $2$--columns by multiplication on the left.
\end{conv}

\begin{figure}[ht!]\small

\begin{center}
\begin{texdraw}
\drawdim{cm}

\setunitscale 0.7

\linewd 0.02

\move(0 -0.7) \lellip rx:1 ry:1.25
\move(0 -0.7) \lellip rx:2.6 ry:3.1

\move(-2.75 -0.55) \lvec(-2.6 -0.85)
\move(-2.45 -0.55) \lvec(-2.6 -0.85)

\move(-3.2 -0.85) \htext{$\beta$}

\move(2 -0.9) \lvec(1.7 -0.75)
\move(2 -0.9) \lvec(1.7 -1.05)

\move(1.8 -1.35) \htext{$\alpha$}

\move(5 -1.5) \lvec(7 -1.5)
\lvec(6.7 -1.35)

\move(7 -1.5) \lvec(6.7 -1.65)

\move(5 -1.5) \lvec(5 0.5)
\lvec(4.85 0.2)

\move(5 0.5) \lvec(5.15 0.2)

\move(6.7 -2) \htext{$x$}
\move(5.2 0.2) \htext{$z$}

\move(1 -0.6) \clvec(1 -1)(2.6 -1)(2.6 -0.6)
\lpatt(0.067 0.1) \move(1 -0.6) \clvec(1 -0.2)(2.6 -0.2)(2.6 -0.6)

\end{texdraw}
\end{center}

\nocolon
\caption{}\label{fig-torus}
\end{figure}

Let us recall the notion of rational surgery on $M$ along $(L,Q)$.
Therefore, let $U_{i}=Q_{i}(B^{2} \times S^{1})$
and let $l_{i}=Q_{i}(e_{1} \times S^{1})$ oriented so that
$[l_{i}]=[L_{i}]$ in $H_{1}(U_{i};\Z)$ where
$L_{i}=Q_{i}(0 \times S^{1})$. Moreover, let
$\mu_{i}=Q_{i}(\partial B^{2} \times 1)$ oriented so that
$(\partial Q_{i})_{*}([\alpha])=[\mu_{i}]$ in $H_{1}(\partial U_{i};\Z)$,
where $\partial Q_{i}$ is the restriction of $Q_{i}$
to $\partial B^{2} \times S^{1}=\Sigma_{(1;)}$.
Let $(p_{i},q_{i})$ be pairs of coprime integers, 
let $h_{i} \co \partial U_{i} \ra \partial U_{i}$ be homeomorphisms such that
\begin{equation}\label{eq:b4}
(h_{i})_{*}([\mu_{i}])= \pm (p_{i}[\mu_{i}]+q_{i}[l_{i}])
\end{equation}
in $H_{1}(\partial U_{i};\Z)$, let $h$ be the union of the $h_{i}$,
and let $U=\amalg_{i=1}^{m} U_{i}$ be the image of $Q$.
Then the $3$--manifold $M' = (M \sm \interior (U)) \cup_{h} U$
is said to be the result of doing surgery on $M$ along the framed link $(L,Q)$ with
surgery coefficients $\{ p_{i}/q_{i} \}_{i=1}^{m}$. If $q_{i}=0$ so
$p_{i}=\pm 1$ we just write $\infty$ for $p_{i}/q_{i}$. Such surgeries
do not change the manifold 
(up to an orientation preserving homeomorphism).
If, in (\ref{eq:b4}), $p_{i}=0$ and $q_{i}=\pm 1$ for all $i$, i.e.\ all
surgery coefficients are $0$, then we call $M'$
the result of doing surgery on $M$ along the framed link $(L,Q)$.
We equip $M'$ with the unique orientation extending the
orientation in $M \sm \interior (U)$. The above generalizes
ordinary rational surgery along links in $S^{3}$, see
Appendix B.
We call a homeomorphism $h$ satisfying (\ref{eq:b4}) an {\it attaching map} for the
surgery. We can and will always choose an orientation preserving attaching
map. Up to an orientation preserving homeomorphism the result of doing surgery
on $M$ along the framed link $(L,Q)$ with surgery coefficients
$\{ p_{i}/q_{i} \}_{i=1}^{m}$ is well defined, independent of the choices 
of representative $Q$ and attaching map $h$.

For $\lambda \in \col(L)$ we let 
$\Gamma(L,\lambda)=\cup_{i=1}^{m} \Gamma(L_{i},\lambda(L_{i}))$,
where $\Gamma(L_{i},j)$ is
the colored ribbon graph equal to the
directed annulus $Q_{i}(([-1/2,1/2]\times 0) \times S^{1})$ with oriented
core $L_{i}$ and color $V_{j}$, $j \in I$.

\begin{thm}\label{surgery-formula}
Let $\mC_{i}=(a_{1}^{(i)},\ldots,a_{m_{i}}^{(i)}) \in \Z^{m_{i}}$ be
a continued fraction expansion of $p_{i}/q_{i}$, $i=1,\ldots,m$. 
Moreover let $\Omega$ be a
colored ribbon graph in $M'$ {\em (}also identified with a colored
ribbon graph in $M \sm L${\em )}.
Then
\begin{eqnarray*}
\tau(M',\Omega) &=& (\Delta \mD^{-1})^{\mu+\sum_{i=1}^{m} c_{i}}\mD^{-\sum_{i=1}^{m} m_{i}} \\
 & & \hspace{.4in} \times \sum_{\lambda \in \col(L)} \tau(M,\Gamma(L,\lambda) \cup \Omega) \left( \prod_{i=1}^{m} G^{\mC_{i}}_{\lambda(L_{i}),0} \right),
\end{eqnarray*}
where $\mu$ is a sum of signs given by {\em (\ref{eq:e7})} and 
$c_{i}=\frac{1}{3}\left( \sum_{j=1}^{m_{i}} a_{j}^{(i)} - \Phi(B^{\mC_{i}}) \right)$,
$i=1,\ldots,m$.
\end{thm}

This theorem obviously generalizes \refthm{surgery-formula3sphere}.
\refthm{surgery-formula} follows by \reflem{surgery-formula1}
and \reflem{E5} below.
To prove these lemmas we use
the machinery of the $2+1$--dimensional TQFT $(\tau,\mT)$ of Reshetikhin
and Turaev, see \cite[Chap.~II and IV]{Turaev}.

\subsection{The TQFT $(\tau,\mT)$}\label{sec-The-TQFT}

The modular functor $\mT$ for the TQFT $(\tau,\mT)$ is a functor from
parametrized decorated surfaces (see below) to the category
of finitely generated projective $K$--modules. Decorated surfaces will be
denoted d-surfaces in the following.
We begin by recalling the concepts and notation from
\cite{Turaev} needed. 

A {\it decorated type} or just a type is a tuple
$t=(g;(W_{1},\nu_{1}),\ldots,(W_{m},\nu_{m}))$,
where $g$ is a non-negative integer, $W_{1},\ldots,W_{m}$ are objects of $\mV$,
and $\nu_{1},\ldots,\nu_{m} \in \{ \pm 1 \}$. The number $m$ of pairs 
$(W_{j},\nu_{j})$ is allowed to be zero.
For a type $t$ as above we let
\begin{equation}\label{eq:d2}
\Psi_{t}= \bigoplus_{i \in I^{g}} \Hom(\I,\Phi(t;i)),
\end{equation}
where
$\Phi(t;i)=W_{1}^{\nu_{1}} \otimes W_{2}^{\nu_{2}} \otimes \ldots \otimes W_{m}^{\nu_{m}} \otimes \bigotimes _{r=1}^{g} (V_{i_{r}} \otimes V_{i_{r}}^{*})$
for every $i=(i_{1},\ldots,i_{g}) \in I^{g}$.
Here $W^{+1}=W$ and $W^{-1}=W^{*}$. Note that $\Psi_{t}$
is a finitely generated projective $K$--module as a finite direct sum of
such modules, cf.\ \cite[Lemma II.4.2.1]{Turaev}.

A connected {\it d-surface} is a connected closed oriented surface $\Sigma$ of genus $g$
with $m\geq0$ distinguished ordered and oriented arcs $\gamma_{1},\ldots,\gamma_{m}$, such
that $\gamma_{j}$ is marked with a pair $(W_{j},\nu_{j})$, where $W_{j}$ is an
object of $\mV$ and $\nu_{j} \in \{ \pm 1 \}$, $j =1,2,\ldots,m$. The tuple
$t(\Sigma)= (g;(W_{1},\nu_{1}),\ldots,(W_{m},\nu_{m}))$ is called the {\it type} of the
d-surface. A non-connected closed oriented surface is said to be decorated if its
connected components are decorated. A {\it d-homeomorphism of d-surfaces}
is an orientation
preserving homeomorphism of the underlying surfaces preserving the distinguished arcs
together with their orientations, marks, and order (on each component).

\begin{figure}[ht!]
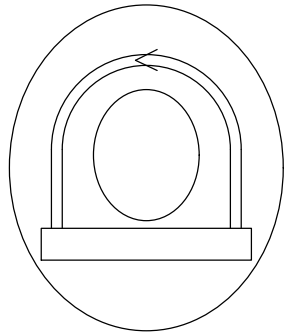
\small

\begin{center}
\begin{texdraw}
\drawdim{cm}

\setunitscale 0.7

\linewd 0.02

\move(-2 -2.6) \lvec(2 -2.6) \lvec(2 -2) \lvec(-2 -2) \lvec(-2 -2.6)  
\move(-1.8 -2) \lvec(-1.8  -0.5) 
\move(-1.6 -2) \lvec(-1.6 -0.5)

\move(1.6 -2) \lvec(1.6 -0.5)
\move(1.8 -2) \lvec(1.8 -0.5)

\move(0 -0.5) \larc r:1.8 sd:0 ed:180
\move(0 -0.5) \larc r:1.6 sd:0 ed:180

\move(-0.2 1.2) \lvec(0.2 1)
\move(-0.2 1.2) \lvec(0.2 1.4)

\move(0 -0.6) \lellip rx:1 ry:1.25
\move(0 -0.85) \lellip rx:2.6 ry:3.1

\end{texdraw}
\end{center}

\caption{Projection of the standard handlebody $U_{(1;)}$}\label{fig-E1}
\end{figure}

For every type $t=(g;(W_{1},\nu_{1}),\ldots,(W_{m},\nu_{m}))$ there is a certain
{\it standard d-surface} of type $t$, denoted $\Sigma_{t}$, which
is the boundary of an oriented handlebody $U_{t}$ of genus $g$ with a certain 
partially colored ribbon graph $R_{t}$ sitting inside, see
\cite[Sect.~IV.1.2]{Turaev}.
In particular, $\Sigma_{(1;)}=S^{1} \times S^{1}$ is an
ordinary oriented torus, see Fig.~\ref{fig-E1}.
The ribbon graph $R_{(1;)}$ lies in the interior of $U_{(1;)}$
and consists of an uncolored coupon with a cap-like uncolored, untwisted, 
and directed band attached to its top base.
A {\it non-connected standard d-surface}
is a disjoint union of a finite number of connected
standard d-surfaces.

In the proof of \refthm{surgery-formula} we will only need
the standard surface $\Sigma_{(1;)}$. For this proof one can therefore
ignore everything about the decoration with destinguished marked arcs.
However, in Sect.~\ref{sec-A-second} we will need such decorations. To avoid
saying things twice, we continue by presenting the concepts
using arbitrary types.

A connected {\it parametrized d-surface} is a connected d-surface $\Sigma$
together with a d-homeo\-morphism $\Sigma_{t} \ra \Sigma$
called the parametrization of $\Sigma$, where
$t=t(\Sigma)$ is the type of $\Sigma$. A {\it non-connected parametrized d-surface} is defined
similarly by using non-connected standard d-surfaces. A morphism in the category of
parametrized d-surfaces, denoted a {\it d-morphism},
is a d-homeomorphism commuting with the parametrizations.

For a connected parametrized d-surface $\Sigma$ of type $t$ we have
$\mT(\Sigma)=\Psi_{t}$. If $\Sigma$ is a non-connected parametrized d-surface with
components $\Sigma_{1},\ldots,\Sigma_{n}$, then $\mT(\Sigma)$ is equal to the non-ordered
tensor product of the $\Psi_{t_{j}}$, $j=1,\ldots,n$, where $t_{j}$ is the type of
$\Sigma_{j}$. Moreover $\mT(\emptyset)=K$.
The modular functor $\mT$ assigns the identity endomorphism to any d-morphism.
By \cite[Lemma IV.1.4.1]{Turaev}, $\mT$ is a modular functor in the sense
of \cite[Sect.~III.1.2]{Turaev}.

A {\it decorated $3$--manifold} is a compact oriented $3$--manifold $M$
with parametrized decorated
boundary $\partial M$ and with an embedded colored ribbon graph $\Omega$,
which is compatible with the
decoration of $\partial M$, see \cite[p.~157]{Turaev}. In this paper
we will only meet decorated $3$--manifolds with empty boundary
or boundary equal to a torus of type $(1;)$. In general, if $\partial M$ contains
no distinguished marked arcs, then $\Omega$ is a colored ribbon graph in
the interior of $M$ with all bases of bands lying on bases of coupons.
A {\it d-homeomorphism of
decorated $3$--manifolds} is an orientation preserving homeomorphism of the underlying
oriented $3$--manifolds preserving all additional structure such as the decoration of
the boundaries and the colored ribbon graphs. Such a d-homeomorphism restricts to a
d-morphism of the boundaries.

A {\it decorated $3$--cobordism} is a triple $(M,\partial_{-}M,\partial_{+}M)$, where
$\partial_{-}M$ and $\partial_{+}M$ (denoted the bottom and top
base respectively) are parametrized d-surfaces and $M$ is a decorated
$3$--manifold with boundary $\partial M=(-\partial_{-}M) \amalg \partial_{+}M$. Here
$-N$ denotes the manifold $N$ with the opposite orientation, where $N$ is an
oriented manifold. (To be precise $-$ is an involution in the space-structure of
parametrized d-surfaces, see \cite[Sect.~IV.1.3 and Sect.~III.1.1]{Turaev}.)
A {\it d-homeomorphism of decorated $3$--cobordisms} is a d-homeomorphism
of the underlying decorated $3$--manifolds which preserves the bases.

A $K$--homomorphism
$\tau(M)=\tau(M,\partial_{-}M,\partial_{+}M) \co \mT(\partial_{-}M) \ra \mT(\partial_{+}M)$
is constructed in \cite[Sect.~IV.1.8]{Turaev} making
$(\tau,\mT)$ a topological quantum field theory (TQFT) based on
decorated $3$--cobordisms and parametrized d-surfaces in the sense of
\cite[Sect.~III.1.4]{Turaev},
cf.\ \cite[Theorem IV.1.9]{Turaev}. If $\partial_{-}M=\emptyset$, then $\tau(M)$
is determined by the element $\tau(M)(1_{K}) \in \mT(\partial_{+}M)$, and it is
common practise in this case to identify $\tau(M)$ with this element.
If $M$ is closed with an embedded colored ribbon graph $\Omega$, then
$\tau(M)$ is the RT--invariant of the pair
$(M,\Omega)$ as defined in (\ref{eq:A4}).
The map $\tau$ is called the {\it operator invariant} of
decorated $3$--cobordisms.

\subsection{Gluing anomalies in the TQFT $(\tau,\mT)$}

The TQFT $(\tau,\mT)$ has so-called (gluing) anomalies,
see \cite[Sect.~III.1.4 and Sect.~IV.4]{Turaev}.
(There is a way to get rid of these anomalies
by changing the TQFT sligthly, cf.\ \cite[Sect.~IV.9]{Turaev}.
However from a computational point of view this `killing' of anomalies
does not make things easier.) To describe these anomalies we need
some concepts from the theory of symplectic vector spaces.

If $H_{1}$ and $H_{2}$ are non-degenerate symplectic
vector spaces, then a {\it Lagrangian relation} between $H_{1}$ and $H_{2}$ is a
Lagrangian subspace of $(-H_{1}) \oplus H_{2}$. For a Lagrangian relation
$N \subseteq (-H_{1}) \oplus H_{2}$ we write $N \co H_{1} \Lra H_{2}$. 
Let $H_{1}$, $H_{2}$
be non-degenerate symplectic vector spaces, and let $\Lambda(H_{i})$
be the set of Lagrangian subspaces of $H_{i}$, $i=1,2$.
A Lagrangian relation $N \co H_{1} \Lra H_{2}$ induces two mappings
$N_{*} \co \Lambda(H_{1}) \ra \Lambda(H_{2})$ and
$N^{*} \co \Lambda(H_{2}) \ra \Lambda(H_{1})$ given by
\begin{displaymath}
N_{*}(\lambda) = \{ h_{2} \in H_{2} \; | \; \exists h_{1} \in \lambda : (h_{1},h_{2}) \in N \}
\end{displaymath}
for $\lambda \in \Lambda(H_{1})$ and
\begin{displaymath}
N^{*}(\lambda) = \{ h_{1} \in H_{1} \; | \; \exists h_{2} \in \lambda : (h_{1},h_{2}) \in N \}
\end{displaymath}
for $\lambda \in \Lambda(H_{2})$.
If $f \co H_{1} \ra H_{2}$ is a symplectic isomorphism and
$\lambda_{i} \in \Lambda(H_{i})$, $i=1,2$, then 
$(N_{f})_{*}(\lambda_{1})=f(\lambda_{1})$
and $(N_{f})^{*}(\lambda_{2})=f^{-1}(\lambda_{2})$, where $N_{f}$ is the graph
of $f$.

For Lagrangian subspaces $\lambda_{1},\lambda_{2},\lambda_{3}$ of a symplectic
vector space $(H,\omega)$, let $W=(\lambda_{1} + \lambda_{2}) \cap \lambda_{3}$
and let $\langle .\, , . \rangle$ be the bilinear form on $W$ defined by
\begin{equation}\label{eq:Maslov index}
\langle a , b \rangle = \omega(a_{2},b)
\end{equation}
for $a,b \in W$ with $a=a_{1}+a_{2}$, $a_{i} \in \lambda_{i}$. This is a
well-defined symmetric form, see e.g.\ \cite[Sect.~IV.3.5]{Turaev}.
The {\it Maslov index} $\mu(\lambda_{1},\lambda_{2},\lambda_{3}) \in \Z$ is the signature
of this bilinear form.
It is invariant under cyclic
permutations of the triple $(\lambda_{1},\lambda_{2},\lambda_{3})$ and changes
sign if we exchange $\lambda_{i}$ and $\lambda_{j}$, $i \neq j$.

If $\Sigma$ is a closed
oriented surface, the real vector space $H_{1}(\Sigma;\R)$ together with the
intersection pairing
\begin{equation}\label{eq:e2}
H_{1}(\Sigma;\R) \times H_{1}(\Sigma;\R) \ra \R
\end{equation}
is a non-degenerate symplectic vector space. For $\Sigma =\emptyset$ we let
$H_{1}(\Sigma;\R)=0$.
For a parametrized 
d-surface $\Sigma$ there is a certain Lagrangian subspace
$\lambda(\Sigma) \subseteq H_{1}(\Sigma;\R)$. For the standard d-surface
$\Sigma_{t}$ of type $t$, $\lambda_{t}=\lambda(\Sigma_{t})$ is the
kernel of the inclusion homomorphism
$H_{1}(\Sigma_{t};\R) \ra H_{1}(U_{t};\R)$.   
For any connected parametrized d-surface $\Sigma$,
$\lambda(\Sigma)=f_{*}(\lambda_{t})$ where $f \co \Sigma_{t} \ra \Sigma$ is the
parametrization. For a non-connected parametrized d-surface $\Sigma$,
$\lambda(\Sigma)$ is the subspace of $H_{1}(\Sigma;\R)$ generated by the
Lagrangian subspaces of the connected components.

For any decorated $3$--cobordism $(M,\partial_{-}M,\partial_{+}M)$ we have
\begin{displaymath}
H_{1}(\partial M;\R)=(-H_{1}(\partial_{-}M;\R)) \oplus H_{1}(\partial_{+}M;\R),
\end{displaymath}
and the kernel of the inclusion homomorphism
$H_{1}(\partial M;\R) \ra H_{1}(M;\R)$
yields a Lagrangian relation 
$H_{1}(\partial_{-}M;\R) \Lra H_{1}(\partial_{+}M;\R)$ which is denoted $N(M)$.
(Note that $N(M)$ does not depend on the
parametrizations and marks of $\partial_{\pm}M$ and the colored ribbon
graph in $M$.)
We let $\lambda_{-}(M)=\lambda(\partial_{-}M)$ and
$\lambda_{+}(M)=\lambda(\partial_{+}M)$.

The anomalies of the TQFT $(\tau,\mT)$ are calculated in
\cite[Theorem IV.4.3]{Turaev}:
Let $M=M_{2}M_{1}$ be a decorated $3$--cobordism obtained from decorated
$3$--cobordisms $M_{1}$ and $M_{2}$ by gluing along a d-morphism
$p \co \partial_{+}(M_{1}) \ra \partial_{-}(M_{2})$.
Set
\begin{displaymath}
N_{r}=N(M_{r}) \co H_{1}(\partial_{-}(M_{r});\R) \Longrightarrow H_{1}(\partial_{+}(M_{r});\R)
\end{displaymath}
for $r=1,2$. Then
\begin{equation}\label{eq:ano}
\tau(M)=(\mD \Delta^{-1})^{m} \tau(M_{2})\tau(M_{1})
\end{equation}
with 
$m=\mu(p_{*}(N_{1})_{*}(\lambda_{-}(M_{1})),\lambda_{-}(M_{2}),N_{2}^{*}(\lambda_{+}(M_{2})))$.
If $\partial_{-}M_{1}=\partial_{+}M_{2}=\emptyset$, then
\begin{equation}\label{eq:ano1}
m=\mu(p_{*}(N(M_{1})),\lambda(-\partial M_{2}), N(M_{2})),
\end{equation}
a Maslov index for Lagrangian subspaces of
$H_{1}(-\partial M_{2};\R)=-H_{1}(\partial M_{2};\R)$.
These anomalies do not depend on the colored
ribbon graphs inside the decorated $3$--cobordisms.
By definition,
$\mT(p) \co \mT(\partial_{+}(M_{1})) \ra \mT(\partial_{-}(M_{2}))$
is the identity and is therefore left out in (\ref{eq:ano}).

\subsection{The projective actions of the modular groups}

We will need to know how the operator invariant $\tau$ of decorated $3$--cobordisms
changes when changing the parametrizations of the parametrized boundary
d-surfaces. To this end we need the projective action of the modular group
$\Mod_{t}$, $t$ a decorated type, cf.\ \cite[Sect.~IV.5]{Turaev}. Here
$\Mod_{t}$ is the group of isotopy classes of d-homeomorphisms
$\Sigma_{t} \ra \Sigma_{t}$. For
$t=(g;)$, $\Mod_{t}=\Mod_{g}$ is the usual modular group of genus $g$.
In this paper we only consider the modular
group $\Mod_{1}$ of genus $1$. The reader can therefore concentrate on this
case if he/she prefers that. We will however state the following
results using arbitrary types since it is not any longer. 
For a decorated type $t$, let $\Sigma=\Sigma_{t}$ and let
$M(\id)=(\Sigma \times [0,1],\Sigma,\Sigma)$, where $\Sigma$ is parametrized by
the identity and $\Sigma \times [0,1]$ is the standard decorated cylinder over
$\Sigma$, cf.\ \cite[p.~158]{Turaev}. For $t=(1;)$ this is just an ordinary
oriented cylinder cobordism without any ribbon graph inside (and without any marked
arcs on the boundary tori, but with identity parametrizations attached to these
tori).
For an arbitrary d-homeomorphism
$g \co \Sigma \ra \Sigma$, let $M(g)$ be as the decorated $3$--cobordism $M(\id)$
except that the bottom base is parametrized by $g$. Let
\begin{displaymath}
\ep(g)=\tau(M(g)) \co \Psi_{t} \ra \Psi_{t}.
\end{displaymath}
Then $g \mapsto \ep(g)$ is a projective linear action of $\Mod_{t}$ on
$\mT(\Sigma_{t})=\Psi_{t}$. In fact we have
\begin{equation}\label{eq:e4}
\ep(gh)=(\mD\Delta^{-1})^{\mu(h_{*}(\lambda_{t}),\lambda_{t},g_{*}^{-1}(\lambda_{t}))} \ep(g)\ep(h),
\end{equation}
cf.\ \cite[Formula (IV.5.1.a)]{Turaev}. 
By the axioms of a TQFT $\ep(\id)=\id$, so
$\ep(g^{-1})=(\ep(g))^{-1}$ by (\ref{eq:e4}).

We use the action $\ep$ to describe the dependency of the operator invariant
$\tau$ on the choice of parametrizations of bases. Let
$(M,\partial_{-}M,\partial_{+}M)$ be a decorated $3$--cobordism with
parametrizations $f_{\pm} \co \Sigma_{t_{\pm}} \ra \partial_{\pm}M$. Let
$g_{\pm} \co \Sigma_{t_{\pm}} \ra \Sigma_{t_{\pm}}$ be d-homeomorphisms.
Provide $\partial_{-}M$ and $\partial_{+}M$ with the structure of
parametrized d-surfaces via $f_{-}'=f_{-}(g_{-})^{-1} \co \Sigma_{t_{-}} \ra \partial_{-}M$
and $f_{+}'=f_{+}(g_{+})^{-1} \co \Sigma_{t_{+}} \ra \partial_{+}M$. Denote
the resulting parametrized d-surfaces by $\partial_{-}'M$ and $\partial_{+}'M$
respectively. These are the same oriented surfaces as $\partial_{-}M$, $\partial_{+}M$
with the same (sets of totally ordered) distinguished marked arcs
but with different parametrizations. 
The cobordism $M$ with the newly parametrized bases is a decorated $3$--cobordism,
say $M'$, between $\partial_{-}(M')=\partial_{-}'M$ and 
$\partial_{+}(M')=\partial_{+}'M$. By definition of
the modular functor $\mT$ we have
$\mT(\partial_{\pm}M)=\mT(\Sigma_{t_{\pm}})=\mT(\partial_{\pm}M')$, and
by \cite[Formula (IV.5.3.b)]{Turaev} we have
\begin{equation}\label{eq:e5}
\tau(M')\ep(g_{-})=(\mD \Delta^{-1})^{\mu_{+}-\mu_{-}} \ep(g_{+}) \tau(M) \co \mT(\Sigma_{t_{-}}) \ra \mT(\Sigma_{t_{+}}),
\end{equation}
where, with $N=N(M)=N(M')$,
\begin{eqnarray*}
\mu_{+} &=& \mu(N_{*}(\lambda_{-}(M)),\lambda_{+}(M),\lambda_{+}(M')), \\
\mu_{-} &=& \mu(\lambda_{-}(M),\lambda_{-}(M'),N^{*}(\lambda_{+}(M'))).
\end{eqnarray*}

\subsection{The proof of \refthm{surgery-formula}}

The standard handlebody $U_{(1;)}$ is a solid torus
with a uncolored ribbon graph $R_{(1;)}$ inside consisting of one
coupon and one band, see Fig.~\ref{fig-E1}. By coloring the band with $V_{i}$
and the coupon with $b_{i}=b_{V_{i}}$, $i \in I$, where 
$b$ is part of the duality of $\mV$,
we get a decorated $3$--manifold
which is the oriented standard solid torus
$B^{2} \times S^{1}$ with a directed untwisted annulus with oriented core
$0 \times S^{1}$ and color $V_{i}$.
Let $Y_{i}$ be this decorated
$3$--manifold considered as a decorated $3$--cobordism
between the empty surface and
$\partial U_{(1;)}=\Sigma_{(1;)}$, where $\Sigma_{(1;)}$
is parametrized by the identity. By \cite[Lemma IV.2.1.3]{Turaev} we have
\begin{equation}\label{eq:handle}
\tau(Y_{i})=b_{i}.
\end{equation}
Let $\Omega$ be a colored ribbon graph in
$S^{3}$ containing an annulus component or a band of color $\I$, and let
$\Omega'$ be the colored ribbon graph obtained from $\Omega$ by eliminating
this annulus (resp.\ band). Then it is a well-known fact that
$F(\Omega')=F(\Omega)$, cf.\ \cite[Exercise I.2.9.2]{Turaev}. Now let $Y$ be the
decorated $3$--cobordism
$(B^{2} \times S^{1},\emptyset,\Sigma_{(1;)})$ with the empty
ribbon graph inside and with $\Sigma_{(1;)}$ parametrized by
the identity. Then
\begin{equation}\label{eq:handle1}
\tau(Y)=\tau(Y_{0})=b_{0}.
\end{equation}
The first equality follows by the just mentioned fact about $F$
together with the technique of presenting $3$--cobordisms by ribbon graphs
in $\R^{3}$, see \cite[Sect.~IV.2]{Turaev} and in particular
\cite[Formula (IV.2.3.a)]{Turaev}. (Alternatively, (\ref{eq:handle1}) follows directly
from the definition of $\tau(Y)$, cf.\ \cite[p.~160]{Turaev}.)

By (\ref{eq:d2}) and the definition of $\mT$,
\begin{displaymath}
\mT(\Sigma_{(1;)})=\Psi_{(1;)} = \bigoplus_{i \in I} \Hom(\I,V_{i} \otimes V_{i}^{*})
\end{displaymath}
is a free $K$--module of rank $\card(I)$ with
basis $\{ b_{i} \co \I \ra V_{i} \otimes V_{i}^{*} \}_{i \in I}$
($\Hom(\I,V_{i} \otimes V_{i}^{*}) \cong \Hom(V_{i},V_{i}) \cong K$
since $V_{i}$ is a simple object).

\begin{lem}\label{surgery-formula1}
Let the situation be as in {\em \refthm{surgery-formula}}. 
Let $g_{i} \co \Sigma_{(1;)} \ra \Sigma_{(1;)}$
be a homeomorphism such that 
$(g_{i})_{*} \co H_{1}(\Sigma_{(1;)},\Z) \ra H_{1}(\Sigma_{(1;)},\Z)$
has 
matrix $B^{\mC_{i}}$ with respect to the basis $\{ [\alpha],[\beta] \}$,
$i=1,\ldots,m$. Then
$$
\tau(M',\Omega)=(\Delta \mD^{-1})^{\mu} \sum_{\lambda \in \col(L)} \tau(M,\Gamma(L,\lambda) \cup \Omega) \left( \prod_{i=1}^{m} A^{(i)}_{\lambda(L_{i}),0} \right),
$$
where $A^{(i)}=\left( A^{(i)}_{k,l}\right)_{k,l \in I}$
is the matrix of $\ep(g_{i}) \co \Psi_{(1;)} \ra \Psi_{(1;)}$ with respect to the
basis $\{ b_{i} \}_{i \in I}$. The integer $\mu$ is given by the sum of Maslov
indices
\begin{equation}\label{eq:e7}
\mu=\sum_{i=1}^{m} \mu((\partial Q_{i})_{*}(\lambda_{(1;)}),(\partial Q_{i} \circ g_{i})_{*}(\lambda_{(1;)}),N(X_{i})),
\end{equation}
where $X_{i}=(M_{i-1} \sm \interior (U_{i}),\partial U_{i}, \emptyset)$. Here $M_{i}$
is the manifold obtained by
doing surgery on
$M$ along $\left( \amalg_{j=1}^{i} L_{j},\amalg_{j=1}^{i} Q_{j} \right)$
with surgery
coefficients $\{ p_{j}/q_{j} \}_{j=1}^{i}$, $i=1,2,\ldots,m$, and
$M_{0}=M$. 
\end{lem}

The $N(X_{i})$ are here subspaces of the
$H_{1}(\partial U_{i};\R)$. The integer $\mu$ in (\ref{eq:e7})
does not depend on the colored
ribbon graph $\Omega$. Moreover $\mu$ is independent of
the choice of the $g_{i}$ since
$(g_{i})_{*}(\lambda_{(1;)})=\Span_{\R} \{ p_{i}[\alpha]+q_{i}[\beta] \}$.

\begin{proof}
Let $h_{i} : \partial U_{i} \ra \partial U_{i}$ be the orientation
preserving homeomorphisms determined by the commutative diagrams

\centerline{
\xymatrix{
  { \Sigma_{(1;)}} \ar[d]_{g_{i}} \ar[r]^{\partial Q_{i}} & { \partial U_{i} } \ar[d]^{h_{i}} \\
  { \Sigma_{(1;)}} \ar[r]^{\partial Q_{i}} & {\partial U_{i}. }
}
}

\noindent The disjoint union $h$ of the $h_{i}$ is an attaching map for the
surgery considered in \refthm{surgery-formula}.
According to the axioms for a TQFT (actually a cobordism theory), see
\cite[Sect.~III.1.3]{Turaev}, we can perform this surgery by
consecutive gluings of the $U_{i}$ to the corresponding boundary components in
$M \sm \interior (U)$ along 
$h_{i} \co \partial U_{i} \ra \partial U_{i} \subseteq M \sm \interior (U)$,
and we see that the general result follows from the case $m=1$.
Therefore, assume $m=1$ and let $g=g_{1}$.
Denote by $X$ the decorated $3$--cobordism $(M \sm \interior (U),\partial U, \emptyset)$,
where $\partial U =-\partial(M \sm \interior (U))$ is parametrized by $\partial Q$. Denote by $X'$ the
decorated $3$--cobordism equal to $X$, except that the base is parametrized by
$\partial Q \circ g$. We identify $U$ with the decorated $3$--cobordism $(U,\emptyset,\partial U)$
with the empty ribbon graph, where $\partial U$ is parametrized by $\partial Q$.
Then $h \co \partial_{+}U \ra \partial_{-}X'$
is a d-morphism of parametrized d-surfaces and
\begin{displaymath}
\tau(M',\Omega)=k^{m_{1}}\tau(X')\tau(U)
\end{displaymath}
by (\ref{eq:ano}), where $k=\mD \Delta^{-1}$ and $m_{1}$
is determined by (\ref{eq:ano1}). By (\ref{eq:e5}) we get
$\tau(X') \ep(g^{-1})=k^{-\mu_{-}} \tau(X)$,
and by the remarks following (\ref{eq:e4}) we have
$\ep(g^{-1})=\ep(g)^{-1}$, so
\begin{displaymath}
\tau(M',\Omega)=k^{m_{1}-\mu_{-}} \tau(X) \ep(g) \tau(U).
\end{displaymath}
Let $Y$ be the decorated $3$--cobordism in (\ref{eq:handle1}).
By \refconv{conventions}, $Q \co Y \ra U$ is a d-homeomorphism, so
$\tau(U)=\tau(Y)$ by the axioms for a TQFT. 
Since $\tau(X)$ is $K$--linear we therefore have (use also (\ref{eq:handle1}))
\begin{displaymath}
\tau(M',\Omega)=k^{m_{1}-\mu_{-}} \tau(X) \ep(g) b_{0}=k^{m_{1}-\mu_{-}} \sum_{j \in I} A_{j,0} \tau(X) b_{j},
\end{displaymath}
where $A=\left( A_{i,j}\right)_{i,j \in I}$ is the matrix of $\ep(g)$ with respect
to the basis $\{ b_{j} \}_{j \in I}$.
The set $\col(L)$ is identified with $I$ since $L$ has only one component.
For $j \in I$ we let $U_{j}=((U,\Gamma(L,j)),\emptyset,\partial U)$
be the decorated
$3$--cobordism identical with $U$, except that $U_{j}$ has the colored ribbon
graph $\Gamma(L,j)$ sitting inside. 
The pair $(M,\Gamma(L,j) \cup \Omega)$ can be obtained
by gluing of $U_{j}$ to $X$ along
$\id_{\partial U} \co \partial_{+}U_{j}=\partial U \ra \partial U=\partial_{-} X$
(which is a d-morphism).
By (\ref{eq:ano}) we therefore get
\begin{displaymath}
\tau(M,\Gamma(L,j) \cup \Omega) = k^{m_{2}} \tau(X) \tau(U_{j}),
\end{displaymath}
where the integer $m_{2}=\mu((N(U_{j})),\lambda(-\partial X),N(X))$ by (\ref{eq:ano1}).
Here $N(U_{j})=N(U)$ since the Lagrangian relation of a decorated $3$--cobordism
does not depend on the colored ribbon graph sitting inside.
By the commutative diagram

\bigskip
\centerline{
\xymatrix{
  { \Sigma_{(1;)} } \ar[d]_{\partial Q} \ar[r]^{i} & { U_{(1;)} } \ar[d]^{Q} \\
  { \partial U } \ar[r]^{j} & { U, }
}
}
\bigskip

\noindent $i$ and $j$ being inclusions, we get 
$N(U)= \ker(j_{*}:H_{1}(\partial U;\R) \ra H_{1}(U;\R)) =\partial Q_{*}(\lambda)$,
where $\lambda=\lambda_{(1;)}$.
Moreover, $\lambda(-\partial X)=\lambda(\partial U)=\partial Q_{*}(\lambda)$,
so $m_{2}=0$.
We have $\tau(U_{j})=\tau(Y_{j})$ 
since $Q \co Y_{j} \ra U_{j}$ is a d-homeomorphism of decorated
$3$--cobordisms by \refconv{conventions}.
By (\ref{eq:handle}) we therefore get
\begin{displaymath}
\tau(M',\Omega)=k^{m_{1}-\mu_{-}} \sum_{j \in I} \tau(M,\Gamma(L,j) \cup \Omega) A_{j,0}.
\end{displaymath}
From (\ref{eq:e5}) we have
$\mu_{-}=\mu(\lambda_{-}(X),\lambda_{-}(X'),N(X))$
since $\lambda_{+}(X')=0$. Here $\lambda_{-}(X)=\lambda(\partial U)=\partial Q_{*}(\lambda)$
and $\lambda_{-}(X')=(\partial Q \circ g)_{*}(\lambda)$, so
\begin{displaymath}
\mu_{-}=\mu(\partial Q_{*}(\lambda),(\partial Q \circ g)_{*}(\lambda),N(X)).
\end{displaymath}
By (\ref{eq:ano1}),
$m_{1}=\mu(h_{*}(N(U)),\lambda(-\partial X'),N(X'))$.
Here $N(U)=\partial Q_{*}(\lambda)$, see above,
so $h_{*}(N(U))=(h \circ \partial Q)_{*}(\lambda)=(\partial Q \circ g)_{*}(\lambda)$.
Moreover, $\lambda(-\partial X')=(\partial Q \circ g)_{*}(\lambda)$,
so $m_{1}=0$.
\end{proof}

To express the matrices $A^{(i)}$
in terms of the $S$-- and $T$--matrices we use the description
of $\ep \co \Psi_{(1;)} \ra \Psi_{(1;)}$ given in
\cite[Sect.~IV.5.4]{Turaev}.
We have an
isomorphism $[f] \mapsto M(f) \co \Mod_{1} \ra SL(2,\Z)$, where $[f]$ is the isotopy
class represented by $f \co \Sigma_{(1;)} \ra \Sigma_{(1;)}$ and $M(f)$ is the
matrix of the induced automorphism on $1$--homologies
$f_{*} \co H_{1}(\Sigma_{(1;)};\Z) \ra H_{1}(\Sigma_{(1;)};\Z)$ 
with respect to the basis $\{[\alpha],[\beta]\}$, see \refconv{conventions}.
Let $[f_{A}]$ be the element in
$\Mod_{1}$ corresponding to the matrix $A \in SL(2,\Z)$
under this isomorphism 
and let $\Theta$ and $\Xi$ be the generators of $SL(2,\Z)$ given in
(\ref{eq:generators}). By \cite[pp.~193-195]{Turaev},
the matrices of the $K$--module automorphisms
$\ep(f_{\Xi}),\ep(f_{\Theta}) \co \Psi_{(1;)} \ra \Psi_{(1;)}$ with respect to the
basis $\{b_{i} \}_{i \in I}$ are given by $\mD S^{-1}$ and $T$ respectively.
(Note here that $\Xi$ and $\Theta$ correspond to respectively $s$ and
$t$ in \cite{Turaev}. Moreover our basis $\{[\alpha],[\beta]\}$ corresponds
to $\{-[\alpha],[\beta]\}$ in \cite[Fig.~IV.5.1]{Turaev}.)

\begin{lem}\label{E5}
Let $\mC=(a_{1},\ldots,a_{n}) \in \Z^{n}$ and let
$g=f_{B^{\mC}}=f_{\Theta}^{a_{n}}f_{\Xi}f_{\Theta}^{a_{n-1}} f_{\Xi} \cdots$ 
$f_{\Theta}^{a_{1}}f_{\Xi}$.
The matrix of $\ep(g):\Psi_{(1;)} \ra \Psi_{(1;)}$ with respect to
the basis $\{b_{i} \}_{i \in I}$ is given by
$$
G = \mD^{-n}(\Delta \mD^{-1})^{m} G^{\mC}S^{-1}\hat{S},
$$
where $\hat{S}=S$ if $n$ is even and $\hat{S}=\bar{S}$ if
$n$ is odd. Here $\bar{S}$ is the $S$--matrix for the mirror
of $\mV$. Moreover,
$m=\sum_{i=1}^{n-1} \sign(\alpha^{\mC}_{i}\beta^{\mC}_{i})=\frac{1}{3}\left( \sum_{i=1}^{n} a_{i} - \Phi(B^{\mC}) \right)$.
\end{lem}

\begin{proof}
Let $h \co \Sigma_{(1;)} \ra \Sigma_{(1;)}$ be
an arbitrary orientation preserving diffeomorphism, and let
$\left( \begin{array}{cc}
                a & b \\
		c & d
		\end{array}
\right) \in SL(2,\Z)$ be the matrix of the induced automorphism
$h_{*} \co H_{1}(\Sigma_{(1;)};\Z) \ra H_{1}(\Sigma_{(1;)};\Z)$ 
with respect to the basis $\{ [\alpha],[\beta] \}$. We have
$\lambda_{1}=\lambda_{(1;)}=\Span_{\R}\{ [\alpha] \}$.
Therefore $(f_{\Theta})_{*}(\lambda_{1})=\lambda_{1}$, and we get directly
from (\ref{eq:e4}) that
\begin{displaymath}
\ep(hf_{\Theta}^{m})=\ep(h)\ep(f_{\Theta}^{m}),\hspace{.2in} \ep(f_{\Theta}^{m}h)=\ep(f_{\Theta}^{m})\ep(h)
\end{displaymath}
and moreover $\ep(f_{\Theta}^{m})=(\ep(f_{\Theta}))^{m}$ for all $m \in \Z$.
Next consider composition with $f_{\Xi}$. We have
that $(f_{\Xi})_{*}(\lambda_{1})=\lambda_{2}$, where $\lambda_{2}=\Span_{\R}\{ [\beta] \}$.
Let $\lambda_{3}=h_{*}(\lambda_{1})=\Span_{\R}\{ a[\alpha]+c[\beta] \}$, and let
$\omega$ be the intersection pairing (\ref{eq:e2}) with $\Sigma=\Sigma_{(1;)}$.
Let $W=(\lambda_{1}+\lambda_{2})\cap\lambda_{3}=\lambda_{3}$
and let $\langle .\, , . \rangle$ be
the bilinear form on $W$ defined
in (\ref{eq:Maslov index}). For $x=a[\alpha]+c[\beta]$ we have
$\langle x,x \rangle=\omega(c[\beta],a[\alpha]+c[\beta])=ac\omega([\beta],[\alpha])$.
By definition, the Maslov index
$\mu(h_{*}(\lambda_{1}),\lambda_{1},(f_{\Xi})_{*}^{-1}(\lambda_{1}))=\mu(\lambda_{1},\lambda_{2},\lambda_{3})$
is equal to the signature of $\langle .\, , . \rangle$
which again is equal to $-\sign(ac)$ since
$\omega([\beta],[\alpha])=-\omega([\alpha],[\beta])=-1$. Therefore
\begin{displaymath}
\ep(f_{\Xi}h)=(\Delta \mD^{-1})^{\sign(ac)}\ep(f_{\Xi})\ep(h).
\end{displaymath}
Let
$g_{i}=f_{\Theta}^{a_{i}}f_{\Xi}f_{\Theta}^{a_{i-1}}\cdots f_{\Theta}^{a_{1}}f_{\Xi}$.
In particular $g=g_{n}$, and
$(g_{i})_{*}:H_{1}(\Sigma_{(1;)};\Z) \ra H_{1}(\Sigma_{(1;)};\Z)$ 
has the matrix $B^{\mC}_{i}$ with respect to the basis $\{ [\alpha],[\beta] \}$.
For $i \geq 1$ we have
\begin{displaymath}
\ep(g_{i+1})=\ep(f_{\Theta}^{a_{i+1}} f_{\Xi}g_{i})=(\Delta \mD^{-1})^{\sign(\alpha^{\mC}_{i}\beta^{\mC}_{i})}(\ep(f_{\Theta}))^{a_{i+1}}\ep(f_{\Xi})\ep(g_{i}).
\end{displaymath}
Also note that
$\ep(g_{1})= (\ep(f_{\Theta}))^{a_{1}} \ep(f_{\Xi})$.
We therefore get
\begin{displaymath}
G=(\Delta \mD^{-1})^{m} T^{a_{n}}(\mD S^{-1})T^{a_{n-1}}(\mD S^{-1}) \cdots T^{a_{1}} (\mD S^{-1}),
\end{displaymath}
where
$m=\sum_{i=1}^{n-1} \sign(\alpha^{\mC}_{i}\beta^{\mC}_{i})$.
We also have $m=\frac{1}{3}\left( \sum_{i=1}^{n} a_{i} - \Phi(B^{\mC}) \right)$
by (\ref{eq:Jeffrey}).
By (\ref{eq:Ssquared}) we have that $S^{-1}=\mD^{-2} \bar{S}$.
Recall here that $\bar{S}_{i,j}=S_{i^{*},j}$. The result now
follows by using that $v_{i^{*}}=v_{i}$ and $S_{i^{*},j^{*}}=S_{i,j}$
for all $i,j \in I$, see \cite[Formulas (II.3.3.a-b)]{Turaev}, and by using
that $i \mapsto i^{*}$ is an involution in $I$.
\end{proof}

\noindent Note that $(G^{\mC}S^{-1}\hat{S})_{j,0}=G^{\mC}_{j,0}$ for all
$j \in I$ since $0^{*}=0$.

\begin{rem}
We have chosen in this paper to work with the generators
$\Xi$ and $\Theta$ for $SL(2,\Z)$ since it seems to
be the standard. However, the above result suggests that
in the above setting it is more natural to work with $\Xi^{-1}=-\Xi$
and $\Theta$. If we do this we will not need the strange factor
$S^{-1}\hat{S}$ in the formula for $G$ in \reflem{E5}.
Note also that the use of $\Xi^{-1}$ instead of $\Xi$ 
causes no difficulties with respect to the
matrices $B^{\mC}_{k}$ in (\ref{eq:Bmatrix}) since we actually only need these
as elements of $PSL(2,\Z)$ in any case. Another (more radical)
way to avoid a factor such as $S^{-1}\hat{S}$ is to use $\bar{S}$
as the $S$--matrix for a modular category instead of $S$, see
\cite{Kirillov}, \cite{BakalovKirillov}.
\end{rem}

\section{A second proof of formula (\ref{eq:osurgery})}\label{sec-A-second}

In this section we use the surgery formula in \refthm{surgery-formula}
to calculate the invariant of
$M=(\os;g \, | \, b;$ $(\alpha_{1},\beta_{1}),\ldots,(\alpha_{n},\beta_{n}))$.
First assume that $b \neq 0$.
Let $\alpha_{n+1}=1$ and $\beta_{n+1}=b$,
let $\Sigma_{g}$ be a closed oriented surface of genus $g$, let
$D_{1},\ldots,D_{n+1}$ be disjoint closed disks in $\Sigma_{g}$,
and let
$Q_{i}' \co D_{i} \times S^{1} \hookrightarrow \Sigma_{g} \times S^{1}$
be the inclusion, $i=1,\ldots,n+1$.
Let $B^{2} \times S^{1}$ be the oriented standard solid torus in $\R^{3}$,
see \refconv{conventions}, and
let $k_{i} \co B^{2} \ra D_{i}$ be orientation preserving homeomorphisms,
$i=1,\ldots,n+1$. Moreover let 
$Q_{i}=Q_{i}' \circ (k_{i}\times id_{S^{1}}) \co B^{2} \times S^{1} \ra \Sigma_{g} \times S^{1}$
and $L_{i}=Q_{i}(0\times S^{1})$.
The manifold $M$ is given by surgery on $\Sigma_{g} \times S^{1}$
along the link $L=\amalg_{i=1}^{n+1} L_{i}$
with framing $Q=\amalg_{i=1}^{n+1}Q_{i}$ and surgery coefficients 
$\{ \alpha_{i}/\beta_{i} \}_{i=1}^{n+1}$. (The orientation of
$\Sigma_{g} \times S^{1}$ is given by the orientation of $\Sigma_{g}$
followed by the orientation of $S^{1}$, where $S^{1}$ is oriented as
in the oriented standard solid torus $B^{2} \times S^{1}$.)
Let $\mC_{i}$ be as above \refthm{invariants}, $i=1,\ldots,n$, let
$\mC_{n+1}=(-b,0)$ (a continued fraction expansion of
$\alpha_{n+1}/\beta_{n+1}$),
and let $m_{n+1}=2$.
By \refthm{surgery-formula} we have
\begin{eqnarray*}
\tau(M) &=& (\Delta \mD^{-1})^{\mu+\sum_{i=1}^{n+1} c_{i}} \mD^{-\sum_{i=1}^{n+1} m_{i}} \\
 & & \hspace{.7in} \times \sum_{\lambda \in \col(L)} \tau(\Sigma_{g} \times S^{1},\Gamma(L,\lambda)) \left( \prod_{i=1}^{n+1} G^{\mC_{i}}_{\lambda(L_{i}),0} \right).
\end{eqnarray*}
Here $(\Sigma_{g} \times S^{1},\Gamma(L,\lambda))=\Sigma_{t} \times S^{1}$
with $t=(g;(V_{\lambda(L_{1})},1),\ldots,(V_{\lambda(L_{n+1})},1))$.
For an arbitrary type $t$ we have
$$
\tau(\Sigma_{t} \times S^{1}) = \Dim(\Psi_{t}),
$$
where $\Psi_{t}$ is the projective $K$--module
given in (\ref{eq:d2}) and
$\Dim$ is the dimension in the (ribbon) category of finitely
generated projective $K$--modules, see \cite[Sect.~I.1.7.1 and Appendix I]{Turaev}.
This follows by \cite[Theorem IV.7.2.1]{Turaev},
the remarks following this theorem,
and \cite[Sect.~IV.6.7]{Turaev}.
The dimension $\Dim(\Psi_{t})$ is calculated for an arbitrary type
in \cite[Sect.~IV.12]{Turaev}. (The formula for this dimension
is a generalization of Verlinde's well-known formula \cite{Verlinde}
to the setting of modular
categories.)
For $t=(g;(V_{i_{1}},1),\ldots,(V_{i_{m}},1))$,
$i_{1},\ldots,i_{m} \in I$, we have
\begin{equation}\label{eq:Verlinde}
\Dim (\Psi_{t})= \mD^{2g-2} \sum_{j \in I} (\dim(j))^{2-2g-m} \left( \prod_{k=1}^{m} S_{i_{k},j} \right)
\end{equation}
by \cite[Theorem IV.12.1.1]{Turaev}.
Putting the above together we get
\begin{eqnarray*}
\tau(M) &=& (\Delta \mD^{-1})^{\mu+\sum_{i=1}^{n+1} c_{i}} \mD^{2g-2-\sum_{i=1}^{n+1} m_{i}} \\
 & & \hspace{.7in} \times \sum_{j \in I} (\dim(j))^{2-2g-n-1} \left( \prod_{i=1}^{n+1} (SG^{\mC_{i}})_{j,0} \right).
\end{eqnarray*}
Here $SG^{\mC_{n+1}}= S^{2}T^{-b}S$, so by (\ref{eq:Ssquared}) we get 
\begin{displaymath}
(SG^{\mC_{n+1}})_{j,0}=\mD^{2} \sum_{k,l \in I} \delta_{j^{*},k}\delta_{k,l}v_{k}^{-b}S_{l,0}=\mD^{2} v_{j}^{-b}\dim(j),
\end{displaymath}
where we use that $v_{j^{*}}=v_{j}$, cf.\ \cite[p.~90]{Turaev}, and that
$\dim(j^{*})=\dim(j)$ by \cite[Corollary I.2.8.2]{Turaev} and
the definition of a modular category. Therefore
\begin{eqnarray*}
\tau(M) &=& (\Delta \mD^{-1})^{\mu+\sum_{i=1}^{n+1} c_{i} } \mD^{2g-2-\sum_{i=1}^{n}m_{i}} \\
 & & \hspace{.7in} \times \sum_{j \in I} v_{j}^{-b}(\dim(j))^{2-2g-n} \left( \prod_{i=1}^{n} (SG^{\mC_{i}})_{j,0} \right).
\end{eqnarray*}
This expression is identical with (\ref{eq:osurgery}) and
(\ref{eq:osignature}) if
\begin{equation}\label{eq:signatures}
\mu+c_{n+1}=n+\sign(e).
\end{equation}
Since $\mC_{n+1}=(-b,0)$ and $\Theta^{-b} \Xi
= \left( \begin{array}{cc}
         b & 1 \\
         -1 & 0
         \end{array}
\right)$ we immediately get $c_{n+1}=-\sign(b)$
by \reflem{E5}. By
notation from \reflem{surgery-formula1} we have $M_{0}=\Sigma_{g} \times S^{1}$ and
\begin{displaymath}
M_{i}=(\os;g \, | \, 0;(\alpha_{1},\beta_{1}),\ldots,(\alpha_{i},\beta_{i})), \hspace{.2in} i=1,2,\ldots,n.
\end{displaymath}
Moreover $X_{i}=M_{i-1} \sm \interior (U_{i})$ is obtained from
$Y_{i}=(\Sigma_{g} \sm \interior (D_{1} \cup \ldots \cup D_{i})) \times S^{1}$ by
pasting in $i-1$ solid tori $U_{1},\ldots,U_{i-1}$ as explained above
leaving one torus shaped cave. We have
\begin{eqnarray*}
\pi_{1}(Y_{i}) &=& < a_{1},b_{1},\ldots,a_{g},b_{g},q_{1},\ldots,q_{i},h \; | \; \prod_{j=1}^{i} q_{j} \prod_{j=1}^{g} [a_{j},b_{j}]=1, \\
 & & \hspace{.2in} [h,a_{k}]=[h,b_{k}]=[h,q_{l}]=1,\; k=1,\ldots,g,\; l=1,\ldots,i >,
\end{eqnarray*}
cf.\ \cite[Sect.~6 pp.~34--35]{JankinsNeumann}, \cite[Sect.~10]{Seifert2}. Here
$q_{j}$ corresponds to the `partial cross-section' $\partial D_{j} \times \{1 \}$,
and $h$ is a fiber. The generators $a_{1},b_{1},\ldots,a_{g},b_{g}$ are induced
by the usual generators of
$\pi_{1}(\Sigma_{g})= < a_{1},b_{1},\ldots,a_{g},b_{g} \; | \; \prod_{j=1}^{g}[a_{j},b_{j}]=1 >$.
By the theorem of Seifert and Van Kampen, gluing in the torus $U_{j}$ adds a
new generator $t$ and two new relations $q_{j}^{\alpha_{j}}h^{\beta_{j}}=1$ and
$q_{j}^{\rho_{j}}h^{\sigma_{j}}=t$. The generator $t$ and the last relation
can be deleted by a Tietze transformation, so we get
\begin{eqnarray*}
\pi_{1}(X_{i}) &=& < a_{1},b_{1},\ldots,a_{g},b_{g},q_{1},\ldots,q_{i},h \; | \; \prod_{j=1}^{i} q_{j} \prod_{j=1}^{g} [a_{j},b_{j}]=1, \\
 & & \hspace{1.0in} [h,a_{k}]=[h,b_{k}]=[h,q_{l}]=q_{s}^{\alpha_{s}}h^{\beta_{s}}=1, \\
 & & \hspace{.85in} k=1,\ldots,g,\; l=1,\ldots,i,\; s=1,\ldots,i-1  >.
\end{eqnarray*}
By abelianizing we see that
$H_{1}(X_{i};\Z)=\Z^{2g} \oplus T$,
where
\begin{displaymath}
T = < q_{1},\ldots,q_{i},h \; | \; \sum_{j=1}^{i}q_{j}=\alpha_{s}q_{s}+\beta_{s}h=0,\;s=1,\ldots,i-1 >,
\end{displaymath}
and by the universal coefficient theorem we have
$H_{1}(X_{i};\R)= \R^{2g} \oplus (T \otimes_{\Z} \R)$.
Let
\begin{eqnarray*}
\lambda_{1}^{(i)} &=& (\partial Q_{i})_{*}(\lambda_{(1;)})=\Span_{\R} \{ q_{i} \}, \\
\lambda_{2}^{(i)} &=& (\partial Q_{i})_{*}(\Span_{\R} \{ \alpha_{i}[\alpha]+\beta_{i}[\beta] \}) = \Span_{\R} \{ \alpha_{i}q_{i}+\beta_{i}h \}, \\
\lambda_{3}^{(i)} &=& \ker \left( i_{*} \co H_{1}(\partial U_{i} ; \R) \ra H_{1}(X_{i};\R) \right),
\end{eqnarray*}
$i=1,2,\ldots,n+1$.
Here $H_{1}(\partial U_{i}; \R)=\Span_{\R}\{q_{i},h \}$.
A small calculation shows that $\lambda_{3}^{(i)}=\Span_{\R}\{ y_{i} \}$,
where $y_{1}=q_{1}$ and
$y_{i}=q_{i}- \left( \sum_{j=1}^{i-1} \frac{\beta_{j}}{\alpha_{j}} \right) h$
for $i=2,\ldots,n+1$. Let $\langle \; , \; \rangle_{i}$ be the bilinear form on 
$(\lambda_{1}^{(i)} + \lambda_{2}^{(i)}) \cap \lambda_{3}^{(i)} = \lambda_{3}^{(i)}$
defined by (\ref{eq:Maslov index}) with $H=H_{1}(\partial U_{i};\R)$,
$\omega$ equal to the intersection pairing $\omega_{i}$ on $H_{1}(\partial U_{i};\R)$,
and $\lambda_{j}=\lambda_{j}^{(i)}$, $j=1,2,3$.
Moreover, let
$\mu_{i}=\mu(\lambda_{1}^{(i)},\lambda_{2}^{(i)},\lambda_{3}^{(i)})$ be the
Maslov index equal to the signature of $\langle .\, , . \rangle_{i}$.
We get immediately that $\mu_{1}=0$. Let $i \in \{2,\ldots,n+1\}$ and
let $x_{i}=\alpha_{i}q_{i}+\beta_{i}h$ and 
$t_{i}=\sum_{j=1}^{i-1} \frac{\beta_{j}}{\alpha_{j}}$. Then
$y_{i}=\left( 1 +\frac{\alpha_{i}}{\beta_{i}} t_{i} \right)q_{i} - \frac{t_{i}}{\beta_{i}} x_{i}$.
Therefore
\begin{displaymath}
\langle y_{i},y_{i} \rangle_{i} = \omega_{i}(-\frac{t_{i}}{\beta_{i}} x_{i} ,q_{i}-t_{i}h)=\frac{\alpha_{i}}{\beta_{i}} t_{i+1} t_{i},
\end{displaymath}
where $t_{n+2}=-e$. 
Here $t_{i}>0$, and for $i \leq n$ we have $\alpha_{i}/\beta_{i}>0$ and $t_{i+1}>0$.
Therefore $\mu_{i}=1$ for $i=2,\ldots,n$. Finally, $\mu_{n+1}=-\sign(b)\sign(e)$,
so $\mu=\sum_{i=1}^{n+1} \mu_{i}= n-1-\sign(b)\sign(e)$. The
identity (\ref{eq:signatures}) is therefore equivalent with the
identity $\sign(e)+\sign(b)+\sign(b)\sign(e)+1=0$ which
is true. Note that the above also holds in case $n=0$
(no exceptional fibers). In this case $e=-b$.
In case $b=0$ we ignore everything concerning the surgery along the
component $L_{n+1}$.
If $n>0$ we have to show that $\mu=n+\sign(e)$, where $\mu=\sum_{i=1}^{n} \mu_{i}$
and $e=-\sum_{i=1}^{n} \frac{\beta_{i}}{\alpha_{i}} <0$. Since $\mu=\sum_{i=1}^{n} \mu_{i}=n-1$
this identity is true. If $n=0$ the surgery formula is of no use.
In this case
$\tau(M)=\tau(\Sigma_{g} \times S^{1})= \mD^{2g-2}\sum_{j \in I} \left( \dim(j) \right)^{2-2g}$
by (\ref{eq:Verlinde}) in accordance with \refthm{invariants}.

\section{A third proof of formula (\ref{eq:osurgery})}\label{sec-A-third}

In this section we will use the formula in
\cite[Theorem X.9.3.1]{Turaev}
for the RT--invariant of graph manifolds
to derive (\ref{eq:osurgery}). This formula is valid for
unimodular categories with a rank,
see Sect.~\ref{sec-Modular-categories}.
 
Let $M = (\os;g \, | \, b;(\alpha_{1},\beta_{1}),\ldots,(\alpha_{n},\beta_{n}))$.
It turns out to be an advantage to work with
$-M$ instead of $M$. According to \refthm{SeifertClassification},
$-M=(\os;g \, | \, -b-n;(\alpha_{1},\alpha_{1}-\beta_{1}),\ldots,(\alpha_{n},\alpha_{n}-\beta_{n}))$.
Let $\mC_{j}=(a_{1}^{(j)},\ldots,a_{m_{j}}^{(j)})$ be a continued
fraction expansion of $\alpha_{j}/\beta_{j}$
with $a_{l}^{j} \geq 2$ and let $A$ be the $m \times m$--matrix in
(\ref{eq:matrix}), $m=1+\sum_{j=1}^{n} m_{j}$. By
\cite[Corollary 5 p.~30]{Orlik}  and \cite[Sect.~X.9.2]{Turaev}, the
$3$--manifold $-M$ is the $3$--dimensional
graph manifold determined by the matrix
$-A$ and the integers $g_{1},g_{2},\ldots,g_{m}$, where $g_{1}=g$ and
$g_{2}=\ldots=g_{m}=0$.

Let $\left( \mV,\{V_{i}\}_{i \in I} \right)$ be a unimodular category with a fixed
rank $\mD$. By \cite[Theorem X.9.3.1]{Turaev},
the RT--invariant of the graph manifold $N$, determined by the
symmetric square matrix $B=(a_{p,q})_{p,q=1}^{m}$ over $\Z$ and a
sequence of non-negative integers $g_{1},g_{2},\ldots,g_{m}$, is given by
\begin{eqnarray}\label{eq:graph}
&&\tau_{(\mV,\mD)}(N)=\Delta^{\sigma(B)} \mD^{b} \\
 && \hspace{.9in} \times \sum_{\varphi \in I^{m}} \left( \prod_{p=1}^{m} v_{\varphi(p)}^{a_{p,p}} (\dim(\varphi(p)))^{2-2g_{p}-a_{p}} \prod_{p<q} (s_{p,q}^{\varphi})^{|a_{p,q}|} \right), \nonumber
\end{eqnarray}
where $b=b_{1}(N)-b_{0}(N)-m-\nul(B)-\sigma(B)$, where $\nul(B)$ and
$\sigma(B)$ are the nullity and signature of $B$ respectively, and
$a_{p} = \sum_{q \neq p} |a_{p,q}|$. Moreover,
$s_{p,q}^{\varphi}=S_{\varphi(p),\varphi(q)}$ if $a_{p,q} \geq 0$ and
$s_{p,q}^{\varphi}=S_{\varphi(p)^{*},\varphi(q)}$ if $a_{p,q} < 0$.

We have $b_{0}(-M)=1$ and $\sigma(-A)=-\sigma_{\os}$,
where $\sigma_{\os}$ is given by (\ref{eq:osignature0}),
(\ref{eq:osignature}). By \cite[Corollary 6.2]{JankinsNeumann},
$b_{1}(-M)=b_{1}(M)=2g+\delta_{e,0}$,
where $\delta_{e,0}=1$
if $e=0$ and $0$ otherwise, and by the proof of
(\ref{eq:osignature0}) we have
$\nul(-A)=\nul(A)=\delta_{e,0}$.
If we write
$-A=\left( a_{p,q} \right)_{p,q=1}^{m}$ then
\begin{displaymath}
a_{p} = \left\{ \begin{array}{ll}
		 n  & ,p=1 \\
		 2  & ,p \in \{1,2,\ldots,m\} \sm \{1,m_{1}+1,m_{1}+m_{2}+1,\ldots,m\} \\
                 1  & ,p \in \{1,m_{1}+1,m_{1}+m_{2}+1,\ldots,m\}.
		\end{array}
\right.
\end{displaymath}
According to \cite[Exercise II.2.5]{Turaev} we have
\begin{displaymath}
\tau_{(\mV,\mD)}(M)=\tau_{(\omV,\mD)}(-M),
\end{displaymath}
where $\omV$ is the mirror
of $\mV$, see Sect.~\ref{sec-Modular-categories}.
The modular category $\left( \omV,\{V_{i}\}_{i \in I} \right)$
is also unimodular, cf.\ \cite[Exercise VI.2.3.1]{Turaev}.
By (\ref{eq:graph}) we get
\begin{eqnarray*}
\tau_{(\omV,\mD)}(-M) &=& (\Delta_{\omV}\mD^{-1})^{-\sigma_{\os}} \mD^{2g-2-\sum_{j=1}^{n} m_{j}} \\
 & & \hspace{.2in} \times \sum_{j,k_{1}^{1},\ldots,k_{m_{1}}^{1},k_{1}^{2},\ldots,k_{m_{n}}^{n} \in I} \bar{v}_{j}^{b} (\dim(j))^{2-2g-n} \\
 & & \hspace{.2in} \times \prod_{j=1}^{n} \prod_{l=1}^{m_{j}-1} \bar{v}_{k_{l}^{j}}^{-a_{m_{j}+1-l}^{(j)}} \bar{S}_{(k_{l}^{j})^{*},k_{l+1}^{j}} \prod_{j=1}^{n} \bar{v}_{k_{m_{j}}^{j}}^{-a_{1}^{(j)}} \dim(k_{m_{j}}^{j}) \bar{S}_{j^{*},k_{1}^{j}}.
\end{eqnarray*}
Here $\bar{v}_{i}=v_{i}^{-1}$ and $\bar{S}_{j^{*},k}=S_{j,k}$.
By (\ref{eq:Dsquared}) we then get
{\allowdisplaybreaks
\begin{eqnarray*}
\tau_{(\mV,\mD)}(M) &=& (\Delta\mD^{-1})^{\sigma_{\os}} \mD^{2g-2-\sum_{j=1}^{n} m_{j}} \\*
 & & \hspace{.4in} \times \sum_{j,k_{1}^{1},\ldots,k_{m_{1}}^{1},k_{1}^{2},\ldots,k_{m_{n}}^{n} \in I} v_{j}^{-b} (\dim(j))^{2-2g-n} \\*
 & & \hspace{.4in} \times \prod_{j=1}^{n} \prod_{l=1}^{m_{j}-1} v_{k_{l}^{j}}^{a_{m_{j}+1-l}^{(j)}} S_{k_{l}^{j},k_{l+1}^{j}} \prod_{j=1}^{n} v_{k_{m_{j}}^{j}}^{a_{1}^{(j)}} \dim(k_{m_{j}}^{j}) S_{j,k_{1}^{j}} \\
 &=& (\Delta\mD^{-1})^{\sigma_{\os}} \mD^{2g-2-\sum_{j=1}^{n} m_{j}} \\*
 & & \hspace{.4in} \times \sum_{j \in I} v_{j}^{-b} (\dim(j))^{2-2g-n} \left( \prod_{i=1}^{n} (SG^{\mC_{i}})_{j,0} \right).
\end{eqnarray*}}\noindent
This expression is identical with (\ref{eq:osurgery}).

\section{The case of $\frsl_{2}(\C)$}\label{sec-The-case}

Let $t=\exp(i\pi/(2r))$, where $r$ is an integer $\geq 2$,
and let $U_{t}$ be the Hopf algebra
considered in \cite[Sect.~8]{ReshetikhinTuraev2}, see below for details.
With notation from
\cite{ReshetikhinTuraev2} we have that
$\left( U_{t},R,v^{-1}, \{ V_{i} \}_{i \in I} \right)$ is a
modular Hopf algebra as defined in \cite[Chap.~XI]{Turaev}. 
Let $\left( \mV_{t} , \{V_{i}\}_{i \in I} \right)$
be the modular category
induced by this modular Hopf algebra, cf.\ \cite[Chap.~XI]{Turaev}. 

Let us recall some notation and results from
\cite{ReshetikhinTuraev2}.
The quantum group $U_{q}(\frsl_{2})$, $q=t^{4}$, is the
$\Q(t)$--algebra with generators $K,K^{-1},X,Y$ subject to the relations
\begin{eqnarray*}
& &XY-YX = \frac{K^{2}-K^{-2}}{t^{2}-t^{-2}}, \\
& &XK = t^{-2} KX, \hspace{.1in} YK=t^{2}KY, \hspace{.1in} KK^{-1}=K^{-1}K=1.
\end{eqnarray*}
The algebra $U_{t}$ is given by the quotient of $U_{q}(\frsl_{2})$ by the
two-sided ideal generated by the elements $X^{r},Y^{r},K^{4r}-1$. 
In the following we will consider $U_{t}$ as an algebra over $\C$.

To determine the RT--invariants of the Seifert manifolds
with non-orientable base we need to determine the signs $\vep_{i}$
in \reflem{selfdual} for all self-dual $i \in I=\{0,1,\ldots,r-2\}$.
It is a well-known fact that all the simple
objects $V_{i}$ in $\mV_{t}$ are self-dual.
Let us provide some details. We let
$$
[k] = \frac{t^{2k}-t^{-2k}}{t^{2}-t^{-2}} = \frac{\sin(\pi k /r)}{\sin(\pi/r)}
$$
for an integer $k$. Let $\alpha \in \{-\sqrt{-1},\sqrt{-1},-1,1\}$.
Then we have irreducible $U_{t}$--modules $\{ V^{i}(\alpha) \}_{i \in I}$
with a basis (over $\C$) of weight vectors $\{e_{n}^{i}(\alpha)\}_{n=0}^{i}$
such that
\begin{eqnarray*}
Ke_{n}^{i}(\alpha) &=& \alpha t^{i-2n} e_{n}^{i}(\alpha), \\
Xe_{n}^{i}(\alpha) &=& \alpha^{2} [n][i+1-n] e_{n-1}^{i}(\alpha), \\
Ye_{n}^{i}(\alpha) &=& e_{n+1}^{i}(\alpha)
\end{eqnarray*}
for $n=0,1,\ldots,i$, where $e_{-1}^{i}(\alpha)=e_{i+1}^{i}(\alpha)=0$.
Let $\{f_{n}^{i}(\alpha)\}_{n=0}^{i}$ be the basis of $V^{i}(\alpha)^{*}$
dual to $\{e_{n}^{i}(\alpha)\}_{n=0}^{i}$ as in
\cite[Sect.~8]{ReshetikhinTuraev2}.
By using the antipode $\gamma$ of $U_{t}$ determined
by $\gamma(K)=K^{-1}$, $\gamma(X)=-t^{2}X$ and $\gamma(Y)=-t^{-2}Y$,
see \cite[(8.1.4)]{ReshetikhinTuraev2}, one gets
\begin{eqnarray*}
Kf_{n}^{i}(\alpha) &=& \alpha^{-1} t^{2n-i} f_{n}^{i}(\alpha), \\
Xf_{n}^{i}(\alpha) &=& -\alpha^{2} t^{2} [n+1][i-n] f_{n+1}^{i}(\alpha), \\
Yf_{n}^{i}(\alpha) &=& -t^{-2} f_{n-1}^{i}(\alpha)
\end{eqnarray*}
for $n=0,1,\ldots,i$, where $f_{-1}^{i}(\alpha)=f_{i+1}^{i}(\alpha)=0$.
We have $V_{i}=V^{i}(1)$, $i \in I$. The following lemma follows by
a straightforward computation using the above $U_{t}$--module
structures.

\begin{lem}\label{isomorphisms}
Let $\alpha,\beta \in \{ -\sqrt{-1},\sqrt{1},-1,1\}$. A 
$\C$--linear map
$h \co V^{i}(\alpha)^{*} \ra V^{i}(\beta)$ is a $U_{t}$--module
isomorphism if and only if $\beta=\alpha^{-1}$ and
\begin{equation}\label{eq:isomorphism}
h(f_{n}^{i}(\alpha))=\delta_{i} (-1)^{n} t^{-2n}e_{i-n}^{i}(\alpha^{-1})
\end{equation}
for a $\delta_{i} \in \C \sm \{0\}$.\HS
\end{lem}

The lemma shows, as claimed above, that the module $V_{i}$ is self-dual
for all $i \in I$.

\begin{lem}\label{signs}
We have $\vep_{i}=(-1)^{i}$ for all $i \in I$. In particular
the modular category $\left( \mV_{t} , \{V_{i}\}_{i \in I} \right)$
is not unimodal.
\end{lem}

\proof
We use \reflem{isomorphisms} together with \reflem{modular-Hopf-algebra}
to determine the signs $\vep_{i}$, $i \in I$. To this end note
that $v=uK^{-2}$, where $u$ is the element $u$ in \reflem{modular-Hopf-algebra}.
As indicated in the beginning of this section we shall use $v^{-1}$ as
the element $v$ in \reflem{modular-Hopf-algebra}. (This is due to different
conventions in \cite{Turaev} and \cite{ReshetikhinTuraev2}.)
We see that $uv^{-1}=K^{2}$ (use that $v^{-1}$ is central).
Fix $i \in I$ and let $\omega =h^{-1}$, where
$h \co V_{i}^{*} \ra V_{i}$ is given by (\ref{eq:isomorphism})
(with $\alpha=1$).
Moreover, let $e_{n}=e_{n}^{i}(1)$, $f_{n}=f_{n}^{i}(1)$. Then
$$
\omega(e_{n})=\delta_{i}^{-1}(-1)^{i-n}t^{2(i-n)}f_{i-n}.
$$
Let $z=G \circ (\omega^{-1})^{*} (f_{i-n}) \in V_{i}$, where
$G \co V_{i}^{**} \ra V_{i}$ is the canonical isomorphism as in
\reflem{modular-Hopf-algebra}. Then
$$
G^{-1}(z)(f_{m})=f_{i-n}(\omega^{-1}(f_{m}))=\delta_{i}(-1)^{m} t^{-2m} f_{i-n}(e_{i-m}) = \delta_{i}(-1)^{n} t^{-2n}\delta_{n,m}.
$$
On the other hand $G^{-1}(z)(f_{m})=f_{m}(z)$, so we see that
$z=\delta_{i}(-1)^{n}t^{-2n}e_{n}$. But then
$$
G \circ (\omega^{-1})^{*} \circ \omega(e_{n}) = \delta_{i}^{-1}(-1)^{i-n}t^{2(i-n)}z=(-1)^{i}t^{2i-4n}e_{n} = (-1)^{i} K^{2} e_{n}.\eqno{\qed}
$$

Let $\kappa(i)$, $i \in I$, be the quantity in (\ref{eq:vector}). 
The $R$--matrix calculation of $\kappa(i)$ gives the result 
$\kappa(i)=(-1)^{i} \Lambda v_{i}^{2} \left( \dim(i) \right)^{-1}$,
$\Lambda=\sum_{u \in I} \left( \dim(u) \right)^{2}$,
so gives together with \reflem{vector} another proof of \reflem{signs}
and the fact that all the simple modules $V_{i}$, $i \in I$, are self-dual.

Let $\bar{v}_{i}$ be equal to the $v_{i}$ in \cite{ReshetikhinTuraev2},
i.e.\ $\bar{v}_{i}\id_{V_{i}}$ is equal to the map $V_{i} \ra V_{i}$ given by
multiplication with $v$, $i \in I$.
Moreover, let $\bar{d}_{i}$ be equal to $d_{i}$
in \cite{ReshetikhinTuraev2}, i.e.\
\begin{displaymath}
\sum_{i \in I} \bar{d}_{i} \bar{v}_{i} S_{i,j} = \bar{v}_{j}^{-1} \dim(j), \hspace{.2in} j \in I.
\end{displaymath}
Since $i = i^{*}$ here, this can also be written
\begin{displaymath}
\sum_{i \in I} \bar{d}_{i} \bar{v}_{i} S_{i^{*},j} = \bar{v}_{j}^{-1} \dim(j), \hspace{.2in} j \in I.
\end{displaymath}
It follows that the $\bar{v}_{i}$ and $\bar{d}_{i}$ are equal to the $v_{i}$ and
$d_{i}$ associated to $\omV_{t}$ in \cite[Sect.~II.3]{Turaev},
where $\omV_{t}$ is the mirror
of $\mV_{t}$, see 
Sect.~\ref{sec-Modular-categories}.
By \cite[Sect.~8.3]{ReshetikhinTuraev2},
\begin{eqnarray*}
\bar{v}_{i} &=& t^{-i(i+2)}, \\
\bar{d}_{i} &=& \sqrt{\frac{2}{r}} \sin(\pi/r) C_{0} \dim(i)
\end{eqnarray*}
for $i \in I$. Here $C_{0} =\exp(\sqrt{-1}d)$ is a square root of 
$C=\sum_{i \in I} \bar{v}_{i}^{-1} \dim(i) \bar{d}_{i}=\exp(2\sqrt{-1}d)$, where
$d=\frac{3\pi(r-2)}{4r}$
and $\dim(i)=[i+1]$.
In particular 
$\bar{d}_{0}=\sqrt{\frac{2}{r}} \sin(\pi/r) \exp(\sqrt{-1}d)$. 
According to \cite[pp.~88--89]{Turaev} we have that 
$\bar{d}_{0}=\Delta_{\omV_{t}} \Lambda^{-1}$ and
$C=\Delta_{\omV_{t}} \bar{d}_{0}$.
(Here we use that the dimensions of any object of $\mV_{t}$ with respect
to $\mV_{t}$ and $\omV_{t}$ are equal, cf.\ \cite[Corollary I.2.8.5]{Turaev},
so $\Lambda$ is the same element in these two categories. By the same
reason $\mD$ is a rank of $\omV_{t}$ if and only if $\mD$ is a rank of $\mV_{t}$.)
We see that
$\Lambda=C \bar{d}_{0}^{-2}= \frac{r}{2} \frac{1}{\sin^{2}(\pi/r)}$.
As a rank we choose
\begin{equation}\label{eq:rank}
\mD=\sqrt{\frac{r}{2}} \frac{1}{\sin(\pi/r)}.
\end{equation}
Let $\{v_{i}\}_{i \in I}$ be the $v_{i}$ associated to $\mV_{t}$. Then
\begin{equation}\label{eq:twist}
v_{i} = \bar{v}_{i}^{-1} = t^{i(i+2)} =t^{-1}t^{(i+1)^2}, \hspace{.2in} i \in I.
\end{equation}
By (\ref{eq:Dsquared}) we get
$\Delta=\Delta_{\mV_{t}}=\mD^{2}\Delta_{\omV_{t}}^{-1}=\mD^{2}\bar{d}_{0} C^{-1}$,
so
\begin{equation}\label{eq:anomalyfactor}
\Delta\mD^{-1}=\mD\bar{d}_{0}C^{-1}=\exp(-\sqrt{-1}d)=\exp\left(\frac{\sqrt{-1}\pi}{4} \frac{3(2-r)}{r}\right).
\end{equation}
The $S$--matrix of $\mV_{t}$ is given by
\begin{equation}\label{eq:Smatrix}
S_{i,j}=\frac{\sin(\pi(i+1)(j+1)/r)}{\sin(\pi/r)}, \hspace{.2in} i,j \in I,
\end{equation} 
cf.\ \cite[Sect.~8]{ReshetikhinTuraev2}.

Inspired by the Chern--Simons path integral
invariant of Witten, see \cite{Witten} and the introduction, the RT--invariant
$\tau_{(\mV_{t},\mD)}$ is also called the quantum $SU(2)$--invariant at
{\it level} $r-2$ and is denoted $\tau_{r}$.
We will take advantage of the fact that the projective action of
$\Mod_{1} = SL(2,\Z)$ considered in
(\ref{eq:e4}) can be normalized to a linear action (in fact to a linear action
of $PSL(2,\Z)$) in the $\frsl_{2}(\C)$--case. Let
$\mR \co PSL(2,\Z) \ra GL(r-1,\C)$
be the unitary representation given by
\begin{equation}\label{eq:C1}
\tilde{\Xi}_{jl}=\sqrt{\frac{2}{r}} \sin\left( \frac{jl\pi}{r} \right),\hspace{.2in}\tilde{\Theta}_{jl}=e^{-\kvi} \exp\left( \frac{i\pi}{2r}j^{2} \right) \delta_{jl}
\end{equation}
for $j,l \in I'=\{1,2,\ldots,r-1\}$.
Here we write $\tilde{M}$ for the matrix $\mR (M)$.
By changing the index set of the basis $\{ b_{i} \}_{i \in I}$
to $I'$ (so that the new $V_{j}$
and $b_{j}$ are equal to the old $V_{j-1}$ and $b_{j-1}$, $j \in I'$)
and comparing (\ref{eq:C1}) with (\ref{eq:rank}),
(\ref{eq:twist}) and (\ref{eq:Smatrix}) we get
\begin{eqnarray}\label{eq:C3b}
S_{jl} &=& \mD\tilde{\Xi}_{jl}, \\
T_{jl} &=& \exp \left( -\frac{i\pi}{2r} \right) \exp\left( \frac{i\pi}{2r} j^{2} \right) \delta_{jl} = e^{\kvi} \exp \left( -\frac{i\pi}{2r} \right) \tilde{\Theta}_{jl} \nonumber
\end{eqnarray}
for $j,l\in I'$.
The representations $\mR$ have been carefully studied by
Jeffrey in \cite{Jeffrey} where she gives rather
explicit formulas for $\tilde{M}$ in terms of
the entries in $M \in SL(2,\Z)$, see also the
proof of \refthm{CS11} below.
These representations are known from the study
of affine Lie algebras, cf.\ \cite{Kac}.

The following corollary is an
$\frsl_{2}(\C)$--version
of \refthm{surgery-formula}.
It simply follows by choosing tuples of integers $\mC_{i}$
such that $B_{i}=B^{\mC_{i}}$, $i=1,2,\ldots,m$. By the first remark
following \refthm{invariants}, $\mC_{i}$ is a continued fraction expansion
of $p_{i}/q_{i}$. Besides note that $\Delta \mD^{-1}=w^{-3}$
with $w=e^{\kvi} \exp \left( -\frac{i\pi}{2r} \right)$ by (\ref{eq:anomalyfactor}),
and that $G^{\mC}=w^{\sum_{k=1}^{m} a_{k} } \mD^{m} \tilde{B}^{\mC}$
for $\mC=(a_{1},\ldots,a_{m}) \in \Z^{m}$
by (\ref{eq:C3b}), (\ref{eq:Bmatrix}), and (\ref{eq:A3}).

\begin{cor}\label{surgery-formula2}
Let the situation be as in {\em \refthm{surgery-formula}} and let $B_{i} \in SL(2,\Z)$
with first column equal to $\pm \left( \begin{array}{c}
		p_{i} \\
		q_{i}
		\end{array}
\right)$, $i=1,2,\ldots,m$. Then
\begin{eqnarray*}
\tau_{r}(M',\Omega) &=& \left( \e^{\kvi} \exp \left( -\frac{i\pi}{2r} \right) \right)^{\sum_{i=1}^{m} \Phi(B_{i}) -3\mu } \\
 & & \hspace{.4in} \times \sum_{\lambda \in \col(L)} \tau_{r}(M,\Gamma(L,\lambda) \cup \Omega) \left( \prod_{i=1}^{m} (\tilde{B}_{i})_{\lambda(L_{i}),1} \right),
\end{eqnarray*}
where $\mu$ is given by {\em (\ref{eq:e7})}.\HS
\end{cor}

The following theorem generalizes results in \cite{Rozansky}
to include the case of Seifert manifolds with non-orientable base. Specifically
the expressions (\ref{eq:C11}) and (\ref{eq:C12}) are equivalent to
\cite[Formulas (2.7) and (2.8)]{Rozansky} for $\ep=\os$.
(Rozansky uses non-normalized Seifert invariants.)
To state the theorem we need some
notation.
Multi-indices are denoted by an underline (e.g.\ $\um$).
For $\underline{k}=(k_{1},\ldots,k_{n}),\underline{l}=(l_{1},\ldots,l_{n}) \in \Z^{n}$,
$\underline{k} < \underline{l}$ if and only if $k_{j} < l_{j}$ for all
$j=1,\ldots,n$. We let
$\underline{1} =(1,\ldots,1)$. For
$\underline{k}=(k_{1},\ldots,k_{n}) \in \Z_{+}^{n}$ we write
$\sum_{\um = \underline{0}}^{\underline{k}}$ for
$\sum_{m_{1}=0}^{k_{1}}\ldots\sum_{m_{n}=0}^{k_{n}}$ etc.
In all expressions below $e$ denotes the Seifert Euler number
(except in factors such as $\e^{\kvi}$).
Let $a_{\os}=2$ and $a_{\ns}=1$.
For a pair of coprime integers $\alpha,\beta$ we let
$\beta^{*}$ be the invers of $\beta$ in the group
of (multiplicative) units in $\Z/\alpha\Z$.

\begin{thm}\label{CS11}
The RT--invariant at level $r-2$ of the Seifert manifold
$M$ with (normalized) Seifert invariants
$(\ep;g \, | \, b;(\alpha_{1},\beta_{1}),\ldots,$ $(\alpha_{n},\beta_{n}))$,
$\ep \in \{ \os, \ns\}$, is
\begin{eqnarray}\label{eq:C11}
\tau_{r}(M) &=& \exp \left( \frac{i \pi}{2r} \left[ 3(a_{\ep}-1)\sign(e) -e - 12\sum_{j=1}^{n} \s (\beta_{j},\alpha_{j}) \right] \right) \\
 && \hspace{.4in} \times (-1)^{a_{\ep}g} \frac{i^{n} r^{a_{\ep}g/2-1}}{2^{n+a_{\ep}g/2-1}} \frac{1}{ \sqrt{\mA}} e^{i\frac{3\pi}{4}(1-a_{\ep})\sign(e)} Z_{\ep}(M;r), \nonumber
\end{eqnarray}
where $\s (\beta_{j},\alpha_{j})$ is given by {\em (\ref{eq:C6})}, $\mA=\prod_{j=1}^{n} \alpha_{j}$, and
\begin{eqnarray}\label{eq:C12}
Z_{\ep}(M;r) &=& \sum_{\gamma =1}^{r-1} (-1)^{\gamma a_{\ep}g} \frac{ \exp \left( \frac{i \pi}{2r} e \gamma^{2} \right) }{ \sin ^{n+a_{\ep}g-2} \left( \frac{\pi}{r} \gamma \right) } \musum \left( \prod_{j=1}^{n} \mu_{j} \right) \\
 & & \hspace{.3in} \times \msum \exp \left( - \frac{i \pi}{r} \gamma \sum_{j=1}^{n} \frac{ 2rm_{j} + \mu_{j}}{\alpha_{j}} \right) \nonumber \\
 & & \hspace{.8in} \times \exp \left( -2 \pi i \sum_{j=1}^{n} \frac{\beta_{j}^{*}}{\alpha_{j}} [ r m_{j}^{2} + \mu_{j} m_{j} ] \right). \nonumber
\end{eqnarray}
The RT--invariant at level $r-2$ of the Seifert manifold $M$
with non-normalized Seifert invariants
$\{\ep;g;(\alpha_{1},\beta_{1}),\ldots,$ $(\alpha_{n},\beta_{n})\}$
is given by the same expression. 
\end{thm}

The theorem is also valid in case $n=0$. In this case one
just has to put $\mA=1$ and 
$\sum_{j=1}^{n} \s (\beta_{j},\alpha_{j})=0$. Moreover, the sum
$\musum \msum$ in $Z_{\ep}(M;r)$ has to be put equal to $1$,
$\ep=\os,\ns$.

\begin{proof}
Let $M=(\ep;g \, | \, b;(\alpha_{1},\beta_{1}),\ldots,$ $(\alpha_{n},\beta_{n}))$,
$\ep \in \{ \os, \ns\}$. Choose tuples of integers
$\mC_{j}=(a_{1}^{(j)},\ldots,a_{m_{j}}^{(j)})$ such that
$B^{\mC_{j}}=\left( \begin{array}{cc}
                \alpha_{j} & \rho_{j} \\
		\beta_{j} & \sigma_{j}
		\end{array}
\right)$ for some $\rho_{j}, \sigma_{j} \in \Z$, $j=1,2,\ldots,n$.
By \refthm{invariants}, the first remark after \refthm{invariants},
\reflem{isomorphisms} and \reflem{signs} we have
\begin{eqnarray*}
\tau_{r}(M) &=& (\Delta\mD^{-1})^{\sigma_{\ep}} \mD^{a_{\ep}g-2-\sum_{j=1}^{n} m_{j}} \\
 && \hspace{.4in} \times \sum_{j=1}^{r-1} (-1)^{(j-1)a_{\ep}g}v_{j}^{-b} \dim(j)^{2-n-a_{\ep}g} \left( \prod_{i=1}^{n} (SG^{\mC_{i}})_{j,1} \right),  
\end{eqnarray*}
where $\sigma_{\ep}$ is given by (\ref{eq:osignature0}) if $\ep=\os$
and by (\ref{eq:nsignature0}) if $\ep=\ns$.
Here $\Delta \mD^{-1}=w^{-3}$, where
$w=e^{\kvi} \exp \left( -\frac{i\pi}{2r} \right)$, see (\ref{eq:anomalyfactor}).
Moreover, $v_{j}=t^{-1}t^{j^{2}}$, see (\ref{eq:twist}),
$\dim(j)=[j]=\sqrt{\frac{2}{r}} \mD \sin \left( \frac{\pi j}{r} \right)$,
and $SG^{\mC_{i}}=w^{\sum_{k=1}^{m_{i}} a_{k}^{(i)}} \mD^{m_{i}+1} \tilde{N}_{i}$,
where $N_{i}=\Xi B^{\mC_{i}}$, so
\begin{equation}\label{eq:in}
\tau_{r}(M) = \alpha_{\ep}(r) \left( \frac{r}{2} \right)^{a_{\ep}g/2+n/2-1} \sum_{j=1}^{r-1} (-1)^{ja_{\ep}g}\frac{t^{-bj^{2}} \prod_{i=1}^{n} (\tilde{N}_{i})_{j,1}}{\sin^{n+a_{\ep}g-2}\left(\frac{\pi j}{r} \right)},
\end{equation}
where
$\alpha_{\ep}(r)=(-1)^{a_{\ep}g}w^{\sum_{i=1}^{n} \sum_{k=1}^{m_{i}} a_{k}^{(i)} -3\sigma_{\ep}} \exp \left( \frac{i\pi}{2r} b\right)$.
By (\ref{eq:osignature}) and the equivalent expression for
$\sigma_{\ns}$ we get
$$
\alpha_{\ep}(r)=(-1)^{a_{\ep}g}w^{\sum_{i=1}^{n} \Phi(B^{\mC_{i}})-3\sum_{j=1}^{n} \sign( \alpha_{j} \beta_{j} )-3(a_{\ep}-1)\sign(e)} \exp \left( \frac{i\pi}{2r} b\right).
$$
If $A = \left( \begin{array}{cc}
                a & b \\
		c & d
		\end{array}
\right) \in SL(2,\Z)$ with $c \neq 0$ 
we have by \cite[Proposition 2.7 (a) and Proposition 2.8]{Jeffrey} that
{\allowdisplaybreaks
\begin{eqnarray*}
\tilde{A}_{j,k} &=& C \sum_{\mu = \pm 1} \sum_{ \stackrel{ \gamma \pmod{2rc} }{ \gamma = j \pmod{2r} }} \mu \exp \left( \frac{i \pi}{2rc} [a \gamma^{2} - 2 \mu \gamma k + d k^{2} ] \right) \\
 &=& C\left\{ \sum_{n=m_{1}}^{|c|-1+m_{1}} \exp \left( \frac{i \pi}{2rc} [a (j + 2rn)^{2} - 2 k (j +2rn) + d k^{2} ] \right) \right. \\*
 & & \hspace{.3in} -\left. \sum_{n=m_{2}}^{|c|-1+m_{2}} \exp \left( \frac{i \pi}{2rc} [a (j + 2rn)^{2} + 2 k (j +2rn) + d k^{2} ] \right) \right\}
\end{eqnarray*}}\noindent
for all $m_{1},m_{2} \in \Z$, where
$C=i \frac{\sign (c)}{\sqrt{2r|c|}} e^{-\kvi \Phi(A)}$. 
For $m_{1}=0$ and $m_{2}=-|c|+1$ we get
$$
\tilde{A}_{j,k} = C\sum_{\mu = \pm 1} \sum_{n=0}^{|c|-1} \mu \exp \left( \frac{i \pi}{2rc} [a (j + 2rn \mu )^{2} - 2 \mu k (j +2rn \mu ) + d k^{2} ] \right).
$$
If
$A_{i} = \left( \begin{array}{cc}
                a_{i} & b_{i} \\
		c_{i} & d_{i}
		\end{array}
\right) \in SL(2,\Z)$
such that $A_{3}=A_{1}A_{2}$ we have
\begin{equation}\label{eq:C7}
\Phi(A_{3})=\Phi(A_{1})+\Phi(A_{2})-3\sign (c_{1}c_{2}c_{3}).
\end{equation}
Since the representation ${\mathcal R}$ is unitary we have ${\mathcal R}(A^{-1})={\mathcal R}(A)^{*}$,
so $\tilde{A}_{j,k} = \overline{(\widetilde{A^{-1}})_{k,j}}$,
where $\bar{\cdot}$ is complex conjugation. Here
$A^{-1} = \left( \begin{array}{cc}
                d & -b \\
		-c & a
		\end{array}
\right)$, and (\ref{eq:C7}) implies that $\Phi(A^{-1})=-\Phi(A)$, so
\begin{eqnarray*}
\tilde{A}_{j,k} &=& i \frac{\sign (c)}{\sqrt{2r|c|}} e^{-\kvi \Phi(A)} \\
 & & \hspace{.2in} \times \sum_{\mu = \pm 1} \sum_{n=0}^{|c|-1} \mu \exp \left( \frac{i \pi}{2rc} [d (k + 2rn \mu)^{2} - 2 \mu j (k +2rn \mu) + a j^{2} ] \right).
\end{eqnarray*}
By this expression we get
\begin{eqnarray*}
(\tilde{N}_{i})_{j,1} &=& i \frac{e^{-\kvi \Phi(N_{i})}}{\sqrt{2r\alpha_{i}}} \sum_{\mu = \pm 1} \sum_{n=0}^{\alpha_{i}-1} \mu \\
 & & \hspace{.3in} \times  \exp \left( \frac{i \pi}{2r\alpha_{i}} [-\beta_{i} j^{2} - 2 j (2rn + \mu) + \rho_{i} (2rn + \mu)^{2} ] \right).
\end{eqnarray*}
By inserting this in (\ref{eq:in})
and using that
$e=-b-\sum_{j=1}^{n} \frac{\beta_{j}}{\alpha_{j}}$ we get $\tau_{r}(M)=\kappa Z_{\ep}(M;r)$, where
$$
\kappa= \frac{i^{n} r^{a_{\ep}g/2-1}}{2^{n+a_{\ep}g/2-1}} \frac{1}{ \sqrt{\mA}} \alpha_{\ep}(r) \exp\left( \frac{i\pi}{2r} \rhoalphasum \right) \exp\left(-\frac{i\pi}{4} \sum_{j=1}^{n} \Phi(N_{j}) \right).
$$
By (\ref{eq:C7}) we have $\Phi(N_{i})=\Phi(B^{\mC_{i}})-3\sign(\alpha_{i}\beta_{i})$ and get
\begin{eqnarray*}
\kappa &=& (-1)^{a_{\ep}g}\frac{i^{n} r^{a_{\ep}g/2-1}}{2^{n+a_{\ep}g/2-1}} \frac{1}{ \sqrt{\mA}}\exp\left( i\frac{3\pi}{4} (1-a_{\ep}) \sign(e) \right) \\
 & & \hspace{0.3in} \times \exp \left( \frac{i\pi}{2r} \left[ 3(a_{\ep}-1)\sign(e) + b + \rhoalphasum - \sum_{j=1}^{n} \Phi(N_{j}) \right]\right).
\end{eqnarray*}
The theorem now follows by using (\ref{eq:C5}) together with the
facts that $\s (a,b)=\s (a',b)$ if $a'a \equiv 1 \pmod{b}$
and $\s (-a,b)=-\s (a,b)$,
cf.\ \cite[Chap.\ 3]{RademacherGrosswald}.
The case with non-normalized Seifert invariants follows as above by
letting $b$ be equal to zero everywhere.
\end{proof}

Let $M=(\ep;g \, | \, b;(\alpha_{1},\beta_{1}),\ldots,(\alpha_{n},\beta_{n}))$.
By (\ref{eq:in}) we have
the following compact formula
\begin{equation}\label{eq:compact}
\tau_{r}(M)=\alpha_{\ep}(r) \sum_{j=1}^{r-1} (-1)^{ja_{\ep}g} \frac{t^{-bj^{2}} \prod_{i=1}^{n} (\tilde{N}_{i})_{j,1}}{\tilde{\Xi}_{j,1}^{n+a_{\ep}g-2}},
\end{equation}
where
$
\alpha_{\ep}(r)=(-1)^{a_{\ep}g}w^{\sum_{i=1}^{n} \Phi(N_{i})-3(a_{\ep}-1)\sign(e)} \exp \left( \frac{i\pi}{2r} b\right)
$, $w=e^{\kvi} \exp \left( -\frac{i\pi}{2r} \right)$,
and $N_{i}=\left( \begin{array}{cc}
                -\beta_{j} & -\sigma_{j} \\
		\alpha_{j} & \rho_{j}
		\end{array}
\right)$ for any integers $\rho_{j}$, $\sigma_{j}$ such that
$\alpha_{j}\sigma_{j}-\beta_{j}\rho_{j}=1$.

Let us finally give a formula for $\tau_{r}(L(p,q))$.
To this end let $b,d$ be any integers such that
$U=\left( \begin{array}{cc}
                q & b \\
		p & d
		\end{array}
\right) \in SL(2,\Z)$. Assume $q \neq 0$, let $V=-\Xi U =\left( \begin{array}{cc}
                p & d \\
		-q & -b
		\end{array}
\right)$, and let $C'=(a_{1},a_{2},\ldots,a_{m-1}) \in \Z^{m-1}$ such that
$B^{\mC'}=V$. Then $\mC'$ is a continued fraction expansion of $-p/q$ and
$U=\Xi V=B^{\mC}$ where $\mC=(a_{1},a_{2},\ldots,a_{m-1},0)$.
By \refcor{lens-spaces} and (\ref{eq:C3b}) we therefore get
\begin{equation}\label{eq:lenssl}
\tau_{r}(L(p,q)) = \left( \e^{\kvi} \exp \left( -\frac{i\pi}{2r} \right) \right) ^{\Phi(U)} \tilde{U}_{1,1}.
\end{equation}
If $q=0$ we have $p=1$ and $L(p,q)=S^{3}$. In this case we have
$U=\Xi \Theta^{d}$ and we immediately find from
(\ref{eq:C3b}) that (\ref{eq:lenssl}) is also true in this case.
The identity (\ref{eq:lenssl}) coincides with \cite[Formula (3.7)]{Jeffrey}.

\begin{rem}
It should not come as a surprise that we find the same result as Rozansky for
the invariants of Seifert manifolds with orientable base. The calculation 
in \cite{Rozansky} of
these invariants follows the very same line as in
the first part of Sect.~\ref{sec-A-second}. He uses a surgery formula
\cite[Formula (1.6)]{Rozansky} which is identical with the surgery formula in
\refcor{surgery-formula2} and a Verlinde formula
\cite[Formula (2.4)]{Rozansky}
which by (\ref{eq:C3b}) is identical with the Verlinde formula 
(\ref{eq:Verlinde}) of Turaev.
\end{rem}

\begin{rem}
In more recent literature the symbol $U_{q}(\frsl_{2}(\C))$
normally refers to a Hopf algebra defined in a slightly different
way than in the above text. It is well known \cite{Kirillov},
\cite{BakalovKirillov}, \cite{Le} that Lusztig's version 
\cite[Part V]{Lusztig}
of quantum deformations of simple complex finite
dimensional Lie algebras at roots of unity is particular well
suited to produce modular categories. Let us specialize to
the $\frsl_{2}$--case. Let $\theta=\exp(i\pi/r)$, $r$ an integer
$\geq 2$, and let $U_{\theta}(\frsl_{2}(\C))$ be Lusztig's version
of the quantum group associated to $\theta$ and $\frsl_{2}(\C)$.
This is a Hopf algebra over $\C$, see the above references
for the definition. (The root of unity $\theta$ is denoted
$q$ in \cite{BakalovKirillov}
and $\vep$ in \cite{Kirillov}.)
The representation theory of $U_{\theta}(\frsl_{2}(\C))$ induces
a modular category $\left( \mV_{\theta}', \{  V_{i}' \}_{i \in I} \right)$,
$I =\{0,1,\ldots,r-2\}$, with $S$-- and $T$--matrices
identical with the $S$-- and $T$--matrices for the modular 
category $\mV_{t}$, $t=\exp(i\pi/(2r))$, considered above, cf.
\cite[Theorem 3.9]{Kirillov}, 
\cite[Theorem 3.3.20]{BakalovKirillov}.
(One should note different notation in \cite{Kirillov} and
\cite{BakalovKirillov}. Note that the $\tilde{s}$--matrix
in \cite{BakalovKirillov} is Turaev's $S$--matrix of the mirror of
$\mV_{\theta}'$, i.e.\ $\tilde{s}_{i,j}=S_{i^{*},j}$, $i,j \in I$,
and that $s$ in \cite{Kirillov} is identical with $\tilde{s}$ in
\cite{BakalovKirillov} and vice versa.) 
The dimension of the simple object $V_{i}'$ of $\mV_{\theta}'$
is equal to the dimension of the simple object $V_{i}$ of $\mV_{t}$.
The categories $\mV_{t}$ and $\mV_{\theta}'$
therefore also have the same ranks and the same $\Delta$.
Similar to the proofs of \reflem{isomorphisms} and \reflem{signs}
we find that $V_{i}'$ is self-dual with associated $\vep_{i}=(-1)^{i}$,
$i \in I$. We conclude that $\mV_{\theta}'$ and
$\mV_{t}$ give the same invariants of the Seifert manifolds.
Probably these two categories are equivalent giving the same
invariants for all closed oriented $3$--manifolds, but we will not check
the details here.

In \cite[Problems, question 8 p.~571]{Turaev}
it is asked whether there
exist unitary (or at least Hermitian) modular categories
that are not unimodal. By combining the above with
\cite{Kirillov}, \cite{Wenzl} we can answer this question by a yes. The
non-unimodal modular categories
$\left( \mV_{\theta}', \{  V_{i}' \}_{i \in I} \right)$
provide such examples.

In \cite[Chap.~XII]{Turaev} Turaev constructs a unimodular category 
$\left( \mV_{p}(a), \{  W_{i} \}_{i \in I} \right)$
for any primitive $4r$'th root of unity $a$ using
Kauffman's skein theoretical approach to the Jones polynomial together
with the Jones--Wenzl idempotens. Here $I=\{0,1,\ldots,r-2\}$ as above
and all the simple objects are self-dual.
(In \cite{Turaev} $W_{i}$ is denoted $V_{i}$.)
In \cite[Problems, question 24 p.~572]{Turaev} it is asked whether
$\mV_{p}(a)$ (with ground ring $\C$) 
is equivalent (as modular category) to the modular category
$\left( \mV_{q}'', \{  V_{i}'' \}_{i \in I} \right)$
induced by the representation theory of $U_{q}(\frsl_{2}(\C))$
for $q=-a^{2}$.
Here $U_{q}(\frsl_{2}(\C))$ is given in
\cite[Sections VI.1 and VII.1]{Kassel}
and differs slightly from the $U_{q}(\frsl_{2}(\C))$ given in
this section and also from Lusztig's version, see above.
For any primitive $2r$'th root of unity $q$ the modular category $\mV_{q}''$
is non-unimodal, so the answer
to the above question is no.
In fact, by using arguments similar to 
the proofs of \reflem{isomorphisms} and \reflem{signs}, one
finds that $V_{i}''$ is self-dual with associated $\vep_{i}=(-1)^{i}$, $i \in I$.
One can construct a non-unimodal 
modular category $\mV_{p}'(a)$ by changing the twist $\theta$ in
$\mV_{p}(a)$ slightly preserving all other structure. In fact one can construct
a new twist $\theta'$ satisfying
$\theta_{W_{i}}'=(-1)^{i}\theta_{W_{i}}$, $i \in I$.
(One simply changes the twist $\theta_{n}$, $n=0,1,2,\ldots$,
in the skein category in
\cite[Sect.~XII.2]{Turaev} into $(-1)^{n}\theta_{n}$.)
By the definition of the elements $\vep_{i}$,
see Sect.~\ref{sec-Modular-categories} above \reflem{selfdual},
we immediately get $\vep_{i}=(-1)^{i}$, $i \in I$, for
$\mV_{p}'(a)$.    
In his thesis \cite{Thys1} H.\ Thys 
shows that the modular category $\mV_{p}'(a)$ is equivalent
to the modular category $\mV_{q}''$, if $q$ is a primitive root
of unity satisfying $q=a^{2}$, see also \cite{Thys2}. By
using the twist $\theta'$ instead of $\theta$ in
the last part of the proof of \cite[Theorem XII.7.1]{Turaev}
and in \cite[Exercise XII.6.10 1)]{Turaev}
one finds that the $S$-- and $T$--matrices for $\mV_{p}'(a)$ are identical
to these matrices for $\mV_{a}$.
 Moreover the dimension of the simple object $W_{i}$ of $\mV_{p}'(a)$
is equal to the dimension of the simple object $V_{i}$ of $\mV_{a}$
, $i \in I$ (use \cite[Sect.~XII.6.8]{Turaev}).
We conclude that these two modular categories
give the same invariants of the Seifert manifolds.
Probably these two categories are equivalent giving the same
invariants for all closed oriented $3$--manifolds, but we will not check
the details here.
\end{rem}

\section{Appendices}

\sh{A. Normalizations of the RT--invariants}

As a convenience to the reader we compare in this appendix
the normalizations of the
RT--invariants used in the literature in particular the ones used in
\cite{ReshetikhinTuraev2}, \cite{KirbyMelvin1},
\cite{TuraevWenzl1}, \cite{Le} and \cite{Turaev}.
The invariants of $3$--manifolds with embedded colored ribbon graphs
constructed in \cite{ReshetikhinTuraev2}, see also
\cite{TuraevWenzl1}, are based on modular Hopf algebras.
The definition of a modular Hopf algebra in
\cite[Chap.~XI]{Turaev} is slightly simplified compared to \cite{ReshetikhinTuraev2},
\cite{TuraevWenzl1}.
If $\left( A,R,v, \{ V_{i} \}_{i \in I} \right)$ is a
modular Hopf algebra as defined in \cite{ReshetikhinTuraev2},
\cite{TuraevWenzl1}, then
$\left( A,R,v^{-1}, \{ V_{i} \}_{i \in I} \right)$ is a
modular Hopf algebra as defined in \cite{Turaev}.
(The definition of the
$v_{i}$'s on p.~557 in \cite{ReshetikhinTuraev2} has to be changed
according to \cite{TuraevWenzl1}.
That is, $v_{i}\id_{V_{i}}$ should be
equal to the map $V_{i} \ra V_{i}$ given by multiplication
with $v^{-1}$ instead of the map given by multiplication with
$v$.)
Let $\left( \mV, \{ V_{i} \}_{i \in I} \right)$
be the modular category induced by the modular Hopf algebra
$\left( A,R,v^{-1}, \{ V_{i} \}_{i \in I} \right)$,
cf.\ \cite[Chap.~XI]{Turaev}, and
let $(M,\Omega)$ be as in (\ref{eq:A4}). The invariant of the pair
$(M,\Omega)$ as defined in \cite{ReshetikhinTuraev2} is given by
\begin{displaymath}
\mF(M,\Omega)=C^{-\sigma_{-}(L)} \sum_{\lambda \in \col(L)} \left( \prod_{i=1}^{m} d_{\lambda(L_{i})} \right) F_{\mV}(\Gamma(L,\lambda) \cup \Omega).
\end{displaymath}
Here $C=\sum_{i \in I} v_{i}^{-1} \dim(i)d_{i}$,
where $\{ d_{i} \}_{i \in I}$ is the
unique solution to
\begin{equation}\label{eq:d}
\sum_{i \in I} v_{i} S_{i,j} d_{i} = v_{j}^{-1} \dim(j) \hspace{.2in}, j \in I,
\end{equation}
where $S$ is the $S$--matrix of $\mV$.
Moreover $\sigma_{-}(L)$ is the number of negative eigenvalues
of the intersection form of $W_{L}$.
By comparing with \cite[Sect.~II.3]{Turaev}
we see that $C=xd_{0}$, where $x=\sum_{i \in I} d_{i}\dim(i)=\Delta=\Delta_{\mV}$, and
$
\mF(M,\Omega)=\tau_{\mV}'(M,\Omega)=(\Delta \mD^{-1})^{b_{1}(M)} \mD \tau_{(\mV,\mD)}(M,\Omega),
$
where $\mD$ is a rank of $\mV$ and $b_{1}(M)$ is the first betti number of $M$.
According to \cite[p.~89]{Turaev}
we also have $C=xd_{0}=(\Delta\mD^{-1})^{2}$. In \cite{TuraevWenzl1}
the invariant $\mF(M,\Omega)$ is slightly changed to
\begin{displaymath}
\tau_{A}(M,\Omega)=C_{0}^{\sigma(L)-m} \sum_{\lambda \in \col(L)} \left( \prod_{i=1}^{m} d_{\lambda(L_{i})} \right) F(\Gamma(L,\lambda) \cup \Omega),
\end{displaymath}
where $C_{0}$ is a square root of $C$. For $C_{0}=\Delta \mD^{-1}$ we
have
\begin{displaymath}
\tau_{A}(M,\Omega)=C_{0}^{-b_{1}(M)} \mF(M,\Omega)=\mD \tau_{(\mV,\mD)}(M,\Omega),
\end{displaymath}
which follows by using that $\sigma_{-}(L)=(m-b_{1}(M)-\sigma(L))/2$.
In case $A=U_{t}$, $t=\exp(i\pi/(2r))$, $r \geq 2$,
see \cite[Sect.~8]{ReshetikhinTuraev2}
and the beginning of Sect.~\ref{sec-The-case} in this paper,
$\tau_{A}(M)=\tau_{A}(M,\emptyset)$ is equal to the invariant
$\tau_{r}(M)$ in \cite{KirbyMelvin1}.

To compare with \cite{Le} we use a more symmetric expression
for $\tau_{(\mV,\mD)}$.
To this end let $\Delta^{-}=\Delta$ and $\Delta^{+}=\Delta_{\omV}$, so
$\Delta^{\pm}=\sum_{i \in I} v_{i}^{\pm 1} \left( \dim(i) \right)^{2}$.
Moreover, let $\sigma_{+}(L)$
be the number of positive eigenvalues
of the intersection form of $W_{L}$. Then, by using (\ref{eq:Dsquared}) and the
above formula for $\sigma_{-}(L)$, one gets
\begin{eqnarray*}
\tau_{(\mV,\mD)}(M,\Omega) &=& \mD^{-b_{1}(M)-1}(\Delta^{+})^{-\sigma_{+}(L)}(\Delta^{-})^{-\sigma_{-}(L)} \\
 && \hspace{.4in} \times \sum_{\lambda \in \col(L)} \left( \prod_{i=1}^{m} \dim(\lambda(L_{i})) \right) F(\Gamma(L,\lambda) \cup \Omega).
\end{eqnarray*}
The invariant $\mD^{b_{1}(M)+1}\tau_{(\mV,\mD)}(M,\Omega)$ is the
invariant considered in \cite{Le} in case the modular categories are
the ones induced by the quantum groups associated to simple finite dimensional
complex Lie algebras.

\sh{B. Framed links in closed oriented $3$--manifolds}

In this appendix we discuss different ways of presenting a
framing of a link $L$ in an arbitrary closed oriented $3$--manifold
$M$. We will here explicitly work in the smooth category so
we can use differential topological concepts.
To simplify writing we restrict to the case of knots.
The generalization to links will be obvious.

\rk{Three ways of defining a framing} 
Let $K$ be a knot in a closed orientable $3$--manifold $M$, let $TM|_{K}$
be the restriction of the tangent bundle of $M$ to $K$, and let
$NK=TM|_{K}/TK$ be the normal bundle of $K$. Since $M$
and $K$ are orientable, $NK$ is an orientable $2$--dimensional real
vector bundle over $K$. Isomorphism classes of oriented
$2$--dimensional real vector bundles 
over $K \cong S^{1}$ are in 1-1 correspondence
with $\pi_{1}(BSO(2)) \cong \pi_{0}(SO(2)) =0$, so
$NK$ is trivializable. We let $\mS_{\tr}(K)$ be the
set of isotopy classes of trivializations of $NK$.
There is a
1-1 correspondence between $\mS_{\tr}(K)$ and
$\pi_{1}(GL(2,\R)) \cong \pi_{1}(O(2)) \cong \Z$
through homotopy classes of transition functions.
However, there is in general no canonical choise
of this bijection.
A normal vector field on $K$ is a
nowhere vanishing section in $NK$. Two normal vector fields
on $K$ are homotopic if they can be deformed into one another
within the class of normal vector fields on $K$. We let
$\mS_{\nvf}(K)$ be the set of
homotopy classes of normal vector fields on $K$.
Finally let $\mS_{\emb}(K)$ be the set of isotopy classes of
embeddings $Q \co B^{2} \times S^{1} \ra M$ with
$Q(0 \times S^{1})=K$ (nothing about orientations here
contrary to \refconv{conventions} in
Sect.~\ref{sec-A-generalized}). A framing of $K$
is an element in one of the sets $\mS_{\nvf}(K)$,
$\mS_{\tr}(K)$, $\mS_{\emb}(K)$.

\rk{Claim} {\sl We have a diagram of maps

\bigskip
\centerline{
\xymatrix{
  {\mS_{\emb}(K)}  \ar@<1ex>[rr]^{J} \ar[dr]_{F} 
   & & {\mS_{\tr}(K)} \ar@<1ex>[ll]^{H_{1},H_{2}} \\
   & {\mS_{\nvf}(K)} \ar[ur]_{G_{1},G_{2}} }
}
\bigskip 

\noindent with $J \circ H_{\nu}$ and $F \circ H_{\nu} \circ G_{\nu}$ the
identity maps, $\nu=1,2$. In particular $H_{\nu}$ and $G_{\nu}$ are injective,
$\nu=1,2$,
and $F$ and $J$ surjective. The images of $G_{1}$ and
$G_{2}$ have the same cardinality and they are disjoint with union $\mS_{\tr}(K)$. 
The union of the images of $H_{1}$ and $H_{2}$ is $\mS_{\emb}(K)$.
Fix an orientation on $K$ and let $-K$ be $K$ with the opposite
orientation. Then $H_{1}=H_{2}$ if $K$ and $-K$ are isotopic,
so in particular this map is an isomorphism (with inverse $J$).
If $K$ and $-K$ are non-isotopic then $H_{1}$ and $H_{2}$ have disjoint
images with the same cardinality.}

\begin{proof}[Proof of claim] 
Let $Q$ be an embedding as above and let $\xi(x)=Q(0,x)$. Define
a normal vector field on $K$ by
$X_{Q}(\xi(x))=\Pi_{K} \circ T_{(0,x)}Q(e_{1},0)$,
where $e_{1}$ is the first standard unit vector in
$\R^{2}$ and $\Pi_{K} \co TM|_{K} \ra NK$ is the projection.
The map $Q \mapsto X_{Q}$ induces a map
$F \co \mS_{\emb}(K) \ra \mS_{\nvf}(K)$. By
$\sigma_{Q}(y,\xi(x))=\Pi_{K} \circ T_{(0,x)}Q(y,0)$
we get a trivialization of $NK$.
The map $Q \mapsto \sigma_{Q}$ induces a map
$J \co \mS_{\emb}(K) \ra \mS_{\tr}(K)$. Let $X$
be a normal vector field on $K$. Fix an
orientation of $NK$ and choose a normal vector
field $Y$ on $K$ so $\{X,Y\}$ is a positively
oriented frame for $NK$. Let $\sigma_{X}$ be the
corresponding trivialization of $NK$,
i.e.\ $\sigma_{X}(ue_{1}+ve_{2},p)=uX(p)+vY(p)$. The
map $X \mapsto \sigma_{X}$ induces a map
$G_{1} \co \mS_{\nvf}(K) \ra \mS_{\tr}(K)$. Let $G_{2}$ be
defined as $G_{1}$ but using the opposite orientation of $NK$. Finally for
a trivialization $\sigma$ of $NK$, a parametrisation
$\xi \co S^{1} \ra K$, and a tubular map $\tau \co NK \ra M$
we get an embedding
$Q_{\sigma} \co B^{2} \times S^{1} \ra M$ by
$Q_{\sigma}(y,x)=\tau(\sigma(y,\xi(x)))$. 
Here $\tau \co NK \ra M$ is an embedding, which on $K$ is the inclusion
$K \subset M$ and for which the differential induces the identity on the
zero-section, cf.\ \cite[p.~123]{BrockerJanich}.
The map $\sigma \mapsto Q_{\sigma}$ induces
a map $H_{1} \co \mS_{\tr}(K) \ra \mS_{\emb}(K)$.
Let $H_{2}$ be defined as $H_{1}$ but using $\xi \circ \kappa$
instead of $\xi$, where $\kappa \co S^{1} \ra S^{1}$ is an 
orientation reversing diffeomorphism. 
By the property of the differential of $\tau$,
we get immediately the first claim
(use that
$\Pi_{K} \circ T_{(0,x)}Q_{\sigma}|_{T_{0}B^{2} \oplus 0} = \sigma(-,\xi(x))$.)
Let $g \co \R^{2} \ra \R^{2}$ be an orientation reversing diffeomorphism.
For an embedding
$Q \co B^{2} \times S^{1} \ra M$ we let
$\bar{Q}=Q \circ (\id_{B^{2}} \times \kappa)$ and
$\hat{Q}=Q \circ (g|_{B^{2}} \times \id_{S^{1}})$.
Similarly for a trivialization
$\sigma$ of $NK$ we let
$\hat{\sigma} = \sigma \circ (g \times \id_{K})$.
Note that $\sigma$ and $\hat{\sigma}$ are non-isotopic.
The claims about $G_{1}$ and $G_{2}$ follows then
basically by the observation
$G_{2}([X])=[\hat{\sigma}_{X}]$.
Fix an orientation of $K$ and let $-K$ be $K$
with the opposite orientation. Then $Q$ and
$\bar{Q}$ are isotopic if and only if $K$ and $-K$
are isotopic.
The claims about $H_{1}$ and $H_{2}$ then basically
follow from the observation
$H_{2}([\sigma])=\bar{Q}_{\sigma}$.
(Use that the isotopy
class of $Q$ is completely determined by $Q|_{0 \times S^{1}}$
and $\Pi_{K} \circ T_{(0,x)}Q|_{T_{0}B^{2} \oplus 0}$.)
\end{proof}

Note that there are oriented knots $K$
for which $K$ and $-K$ are not isotopic, cf.\ \cite{Trotter}. Also note that
for an embedding $Q \co B^{2} \times S^{1} \ra M$, $F$ maps the isotopy classes
of $Q$, $\hat{Q}$, $\bar{Q}$, and $\hat{\bar{Q}}$ to the same point.
Here $Q$ and $\hat{Q}$ are always non-isotopic.

\rk{Integral homology spheres} Let us consider the case where $M=S^{3}$
or more generally where $M$ is an integral homology
sphere (meaning that $H_{*}(M;\Z)=H_{*}(S^{3};\Z)$).
Then we have a well-defined linking number $\lk(.\,,.)$
between knots in $M$. If $Q \co B^{2} \times S^{1} \ra M$
is an embedding with $Q(0 \times S^{1})=K$ we let
$K'$ be the knot $Q(e_{1} \times S^{1})$. Fix an
orientation of $K$ and give $K'$ the induced orientation,
i.e.\ $[K]=[K']$ in $H_{1}(U;\Z)$, $U=Q(B^{2} \times S^{1})$.
Note that we get the same parallel $K'$ (up to isotopy) if we use
$\hat{Q}$ or $\bar{Q}$ instead (see the above proof). We therefore
have an identification $\mS_{\nvf}(K) \cong \Z$ 
by $\lk(K,K')$. The framing corresponding to zero is
sometimes called the {\it preferred framing}, cf.\ \cite[p.~31 and p.~136]{Rolfsen1}.
(We also have $\mS_{\nvf}(K) \cong \Z$ in the
general case, see the claim above, but we do not
in general have a canonical choice of a zero.)

\rk{Notes on surgery} Assume that $M$ is an arbitrary closed
oriented $3$--manifold and
that $Q \co B^{2} \times S^{1} \ra M$ is an
orientation preserving embedding, where $B^{2} \times S^{1}$
is the oriented standard solid torus and
$U=Q(B^{2} \times S^{1})$ is given the orientation induced by that
of $M$, see \refconv{conventions}.
Then $\hat{\bar{Q}}$ is also orientation preserving.
However from a surgical point of view this causes no
problems since rational surgery along
$(K,Q)$, $K=Q(0 \times S^{1})$, with surgery coefficient $p/q$ as
defined in (\ref{eq:b4}) is identical with
rational surgery along $(K,\hat{\bar{Q}})$ with the same surgery
coefficient $p/q$. If we change the orientation of
$M$ ($M$ is connected) then by \refconv{conventions}
we must use $\bar{Q}$ (or $\hat{Q}$) when
doing surgery along $K$. 
This changes the signs of all surgery coefficients
for a given surgery.
If $M=S^{3}=\R^{3} \cup \{ \infty \}$
given the standard right-handed orientation,
then rational surgery on $M$
along $K$ as defined in (\ref{eq:b4}),
where $K$ is given the preferred
framing,
is ordinary rational surgery on $S^{3}$
along $K$ as defined in
e.g.\ \cite[Sect.~9.F]{Rolfsen1}.

\Addresses

\recd


\begin{thebibliography}{[TW2]}
\bibitem[A]{Andersen} J. E. Andersen, {\em The Witten invariant of finite order mapping tori I}, to appear in J. Reine Angew. Math.
\bibitem[At]{Atiyah} M. F. Atiyah, {\em Topological quantum field theories}, Inst. Hautes \'{E}tudes Sci. Publ. Math. {\bf 68} (1989), 175--186.
\bibitem[AS1]{AxelrodSinger1} S. Axelrod, I. M. Singer, {\em Chern--Simons perturbation theory}, Proceedings of the XXth International Conference on Differential Geometric Methods in Theoretical Physics, Vol. {\bf 1}, {\bf 2}, World Sci. Publishing (1992), 3--45.
\bibitem[AS2]{AxelrodSinger2} S. Axelrod, I. M. Singer, {\em Chern--Simons perturbation theory. II}, J. Differential Geom. {\bf 39} (1994), no. 1, 173--213. 
\bibitem[BK]{BakalovKirillov} B. Bakalov, A. Kirillov Jr., {\em Lectures on tensor categories and modular functors}, University Lecture Ser. {\bf 21}, Amer. Math. Soc. (2001).
\bibitem[BJ]{BrockerJanich} T. Br\"{o}cker, K. J\"{a}nich, {\em Introduction to differential topology}, Cambridge University Press (1982).
\bibitem[FG]{FreedGompf} D. S. Freed, R. E. Gompf, {\em Computer calculation of Witten's $3$--manifold invariant}, Comm. Math. Phys. {\bf 141} (1991), no. 1, 79--117.
\bibitem[G]{Garoufalidis} S. Garoufalidis, {\em Relations among $3$--manifold invariants}, D. Phil. thesis, The University of Chicago (1992).
\bibitem[Ha1]{Hansen1} S. K. Hansen, {\em Reshetikhin--Turaev invariants of Seifert $3$--manifolds, and their asymptotic expansions}, D. Phil. thesis, University of Aarhus (1999).
\bibitem[Ha2]{Hansen2} S. K. Hansen, {\em Analytic asymptotic expansions of the Reshetikhin--Turaev invariants of Seifert $3$--manifolds for $\frsl_{2}(\C)$}, in preparation.
\bibitem[JN]{JankinsNeumann} M. Jankins, W. D. Neumann, {\em Lectures on Seifert manifolds}, Brandeis Lecture Notes {\bf 2}, Brandeis University (1981).
\bibitem[J]{Jeffrey} L. C. Jeffrey, {\em Chern--Simons--Witten invariants of lens spaces and torus bundles, and the semiclassical approximation}, Comm. Math. Phys. {\bf 147} (1992), no. 3, 563--604.
\bibitem[JL]{JohnsonLapidus} G. W. Johnson, M. L. Lapidus, {\em The Feynman integral and Feynman's operational calculus}, Oxford mathematical monographs, Clarendon Press, Oxford University Press (2000).
\bibitem[Jo]{Jones} V. F. R. Jones, {\em A polynomial invariant for knots via von Neumann algebras}, Bull. Amer. Math. Soc. {\bf 12} (1985), no. 1, 103--111 
\bibitem[K]{Kac} V. G. Kac, {\em Infinite dimensional Lie algebras}, 3rd ed., Cambridge University Press (1990).
\bibitem[Ka]{Kassel} C. Kassel, {\em Quantum groups}, Grad. Texts in Math. {\bf 155}, Springer-Verlag (1995).
\bibitem[Ki]{Kirby} R. C. Kirby, {\em The topology of $4$--manifolds}, Lecture Notes in Math. {\bf 1374}, Springer-Verlag (1989).
\bibitem[KM1]{KirbyMelvin1} R. C. Kirby, P. M. Melvin, {\em On the $3$--manifold invariants of Witten and Reshetikhin--Tureav for $sl(2, \C)$}, Invent. Math. {\bf 105} (1991), no. 3, 473--545.
\bibitem[KM2]{KirbyMelvin2} R. C. Kirby, P. M. Melvin, {\em Dedekind sums, $\mu$-invariants and the signature cocycle}, Math. Ann. {\bf 299} (1994), no. 2, 231--267.
\bibitem[Kir]{Kirillov} A. A. Kirillov, {\em On an inner product in modular tensor categories}, J. Amer. Math. Soc. {\bf 9} (1996), no. 4, 1135--1169.
\bibitem[LR]{LawrenceRozansky} R. Lawrence, L. Rozansky, {\em Witten--Reshetikhin--Turaev invariants of Seifert manifolds}, Comm. Math. Phys. {\bf 205} (1999), no. 2, 287--314.
\bibitem[Le]{Le} T. T. Q. Le, {\em Quantum invariants of $3$--manifolds: Integrality, splitting, and perturbative expansions}, preprint math.QA/0004099.
\bibitem[Li1]{Lickorish1} W. B. R. Lickorish, {\em A representation of orientable combinatorial $3$--manifolds}, Ann. of Math. {\bf 76} (1962), 531--540.
\bibitem[Li2]{Lickorish2} W. B. R. Lickorish, {\em The Skein method for three-manifold invariants}, J. Knot Theory Ramifications {\bf 2} (1993), no. 2, 171--194.
\bibitem[Lu]{Lusztig} G. Lusztig, {\em Introduction to quantum groups}, Birkh\"{a}user (1993).
\bibitem[Ma]{MacLane} S. MacLane, {\em Categories for the working mathematician}, Grad. Texts in Math. {\bf 5}, 2nd ed., Springer-Verlag (1998).
\bibitem[M]{Montesinos} J. M. Montesinos, {\em Classical tessellations and three-manifolds}, Universitext, Springer-Verlag (1987).
\bibitem[N]{Neil} J. R. Neil, {\em Combinatorial calculations of the various normalizations of the Witten invariants for $3$--manifolds}, J. Knot Theory Ramifications {\bf 1} (1992), no. 4, 407--449.
\bibitem[Ne]{Neumann} W. D. Neumann, {\em $S^{1}$--actions and the $\alpha$--invariant of their involutions}, Bonn. Math. Schr. {\bf 44} (1970).
\bibitem[NR]{NeumannRaymond} W. D. Neumann, F. Raymond, {\em Seifert manifolds, plumbing, $\mu$--invariant and orientation reversing maps}, Algebraic and geometric topology, Proceedings, Santa Barbara 1977, Lecture Notes in Math. {\bf 664}, Springer-Verlag (1978), 163--196. 
\bibitem[O]{Orlik} P. Orlik, {\em Seifert manifolds}, Lecture Notes in Math. {\bf 291}, Springer-Verlag (1972). 
\bibitem[PS]{PrasolovSosinski} V. V. Prasolov, A. B. Sossinsky, {\em Knots, links, braids and $3$--manifolds. An introduction to the new invariants in low-dimensional topology}, Transl. Math. Monogr. {\bf 154}, Amer. Math. Soc. (1997).
\bibitem[Q]{Quinn} F. Quinn, {\em Lectures on axiomatic topological quantum field theory}. In {\em Geometry and quantum field theory} (Parc City, UT, 1991), IAS/Park City Math. Ser. {\bf 1}, Amer. Math. Soc. (1995), 323--453.
\bibitem[RG]{RademacherGrosswald} H. Rademacher, E. Grosswald, {\em Dedekind sums}, The Carus Mathematical Monographs {\bf 16}, The Mathematical Association of America (1972).
\bibitem[RT1]{ReshetikhinTuraev1} N. Reshetikhin, V. G. Turaev, {\em Ribbon graphs and their invariants derived from quantum groups}, Comm. Math. Phys. {\bf 127} (1990), no. 1, 1--26.
\bibitem[RT2]{ReshetikhinTuraev2} N. Reshetikhin, V. G. Tureav, {\em Invariants of $3$--manifolds via link polynomials and quantum groups}, Invent. Math. {\bf 103} (1991), no. 3, 547--597.
\bibitem[Ro1]{Rolfsen1} D. Rolfsen, {\em Knots and links}, 2nd printing with corrections, Mathematics Lecture Series {\bf 7}, Publish or Perish, Inc. (1990).
\bibitem[Ro2]{Rolfsen2} D. Rolfsen, {\em Rational surgery calculus: extension of Kirby's theorem}, Pacific J. Math. {\bf 110} (1984), no. 2, 377--386.
\bibitem[Roz]{Rozansky} L. Rozansky, {\em Residue formulas for the large $k$ asymptotics of Witten's invariants of Seifert manifolds. The case of $SU(2)$}, Comm. Math. Phys. {\bf 178} (1996), no. 1, 27--60.
\bibitem[Sa]{Sawin} S. F. Sawin, {\em Jones--Witten invariants for non-simply connected Lie groups and the geometry of the Weyl alcove}, preprint math.QA/9905010.
\bibitem[Sc]{Scott} P. G. Scott, {\em The geometries of $3$--manifolds}, Bull. London Math. Soc. {\bf 15} (1983), no. 5, 401--487.
\bibitem[Se1]{Seifert1} H. Seifert, {\em Topologie dreidimensionaler gefaserter r\"{a}ume}, Acta Math. {\bf 60} (1933), 147--238.
\bibitem[Se2]{Seifert2} H. Seifert, {\em Topology of $3$--dimensional fibered spaces} (english translation by W. Heil of: {\em Topologie dreidimensionaler gefaserter R\"{a}ume}, Acta Math. {\bf 60} (1933), 147--238). In H. Seifert, W. Threlfall, {\em A textbook of topology}, Pure Appl. Math. {\bf 89}, Academic Press (1980), 359--422.
\bibitem[Ta1]{Takata1} T. Takata, {\em On quantum $PSU(n)$--invariants for lens spaces}, J. Knot Theory Ramifications {\bf 5} (1996), no. 6, 885--901.
\bibitem[Ta2]{Takata2} T. Takata, {\em On quantum $PSU(n)$--invariants for Seifert manifolds}, J. Knot Theory Ramifications {\bf 6} (1997), no. 3, 417--426.
\bibitem[T]{Thurston} W. P. Thurston, {\em The geometry and topology of $3$--manifolds}, Lecture notes, Princeton University (1978).
\bibitem[Th1]{Thys1} H. Thys, {\em Groupes quantiques et cat\'{e}gories de diagrammes planaires}, D. Phil. thesis, Universit\'{e} Louis Pasteur, Strasbourg (2000).
\bibitem[Th2]{Thys2} H. Thys, {\em Description topologique des repr\'{e}sentations de $U_{q}(\frsl_{2})$}, Ann. Fac. Sci. Toulouse Math. {\bf 8} (1999), no. 4, 695--725.
\bibitem[Tr]{Trotter} H. F. Trotter, {\em Non-invertible knots exist}, Topology {\bf 2} (1963), 275--280.
\bibitem[Tu]{Turaev} V. G. Turaev, {\em Quantum invariants of knots and $3$--manifolds}, de Gruyter Stud. Math. {\bf 18}, Walter de Gruyter (1994).
\bibitem[TW1]{TuraevWenzl1} V. G. Turaev, H. Wenzl, {\em Quantum invariants of $3$--manifolds associated with classical simple Lie algebras}, Internat. J. Math. {\bf 4} (1993), no. 2, 323--358.
\bibitem[TW2]{TuraevWenzl2} V. G. Turaev, H. Wenzl, {\em Semisimple and modular categories from link invariants}, Math. Ann. {\bf 309} (1997), no. 3, 411--461.
\bibitem[V]{Verlinde} E. Verlinde, {\em Fusion rules and modular transformations in 2D conformal field theory}, Nuclear Phys. B {\bf 300} (1988), no. 3, 360--376.
\bibitem[Wa]{Wallace} A. H. Wallace, {\em Modifications and cobounding manifolds}, Canad. J. Math. {\bf 12} (1960), 503--528.
\bibitem[W]{Wenzl} H. Wenzl, {\em $C^{*}$ tensor categories from quantum groups}, J. Amer. Math. Soc. {\bf 11} (1998), no. 2, 261--282.
\bibitem[Wi]{Witten} E. Witten, {\em Quantum field theory and the Jones polynomial}, Comm. Math. Phys. {\bf 121} (1989), no. 3, 351--399.
\end{thebibliography}
\end{document}